\documentclass[a4paper,11pt]{article}
\usepackage[T1]{fontenc}
\usepackage[utf8]{inputenc}
\usepackage[francais]{babel}
\usepackage{amsmath}
\usepackage{amsthm}
\usepackage{amssymb}
\usepackage{amsfonts}
\usepackage{graphics}
\usepackage{color}
\newcounter{theorem}

\newtheorem{theo}[theorem]{Th\'eor\`eme}%
\newtheorem{prop}[theorem]{Proposition}%
\newtheorem{defi}[theorem]{D\'efinition}%
\newtheorem{lemm}[theorem]{Lemme}%
\newtheorem{coro}[theorem]{Corollaire}%
\newenvironment{demo}{\trivlist\item[\hskip \labelsep{\rm
D\'emonstration.}\enskip]}%

\def\card{{\rm Card}}

\def\sgn{{\rm sgn}}

\def\eef{{\bf E}}
\def\iif{{\bf I}}
\def\nnf{{\bf N}}
\def\ppf{{\bf P}}
\def\qqf{{\bf Q}}
\def\rrf{{\bf R}}

\def\ttf{{\bf T}}
\def\wwf{{\bf W}}
\def\zzf{{\bf Z}}

\def\ac{{\cal A}}

\def\ec{{\cal E}}
\def\fc{{\cal F}}

\def\hc{{\cal H}}
\def\ic{{\cal I}}

\def\lc{{\cal L}}

\def\wc{{\cal W}}

\def\osc{{\rm Osc}}
\def\amp{{\rm amp}}
\begin{document}

\title{Densit{\'e} des orbites des trajectoires browniennes sous l'action
  de la transformation de L{\'e}vy}
\author{Jean BROSSARD et Christophe LEURIDAN}
\maketitle

\begin{abstract}
  Etant donn{\'e} une transformation mesurable $T$ d'un espace probabilis{\'e}
  $(E,\ec,\pi)$ dans lui-m{\^e}me pr{\'e}servant la mesure, et une partie $B
  \in \ec$, nous donnons une condition suffisante pour que l'orbite de
  $\pi$-presque tout point visite $B$ : il suffit que $B$ soit
  accessible depuis $\pi$-presque tout point pour une cha{\^\i}ne de Markov
  de noyau $K$, o{\`u} $K(\cdot,\cdot)$ est une version r{\'e}guli{\`e}re de la
  loi conditionnelle de $X$ sachant $T(X)$ lorsque $X$ est une
  variable al{\'e}atoire de loi $\pi$. 

  Nous appliquons ensuite ce fait g{\'e}n{\'e}ral {\`a} la transformation 
  de L{\'e}vy,
  qui {\`a} un mouvement brownien $W$ associe le mouvement brownien
  $|W|-L$ o{\`u} $L$ est le temps local en $0$ de $W$. Cela nous permet de
  d\'emontrer le th\'eor\`eme de Malric : l'orbite de presque toute 
  trajectoire visite tout ouvert
  non vide de l'espace de Wiener $\wwf = {\cal C}(\rrf_+,\rrf)$ pour
  la topologie de la convergence uniforme sur les compacts. 
\end{abstract}

\begin{flushleft}
{\it Classification math\'ematique}~: 37A50,28D05,60J65.\\
{\it Mots-cl\'es}~: mouvement brownien, transformation de L{\'e}vy, 
densit{\'e} des orbites, r{\'e}currence, ergodicit{\'e}.
\end{flushleft}

\section*{Introduction}

Soit $W$ un mouvement brownien dans $\rrf$ issu de $0$. Le temps local
en $0$ de $W$, d{\'e}fini par la formule de Tanaka,  
$$|W_t| = \int_0^t \sgn(W_s)\ dW_s + L_t$$
est aussi la densit{\'e} en $0$ du temps d'occupation : 
pour tout $t \in \rrf_+$, 
$$L_t = \lim_{\epsilon \to 0} \frac{1}{2\epsilon} \int_0^t
\iif_{[|W_s| \le \epsilon]}\ ds\ \text{\ p.s.}$$
On v{\'e}rifie que la martingale locale $\hat{W} = \int_0^\cdot \sgn(W_s)\ dW_s$
est un mouvement brownien en remarquant que pour tout $t \in \rrf_+$, 
$$\langle \hat{W},\hat{W}\rangle_t = \int_0^t \sgn(W_s)^2\ ds = t\
\text{\ p.s.}$$
La transformation de L{\'e}vy est la transformation qui au mouvement brownien
$W$ associe le mouvement brownien $\hat{W}$. 

Plus pr{\'e}cis{\'e}ment, soit $\wwf$ l'ensemble des trajectoires continues de
$\rrf_+$ dans $\rrf$, nulles en $0$. La transformation de L{\'e}vy associe {\`a} presque
toute trajectoire $w \in \wwf$ la trajectoire $\ttf(w) = \hat{w}$ d{\'e}finie par  
$$\hat{w}(t) = |w(t)| - l(t,w),$$
avec 
$$l(t,w) = \liminf_{\epsilon \to 0} \frac{1}{2\epsilon} \int_0^t
\iif_{[|w(s)| \le \epsilon]}\ ds.$$
Lorsque l'application $l(\cdot,w)$ ainsi d{\'e}finie n'est pas continue,
on peut poser par exemple $\hat{w}=0$ pour que la transformation de
L{\'e}vy soit bien d{\'e}finie comme application de $\wwf$ dans $\wwf$, mais
cela est sans importance. 

La transformation de L{\'e}vy pr{\'e}serve la mesure de Wiener. Dans 
\cite{Revuz-Yor},
D. Revuz et M. Yor posent la question de savoir si elle est
ergodique. Cette question difficile n'est toujours pas r{\'e}solue, mais
des avanc{\'e}es ont {\'e}t{\'e} faites par M. Malric. 
Dans \cite{Malric 0}, Malric
d{\'e}montre que pour presque toute trajectoire $w \in \wwf$, l'ensemble
des z{\'e}ros de $w$ et de ses images successives par la transformation de
L{\'e}vy est dense dans $\rrf_+$. 

Par la suite, Malric a démontré un résultat nettement plus difficile : 
pour presque toute trajectoire $w \in \wwf$, l'orbite de $w$ sous
l'action de la transformation de L{\'e}vy est dense dans $\wwf$, pour la
topologie de la convergence uniforme sur les compacts. Ce r{\'e}sultat
s'approche de l'ergodicit{\'e} et serait une cons{\'e}quence de 
l'ergodicit{\'e} si elle {\'e}tait d{\'e}montr{\'e}e. Ce résultat fait 
l'objet de l'article \cite{Malric 5} paru à ESAIM PS. 

Entre 2007 et 2009, nous avons {\'e}crit une autre 
d{\'e}monstration du théorème de Malric, dans laquelle les preuves des 
r{\'e}sultats interm{\'e}diaires sont ind{\'e}pendantes les unes des autres, 
et qui fait l'objet du pr{\'e}sent article. 
Grâce au caractère très modulaire de notre d{\'e}monstration, nous avons mis en 
évidence un procédé d'approximation des trajectoires facile à comprendre, 
dont les étapes sont résumées sur la figure~\ref{etapes de l'approximation}. 
Même si le procédé est clair sur le dessin, nous nous sommes attachés 
à écrire en détail les démonstrations des résultats intermédiaires.

Le point de d{\'e}part est une id{\'e}e d{\'e}j{\`a} sous-jacente d{\`e}s 
les premi{\`e}res versions de l'article de 
Malric~\cite{Malric 2,Malric 3,Malric 4}. 
Il nous a paru int{\'e}ressant d'en d{\'e}gager la g{\'e}n{\'e}ralit{\'e} 
en {\'e}tablissant dans la premi{\`e}re partie un raffinement du 
th{\'e}or{\`e}me 
de r{\'e}currence de Poincar{\'e} et son corollaire suivant. 

\noindent{\bf Corollaire~\ref{Condition suffisante de recurrence}}
{\bf (Condition suffisante de r{\'e}currence)}
 
{\it Soient $T$ une transformation mesurable d'un espace probabilis{\'e}
  $(E,\ec,\pi)$ dans lui-m{\^e}me pr{\'e}servant la mesure. 

Supposons que pour toute variable al{\'e}atoire $X$ de loi $\pi$, la loi
conditionnelle de $X$ sachant $T(X)$ poss{\`e}de une version r{\'e}guli{\`e}re
$K(\cdot,\cdot)$.  

  Pour que l'orbite de $\pi$-presque tout point visite une
  infinit{\'e} de fois une partie $B \in \ec$, il suffit
  que $B$ soit accessible depuis $\pi$-presque tout point pour une
  cha{\^\i}ne de Markov de noyau $K$. }

Remarquons que $K$ est le noyau de transition de la cha{\^\i}ne de
Markov stationnaire index{\'e}e par $-\nnf$ obtenue en retournant la suite
de variables al{\'e}atoires $(T^n(X))_{n \in \nnf}$ pour une
variable al{\'e}atoire $X$ de loi $\pi$. 

Dans le cas o{\`u} $T$ est la transformation de L{\'e}vy, l'{\'e}volution 
pas {\`a}
pas de cette cha{\^\i}ne de Markov est simple {\`a} d{\'e}crire : partant d'une
trajectoire $w \in \wwf$, on lui retranche son minimum courant
$\underline{w}$ (d{\'e}fini par $\underline{w}(t) = \min\{w(s)\ ; s \in
[0,t]\}$). Puis on multiplie les excursions de la trajectoire positive
$w-\underline{w}$ par des signes pris au hasard. 

D'apr{\`e}s la condition suffisante de r{\'e}currence 
{\'e}nonc{\'e}e ci-dessus, il
s'agit de d{\'e}montrer que si $U$ est un ouvert non vide de $\wwf$, la
cha{\^\i}ne de Markov issue d'une trajectoire typique de $\wwf$ visite
$U$ avec probabilit{\'e} strictement positive. Concr{\`e}tement, on montre
qu'on atteint $U$ en un nombre fini d'{\'e}tapes en imposant la valeur
d'un nombre fini de signes {\`a} chaque {\'e}tape.   

Cette méthode a déjà été utilisée par Malric avec une formulation
différente, puisque Malric utilise des « remontées de Lévy »
tandis que nous raisonnons avec des chaînes de Markov. 
Nous expliquons le lien entre ces notions à la section~\ref{comparaison}. 
Néanmoins notre d{\'e}monstration 
s'{\'e}carte nettement de celle de Malric dans la mise en {\oe}uvre de 
cette m{\'e}thode : {\`a} l'aide de la propri{\'e}t{\'e} de Markov,
nous d{\'e}composons l'accessibilit{\'e} des ouverts en accessibilit{\'e}s
successives plus faciles {\`a} montrer. La preuve de ces accessibilit{\'e}s
interm{\'e}diaires utilise essentiellement trois id{\'e}es nouvelles : 
l'introduction d'une topologie sur $\wwf$ permettant le contr{\^o}le des
z{\'e}ros, l'utilisation de transformations sur les mouvements browniens
r{\'e}fl{\'e}chis, et l'exploitation de leur continuit{\'e} presque partout. 
Voyons cela un peu plus en d{\'e}tail. 

\paragraph{Introduction d'une nouvelle topologie sur $\wwf$.}
Tout d'abord, nous introduisons une topologie plus fine
que la topologie de la convergence uniforme sur les compacts pour
pouvoir contr{\^o}ler les z{\'e}ros des trajectoires et jouer librement avec
les signes. Cette topologie, que nous appelons {\it topologie de la convergence
uniforme sur les compacts avec contr{\^o}le des z{\'e}ros (CUCZ)} est associ{\'e}e
aux {\'e}carts $d^{CU}_t$ et $d^{CZ}_t$ d{\'e}finis ci-dessous. Pour $f,g \in
\wwf$, on note $Z(f)$, $Z(g)$ les ensembles des z{\'e}ros de $f$ et de
$g$ et $Z_t(f) = Z(f) \cap [0,t]$, $Z_t(g) = Z(g) \cap [0,t]$.
On pose alors 
$$d^{CU}_t(f,g) = ||g-f||_{[0,t]} = \max\{|g(s)-f(s)|\ ;\ s \in [0,t]\},$$
$$d_t^{CZ}(f,g) = \inf\{\delta > 0 : Z_{t-\delta}(f) \subset Z(g) +
]-\delta,\delta[\ \text{ et }\ Z_{t-\delta}(g) \subset Z(f) +
]-\delta,\delta[\}.$$
Autrement dit, l'in{\'e}galit{\'e} $d_t^{CZ}(f,g) < \delta$ signifie que tout
z{\'e}ro de $f$ ant{\'e}rieur {\`a} $t-\delta$ est {\`a} distance inf{\'e}rieure {\`a}
$\delta$ d'un z{\'e}ro de $g$, et inversement. 

Nous montrons que la topologie CUCZ poss{\`e}de une base d{\'e}nombrable 
d'ouverts, de la forme
$$V_t(f,\rho,\delta) = \{g \in V :
d^{CU}(f,g) < \rho\ ;\ d^{CZ}(f,g) < \delta\}.$$ 
avec $f \in \wwf$, $t>0$, $\rho>0$, $\delta>0$. Nous montrons que 
ces ouverts sont accessibles depuis presque toute trajectoire de $\wwf$, 
ce qui permet de d{\'e}montrer le r{\'e}sultat suivant. 

\begin{theo}~\label{Densite pour la topologie CUCZ} {\bf (Densit{\'e} 
des orbites pour la topologie CUCZ)}

L'orbite de presque toute trajectoire sous l'action de la
transformation de L{\'e}vy est dense pour la topologie CUCZ et {\it a
  fortiori} pour la topologie de la convergence uniforme sur les
compacts. 
\end{theo}

\paragraph{Transformations sur les mouvements browniens r{\'e}fl{\'e}chis.}
En fait, nous travaillons davantage sur les mouvements browniens
r{\'e}fl{\'e}chis que sur les mouvements browniens. 

Nous introduisons une action du groupe $E = \{-1,1\}^{\qqf_+^*}$ sur
$\wwf$ : $(e,w) \mapsto e \cdot w$. L'action d'une famille de signes
$e \in E$ sur une trajectoire $w \in \wwf$ consiste, une fois les
excursions de $w$ num{\'e}rot{\'e}es {\`a} l'aide des rationnels strictement
positifs, {\`a} multiplier par $e(q)$ l'excursion num{\'e}rot{\'e}e par le
rationnel $q$. 

Mais au lieu de travailler avec la \og transformation inverse de L{\'e}vy \fg
qui, {\`a} un mouvement brownien $W$ et une famille $\varepsilon$ de
signes choisie au hasard ind{\'e}pendamment de $W$ associe le mouvement
brownien $\varepsilon \cdot (W-\underline{W})$, on regarde plut{\^o}t la
transformation $F$ qui, {\`a} un mouvement brownien r{\'e}fl{\'e}chi $R$ et une
famille $\varepsilon$ de signes choisie au hasard ind{\'e}pendamment de
$R$ associe le mouvement brownien r{\'e}fl{\'e}chi $F(\varepsilon,R) =
\varepsilon \cdot R - \underline{\varepsilon \cdot R}$. 
Nous montrons que l'accessibilit{\'e} des ouverts de $\wwf$ par la
transformation inverse de L{\'e}vy se ram{\`e}ne {\`a} l'accessibilit{\'e} des ouverts
de $\wwf_+$ pour la topologie CUCZ par la transformation $F$.  

Nous utiliserons aussi la \og transformation $F$ apr{\`e}s un instant $a$\fg,
not{\'e}e $F_a$, qui pr{\'e}serve les trajectoires jusqu'{\`a} leur premier z{\'e}ro
apr{\`e}s $a$ et agit comme $F$ apr{\`e}s cet instant. La pr{\'e}servation par
$F_a$ du d{\'e}but des trajectoires (au moins sur l'intervalle $[0,a]$)
est extr{\^e}mement utile et simplifie bien des d{\'e}monstrations.

\paragraph{Utilisation de la continuit{\'e} presque partout.}
Nous montrons que les transformations $F_a$ sont continues presque
partout de $E \times \wwf_+$ dans $\wwf_+$, lorsque $E$ est muni de la
topologie produit et $\wwf_+$ de la topologie CUCZ.   
Cette propri{\'e}t{\'e} de continuit{\'e} est tr{\`e}s utile puisqu'elle
nous permet de remplacer certains arguments de nature probabiliste par
des arguments de nature topologique.  
Elle nous permet notamment de montrer que si $b \ge a \ge 0$,
l'accessibilit{\'e} d'un ouvert de $\wwf_+$ par $F_b$ entra{\^\i}ne son
accessibilit{\'e} par $F_a$ et son accessibilit{\'e} par $F$
(proposition~\ref{Comparaison des accessibilites}). 

\subsection*{Plan de l'article}

Dans la premi{\`e}re partie, nous d{\'e}montrons des r{\'e}sultats g{\'e}n{\'e}raux 
permettant de montrer que les orbites sous l'action d'une transformation 
pr{\'e}servant une probabilit{\'e} sur un espace poss{\'e}dant une base d{\'e}nombrable 
d'ouverts sont presque toutes denses. 

Dans la deuxi{\`e}me partie, nous appliquons ce r{\'e}sultat {\`a} la transformation 
de L{\'e}vy en introduisant la cha{\^\i}ne de Markov stationnaire sur $\wwf$ 
obtenue par inversion de la transformation de L{\'e}vy. Nous introduisons 
aussi une cha{\^\i}ne de Markov stationnaire sur $\wwf_+$ associ{\'e}e {\`a} la 
transformation $F$ d{\'e}crite plus haut et admettant comme mesure invariante 
la loi du mouvement brownien r{\'e}fl{\'e}chi. Nous verrons par la suite qu'il est 
plus commode de travailler avec cette seconde cha{\^\i}ne de Markov. 

Les trois parties suivantes sont consacr{\'e}es {\`a} l'affutage des outils 
permettant de montrer le th{\'e}or{\`e}me~\ref{Densite pour la topologie CUCZ} 
(densit{\'e} de presque toute orbite sous l'action de la transformation 
de L{\'e}vy) : 
\begin{itemize} 
\item dans la troisi{\`e}me partie, nous introduisons la topologie de la 
convergence uniforme avec contr{\^o}le des z{\'e}ros. Nous montrons en particulier 
que les fonctions positives, continues, affines par morceaux, nulles en $0$ 
sont denses dans $\wwf_+$ pour cette topologie. Pour montrer le 
th{\'e}or{\`e}me~\ref{Densite pour la topologie CUCZ}, 
nous sommes ramen{\'e}s (gr{\^a}ce au corollaire~\ref{Ce qu'il faudra montrer})
{\`a} montrer que tout ouvert non vide de 
$\wwf_+$ est accessible depuis presque toute trajectoire de $\wwf_+$ par 
la cha{\^\i}ne de Markov que nous avons construite sur $\wwf_+$. 
\item dans la quatri{\`e}me partie, nous {\'e}tablissons des r{\'e}sultats de continuit{\'e} 
presque partout pour la topologie de la convergence uniforme avec contr{\^o}le 
des z{\'e}ros, notamment la continuit{\'e} des transformations \og $F$ apr{\`e}s $a$\fg. 
\item dans la cinqui{\`e}me partie, nous introduisons la notion d'accessibilit{\'e} 
d'un ouvert ou d'un \og presque-ouvert\fg par les transformations 
\og $F$ apr{\`e}s $a$\fg. Nous utilisons ces r{\'e}sultats de continuit{\'e} presque 
partout pour {\'e}tablir des propri{\'e}t{\'e}s de cette notion qui nous seront 
tr{\`e}s utiles, notamment le r{\'e}sultat de comparaison des accessibilit{\'e}s 
(proposition~\ref{Accessibilites successives et intersection}).
\end{itemize}

\medskip
Dans les trois derni{\`e}res parties, nous {\'e}tudions l'accessibilit{\'e} de 
certains ouverts depuis presque toute trajectoire d'autres ouverts. 
Mises bout-{\`a}-bout, ces accessibilit{\'e}s constituent avec le 
corollaire~\ref{Ce qu'il faudra montrer} la preuve du 
th{\'e}or{\`e}me~\ref{Densite pour la topologie CUCZ} : 
\begin{itemize}
\item dans la sixi{\`e}me partie, nous montrons que tout voisinage de la 
trajectoire nulle est accessible par $F$ depuis presque toute 
trajectoire de $\wwf_+$. L'utilisation de la transformation $F_a$ 
au lieu de la transformation $F$ permet donc de \og remettre {\`a} z{\'e}ro\fg 
une trajectoire apr{\`e}s son premier z{\'e}ro apr{\`e}s $a$ en pr{\'e}servant la 
trajectoire avant cet instant.
\item dans la septi{\`e}me partie, nous montrons comment construire des excursions 
de hauteur voulue dans un intervalle de temps donn{\'e}, et nous {\'e}tablissons 
le lemme du v{\'e}rin, qui permet de soulever une trajectoire d'une hauteur 
donn{\'e}e sur un intervalle donn{\'e}. 
\item dans la huiti{\`e}me et derni{\`e}re partie, nous utilisons le lemme de remise 
{\`a} z{\'e}ro et le lemme du v{\'e}rin pour montrer comment approcher n'importe 
quelle fonction positive, 
continue, affine par morceaux, nulle en $0$, ce qui compl{\`e}te la preuve.
\end{itemize}
\section{R{\'e}sultats g{\'e}n{\'e}raux}

\subsection{R{\'e}currence et accessibilit{\'e}}

Soit $T$ une transformation d'un espace probabilis{\'e} $(E,\ec,\pi)$ dans
lui-m{\^e}me pr{\'e}servant la probabilit{\'e} $\pi$. 

Supposons que pour toute variable al{\'e}atoire $X$ de loi $\pi$, la loi
conditionnelle de $X$ sachant $T(X)$ poss{\`e}de une version r{\'e}guli{\`e}re
$K(\cdot,\cdot)$.  

Pour $B \in \ec$, notons :
\begin{itemize} 
\item $U(B)$ l'ensemble des $x \in E$ tels que l'orbite de $x$ sous
  l'action de $T$ visite $B$,  
$$U(B) = \bigcup_{k \in \nnf} T^{-k}(B)\ ;$$
\item $R(B)$ l'ensemble des $x \in E$ tels que l'orbite de $x$ sous
  l'action de $T$ visite $B$ une infinit{\'e} de fois,  
$$R(B) = \limsup_{k \to +\infty} T^{-k}(B)\ ;$$
\item $A(B)$ l'ensemble des {\'e}tats depuis lesquels $B$ est accessible pour une
cha{\^\i}ne de Markov de noyau $K(\cdot,\cdot)$, 
$$A(B) = \{x \in E : \exists n \in \nnf : K^n(x,B)>0\}.$$
\end{itemize}

Le th{\'e}or{\`e}me de r{\'e}currence de Poincar{\'e} (voir \cite{Petersen}) 
affirme que $B \subset R(B)$ $\pi$-presque s{\^u}rement. 
Comme $A(B)$ contient $B$, la proposition 
ci-dessous renforce ce r{\'e}sultat. 

\begin{prop} {\bf (Lien entre r{\'e}currence et accessibilit{\'e})}

Sous les hypoth{\`e}ses ci-dessus, la probabilit{\'e} $\pi$ est invariante
pour le noyau $K$ et $A(B) \subset U(B) = R(B)$
$\pi$-presque s{\^u}rement.   
\end{prop}

\begin{demo} Soit $(X_n)_{n \in \nnf}$ une cha{\^\i}ne de Markov de loi
initiale $\pi$ et de noyau de transition $K(\cdot,\cdot)$. 

1. Montrons d'abord que pour tout $n \in \nnf^*$, $X_n$ a pour loi $\pi$ et
$T^n(X_n) = X_0$ presque s{\^u}rement. 

Soit $Y$ une variable al{\'e}atoire de loi $\pi$. Alors
$T^n(Y),T^{n-1}(Y),\ldots,T(Y),Y$ ont pour loi $\pi$ et pour tout $k \in
[0,\ldots,n-1]$, on a l'{\'e}galit{\'e} des lois conditionnelles
$$\lc \big(T^{k}(Y) \big| T^n(Y),\ldots,T^{k+1}(Y) \big) =
\lc \big(T^k(Y) \big| T^{k+1}(Y) \big) = K \big(T^{k+1}(Y),\cdot \big).$$
Par cons{\'e}quent, $(T^n(Y),T^{n-1}(Y),\ldots,Y)$ a m{\^e}me loi
que $(X_0,X_1\ldots,X_n)$. 
En particulier, $(X_0,X_n)$ a m{\^e}me loi que $(T^n(Y),Y)$, ce qui entra{\^\i}ne
le r{\'e}sultat annonc{\'e}. 

2. L'{\'e}galit{\'e} presque s{\^u}re $U(B) = R(B)$ vient du fait que $R(B)$ est
l'intersection des parties  
$$R_n(B) = \bigcup_{k \ge n} T^{-k}(B) = \bigcup_{l \in \nnf}
T^{-(n+l)}(B) = T^{-n}(U(B))$$ 
qui sont embo{\^\i}t{\'e}es et de m{\^e}me probabilit{\'e}.

3. Montrons l'inclusion presque s{\^u}re $A(B) \subset U(B)$. Comme
$T^n(X_n) = X_0$ presque s{\^u}rement, 
$$\{X_0 \in U(B)\} \subset \{X_n \in U(B)\}\ \text{p.s.}.$$ 
Mais ces {\'e}v{\'e}nements ont m{\^e}me probabilit{\'e}, puisque $X_n$ a m{\^e}me loi que
$X_0$. Ils sont donc {\'e}gaux presque s{\^u}rement. 
En particulier, 
$$\{X_0 \in U(B)\} \supset \{X_n \in B\}\ \text{p.s.},$$
d'o{\`u} en conditionnant par rapport {\`a} $X_0$ : 
$$\iif_{U(B)}(X_0) \ge P[X_n \in B|\sigma(X_0)] = K^n(X_0,B)\
\text{p.s.}$$ 
Ainsi, 
$$\iif_{U(B)}(X_0) \ge \sup_{n \in \nnf} K^n(X_0,B)\ \text{p.s.},$$ 
ce qui montre le r{\'e}sultat annonc{\'e}. \hfill $\square$
\end{demo}

De la proposition pr{\'e}c{\'e}dente, on d{\'e}duit imm{\'e}diatement une condition
suffisante de r{\'e}currence pour une partie fix{\'e}e. 

\begin{coro}\label{Condition suffisante de recurrence}
{\bf (Condition suffisante de r{\'e}currence)}
 
  Soient $T$ une transformation mesurable d'un espace probabilis{\'e}
  $(E,\ec,\pi)$ dans lui-m{\^e}me pr{\'e}servant la mesure. 
 
Supposons que pour toute variable al{\'e}atoire $X$ de loi $\pi$, la loi
conditionnelle de $X$ sachant $T(X)$ poss{\`e}de une version r{\'e}guli{\`e}re
$K(\cdot,\cdot)$.  

  Pour que l'orbite de $\pi$-presque tout point visite une
  infinit{\'e} de fois une partie $B \in \ec$, il suffit
  que $B$ soit accessible depuis $\pi$-presque tout point par une
  cha{\^\i}ne de Markov de noyau $K$. 
\end{coro}

{\parindent 0cm {\bf Remarques}} 

\begin{itemize}
\item Comme $\pi$ est invariante pour le noyau $K$, on peut construire une 
cha{\^\i}ne de Markov stationnaire index{\'e}e par $\zzf$ de noyau $K$ en 
munissant l'espace canonique $S^\zzf$ d'une probabilit{\'e} ad hoc. 
L'op{\'e}rateur de d{\'e}calage 
$(x_n)_{n \in \zzf} \mapsto (x_{n-1})_{n \in \zzf}$ 
sur cet espace est l'extension naturelle de $(E,\ec,\pi,T)$ : 
voir par exemple \cite{Petersen}.

\item Les \og mouvements browniens remont{\'e}s\fg introduits par Malric 
d{\`e}s la premi{\`e}re version~\cite{Malric 2} sont proches de cette 
cha{\^\i}ne de Markov. Le corollaire 3 ci-dessus joue dans notre article 
le m{\^e}me r{\^o}le que la 
proposition 1 dans celui de Malric \cite{Malric 2}.
\end{itemize}

Lorsque $E$ est un espace topologique poss{\'e}dant une base d{\'e}nombrable
d'ouverts, on obtient une condition suffisante
de densit{\'e} des orbites. 

\begin{coro}\label{Condition suffisante de densite des orbites}
{\bf (Condition suffisante de densit{\'e} des orbites)}
 
Soient $E$ un espace topologique poss{\'e}dant une base
d{\'e}nombrable d'ouverts, $\ec$ sa tribu bor{\'e}lienne, $\pi$ une
probabilit{\'e} sur $(E,\ec)$ et $T$ une transformation mesurable de
$(E,\ec,\pi)$ dans lui-m{\^e}me pr{\'e}servant la mesure. 

Supposons que pour toute variable al{\'e}atoire $X$ de loi $\pi$, la loi
conditionnelle de $X$ sachant $T(X)$ poss{\`e}de une version r{\'e}guli{\`e}re
$K(\cdot,\cdot)$.  

  Pour que l'orbite de $\pi$-presque tout point soit dense dans $E$, il suffit
  que tout ouvert de $E$ soit accessible depuis $\pi$-presque tout
  point par une cha{\^\i}ne de Markov de noyau $K$. 
\end{coro}

\subsection{Utilisation de la propri{\'e}t{\'e} de Markov}

La propri{\'e}t{\'e} de Markov permet d'{\'e}tablir l'accessibilit{\'e} d'une
partie en d{\'e}composant la d{\'e}monstration en plusieurs {\'e}tapes.  

\begin{lemm}~\label{Accessibilites successives} 
{\bf (Accessibilit{\'e}s successives)}
 
  Soit $K(\cdot,\cdot)$ un noyau de transition sur $(E,\ec)$
  ayant comme probabilit{\'e} invariante $\pi$. 

  Soient $B_0,B_1,B_2 \in \ec$ de mesure strictement positive pour
  $\pi$. 

  Si $B_2$ est accessible depuis $\pi$-presque tout point de $B_1$ et 
$B_1$ est accessible depuis $\pi$-presque tout point de $B_0$, alors
$B_2$ est accessible depuis $\pi$-presque tout point de $B_0$. 
\end{lemm}

\begin{demo}
Pour tout probabilit{\'e} $\mu$ sur $E$, notons $\ppf^\mu$ la probabilit{\'e}
sur l'espace canonique $E^\nnf$
faisant du processus canonique $(X_n)_{n \in \nnf}$ une cha{\^\i}ne de
Markov de noyau de transition $K(\cdot,\cdot)$ et de loi initiale
$\mu$. Notons $\tau_1$ et $\tau_2$ les temps d'atteinte de $B_1$ et
$B_2$. Pour $i,j \in \{0,1,2\}$, notons 
$A_{i,j} = \{x \in B_i : \ppf^x[\tau_j < +\infty]=0\}$,
Par hypoth{\`e}se, $\pi(A_{0,1}) = 0$
et $\pi(A_{1,2}) = 0$. Il s'agit de montrer que $\pi(A_{0,2})=0$. 

On commence par v{\'e}rifier que la loi de $X_{\tau_1}$ sous
$\ppf^\pi[\cdot|\tau_1 < +\infty]$ est absolument continue par
rapport {\`a} $\pi$. En effet, si $A \in \ec$ v{\'e}rifie $\pi(A)=0$, alors  
\begin{eqnarray*}
\ppf^\pi [X_{\tau_1} \in A\ ;\ \tau_1 < +\infty] 
\le \sum_{n \in \nnf} \ppf^\pi [X_n \in A] 
= \sum_{n \in \nnf} \pi(A) = 0,
\end{eqnarray*}
ce qui montre l'absolue continuit{\'e}. En particulier, $\ppf^{X_{\tau_1}}
[\tau_2 < +\infty]>0$ $\ppf^\pi$-presque s{\^u}rement sur l'{\'e}v{\'e}nement
$\{\tau_1 < +\infty\}$ puisque $X_{\tau_1}$ est {\`a} valeurs dans $B_1$. 

Par ailleurs, par en utilisant la d{\'e}finition de $A_{0,2}$ et la
propri{\'e}t{\'e} de Markov au temps $\tau_1$, on obtient
\begin{eqnarray*} 
0 
&=& \int_{A_{0,2}} \ppf^x[\tau_1 \le \tau_2 < +\infty]\ \pi(dx)\\
&=& \ppf^\pi [X_0 \in A_{0,2}\ ; \tau_1 \le \tau_2 < +\infty] \\
&=& \eef^{\pi} \Big[\iif_{[X_0 \in A_{0,2}]}\ \iif_{[\tau_1 < +\infty]}\
\ppf^{X_{\tau_1}} [\tau_2 < +\infty] \Big].
\end{eqnarray*}
Comme $\ppf^{X_{\tau_1}} [\tau_2 < +\infty]>0$ $\ppf^\pi$-p.s.
sur l'{\'e}v{\'e}nement $\{\tau_1 < +\infty\}$, on a ainsi 
$$\ppf^\pi [X_0 \in A_{0,2}\ ;\ \tau_1 < +\infty] = 0.$$
Autrement dit, $\pi(A_{0,2} \setminus A_{0,1})=0$, d'o{\`u} $\pi(A_{0,2})=0$. 
\hfill $\square$
\end{demo}

\section{Inversion de la transformation de L{\'e}vy}

\subsection{Rappels sur la transformation de L{\'e}vy}

Soit $\wwf$ l'ensemble des trajectoires continues de $\rrf_+$ dans
$\rrf$, nulles en $0$. On munit $\wwf$ de la tribu $\wc$ engendr{\'e}e par les
applications coordonn{\'e}es et de la mesure de Wiener. La tribu $\wc$ est
aussi la tribu bor{\'e}lienne pour la topologie de la convergence uniforme
sur les compacts. Pour tout $w \in \wwf$, on note $\underline{w}$ 
l'application de $\wwf$ d{\'e}finie par 
$$\underline{w}(t) = \min\{w(s)\ ;\ s \in [0,t]\}.$$

La transformation de L{\'e}vy est une application mesurable de $\wwf$ dans
$\wwf$, pr{\'e}servant la mesure de Wiener : si $W$ est un mouvement
brownien dans $\rrf$ issu de $0$, son image par la 
transformation de L{\'e}vy est le mouvement brownien 
$$\hat{W} = \int_0^\cdot \sgn(W_s) dW_s = |W| - L$$
o{\`u} $L$ est le temps local en $0$ du mouvement brownien $W$, qu'on peut
d{\'e}finir comme densit{\'e} de temps d'occupation : 
$$L_t = \liminf_{\epsilon \to 0} \frac{1}{2\epsilon} \int_0^t
\iif_{[|W_s| \le \epsilon]}\ ds.$$

Pour appliquer les r{\'e}sultats de la premi{\`e}re partie, il nous faut
d{\'e}terminer la loi de $W$ sachant $\hat{W}$. L'{\'e}galit{\'e} de processus
$\hat{W} + L = |W|$ et le fait que le temps local $L$ ne croisse que
sur l'ensemble des z{\'e}ros de $W$ entra{\^\i}ne l'{\'e}galit{\'e} $L_t =
\max\{-\hat{W}_s\ ; 0 \le s \le t\}$ pour tout $t$, d'o{\`u} 
$$|W_t| = \hat{W}_t - \underline{\hat{W}}_t \text{ avec } 
\underline{\hat{W}}_t = \min\{\hat{W}_s\ ; 0 \le s \le t\}.$$ 
Cette {\'e}galit{\'e} montre que le mouvement brownien r{\'e}fl{\'e}chi $|W|$ est une
fonction mesurable de $\hat{W}$. Mais la d{\'e}finition de $\hat{W}$
montre que $\hat{W}$ est une fonction mesurable de $|W|$. Ces
processus engendrent donc la m{\^e}me tribu. 

La transformation de L{\'e}vy perd donc de l'information, plus
pr{\'e}cis{\'e}ment les signes des excursions de $W$. Nous allons montrer que
ces signes sont ind{\'e}pendants et de loi uniforme sur $\{-1,1\}$
conditionnellement {\`a} $W$. Mais pour donner un sens pr{\'e}cis {\`a} cette
affirmation, nous devons num{\'e}roter les excursions.

\subsection{Num{\'e}rotation des excursions d'une trajectoire de $\wwf$}

Dans toute la suite, nous munirons l'ensemble $\qqf_+^*$ des rationnels
strictement positifs d'un ordre tel que 
\begin{itemize} 
\item toute partie non vide de $\qqf_+^*$ poss{\`e}de un premier
{\'e}l{\'e}ment ;
\item avant tout {\'e}l{\'e}ment de $\qqf_+^*$, il n'y a qu'un nombre fini
d'{\'e}l{\'e}ments. 
\end{itemize}
Par exemple, on ordonne les rationnels suivant la somme du num{\'e}rateur
et du d{\'e}nominateur puis, pour une somme fix{\'e}e, par num{\'e}rateurs croissants : 
$$\qqf_+^* = \Big\{\frac{1}{1} ; \frac{1}{2} ; \frac{2}{1} ; \frac{1}{3} ;
\frac{3}{1} ; \frac{1}{4} ; \frac{2}{3} ; \frac{3}{2} ; \frac{4}{1} ;
\frac{1}{5} ; \frac{5}{1} ; \ldots \Big\}.$$
Si $A$ est une partie non vide de $\qqf_+^*$, nous noterons $q(A)$ son
premier {\'e}l{\'e}ment.

Nous pouvons ainsi num{\'e}roter les excursions des trajectoires $w \in
\wwf$ par des rationnels, en posant $Q_t(w) = q(\qqf_+^*
\cap I_t(w))$ si $w(t) \ne 0$, o{\`u} $I_t(w)$ est l'intervalle ouvert
d'excursion enjambant $t$. On pose $Q_t(w) = 0$ si $w(t) =
0$. Remarquons que tous les rationnels ne servent pas dans la num{\'e}rotation. 

\subsection{Action d'une famille de signes sur un trajectoire de $\wwf$}

On d{\'e}finit une action du groupe $E = \{-1,1\}^{\qqf_+^*}$ sur $\wwf$
de la fa{\c c}on suivante : l'action d'une famille de signes $e \in E$ sur
une trajectoire $w$ est de multiplier chaque excursion de $w$ par
$e(q)$ o{\`u} $q$ est le \og num{\'e}ro de l'excursion\fg. Plus pr{\'e}cis{\'e}ment, pour
tout $t \in \rrf_+$, 
$$(e \cdot w)(t) = e(Q_t(w)) w(t),$$
avec la convention $e(0)=0$. On remarque que la trajectoire $e \cdot w$
ne d{\'e}pend que de $w$ et des signes $e(q)$ pour les rationnels $q$
num{\'e}rotant les excursions de $w$. 

Cette action du groupe va nous servir {\`a} construire un mouvement
brownien de valeur absolue $R$ donn{\'e}e {\`a} l'aide d'une famille de signes
de loi uniforme sur $E$, ind{\'e}pendante de $R$. Ce r{\'e}sultat fait l'objet
de la proposition~\ref{action sur le mouvement brownien reflechi}
ci-dessous. Pour le d{\'e}montrer, commen{\c c}ons par {\'e}tablir un lemme simple.  

\begin{lemm} {\bf (Changement de signe d'une excursion)}

Soit $W$ un mouvement brownien dans $\rrf$ issu de $0$. Soit $q>0$
fix{\'e}. Le processus $W'$ obtenu {\`a} partir de $W$ en changeant le signe de
l'excursion enjambant $q$ est encore un mouvement brownien.  
\end{lemm}

\begin{demo}
Nous donnons une d{\'e}monstration {\'e}l{\'e}mentaire qui ne repose pas sur la th{\'e}orie
des excursions. 

Par invariance d'{\'e}chelle du mouvement brownien, on se
ram{\`e}ne au cas o{\`u} $q=1$. Notons $]g_1,d_1[$ l'intervalle d'excursion
enjambant $1$. 

Comme $d_1$ est un z{\'e}ro du mouvement brownien $W$ et un
temps d'arr{\^e}t pour sa filtration naturelle $\fc^W$, le processus
$W_{d_1+\cdot}$ est un mouvement brownien ind{\'e}pendant de
$\fc^W_{d_1}$. Par cons{\'e}quent, le processus obtenu {\`a} partir de $W$ en
changeant le signe apr{\`e}s l'instant $d_1$ est encore un mouvement brownien.

En appliquant ce r{\'e}sultat au mouvement brownien $(tW_{1/t})_{t \ge
  0}$, on voit {\'e}galement que le processus obtenu {\`a} partir de $W$ en
changeant le signe avant l'instant $g_1$ est encore un mouvement
brownien. Par composition, le processus obtenu en changeant les signes
avant $g_1$ et apr{\`e}s $d_1$ est encore un mouvement brownien, et son
oppos{\'e} aussi. On en d{\'e}duit que le processus obtenu {\`a} partir de
$W$ en changeant le signe de l'excursion enjambant $1$ est encore un
mouvement brownien. \hfill $\square$  
\end{demo}

\begin{prop}~\label{action sur le mouvement brownien reflechi} 
{\bf (Action d'une famille de signes al{\'e}atoires sur un mouvement brownien 
r{\'e}fl{\'e}chi ind{\'e}pendant)}

Soit $R$ un mouvement brownien r{\'e}fl{\'e}chi. Soit $\varepsilon$ une variable
al{\'e}atoire {\`a} valeurs dans $E$, ind{\'e}pendante de $R$ et de loi uniforme
sur $E$. Alors $\varepsilon \cdot R$ est un mouvement brownien. 
\end{prop}

\begin{demo}
Soient $W$ un mouvement brownien et $\eta$ une variable al{\'e}atoire 
ind{\'e}pendante de $W$ et de loi uniforme sur $E = \{-1,1\}^{\qqf_+^*}$. 
Pour $q \in \qqf_+^*$, notons $A_q$ l'{\'e}v{\'e}nement \og $q$ num{\'e}rote 
une excursion de $|W|$\fg (autrement dit, $q$ est le premier rationnel de 
$I_q(R)$, l'intervalle d'excursion de $R$ enjambant $q$) et posons 
$$\xi(q) = {\rm sgn}(W_q) \iif_{A_q} + \eta(q) \iif_{A_q^c}.$$
La variable al{\'e}atoire $\xi$ ainsi d{\'e}finie peut s'{\'e}crire sous 
la forme $g(W,\eta)$ o{\`u} $g$ est une fonction mesurable
de $\wwf \times E$ dans $E$. 

Comme $|W|$ est un mouvement brownien r{\'e}fl{\'e}chi et comme 
$\xi \cdot |W| = W$, il suffit de
montrer que la variable al{\'e}atoire $\xi$ est
ind{\'e}pendante de $|W|$ et de loi uniforme sur $E$. 
Il s'agit donc de montrer que pour tout $q \in \qqf_+^*$, la loi
conditionnelle de $\xi$ sachant $|W|$ est invariante par $s_q$, o{\`u}
$s_q : E \to E$ est la sym{\'e}trie d{\'e}finie par $s_q(e)(q) =
-e(q)$ et $s_q(e)(q') = e(q')$ pour $q' \ne q$. Soit $W'$ le
mouvement brownien obtenu {\`a} partir de $W$ en changeant le signe de
l'excursion enjambant $q$. Le r{\'e}sultat vient de l'{\'e}galit{\'e} 
$$s_q(\xi) =  g(W',\eta) \iif_{A_q} + g(W,s_q(\eta)) \iif_{A_q^c},$$
du fait que $W'$ a m{\^e}me valeur absolue que $W$ et du fait que
l'{\'e}v{\'e}nement $A_q$ est mesurable par rapport {\`a} $\sigma(|W|)$. 
\hfill $\square$ 
\end{demo}

\begin{coro} {\bf (Loi de $W$ connaissant $\hat{W}$)}

Soit $W$ un mouvement brownien. Soit $\varepsilon$ une variable
al{\'e}atoire {\`a} valeurs dans $E$, ind{\'e}pendante de $\hat{W}$ et de loi uniforme
sur $E$. Alors $(W,\hat{W})$ a m{\^e}me loi que
$(\varepsilon \cdot (\hat{W}-\underline{\hat{W}}),\hat{W})$. Par
cons{\'e}quent, une version r{\'e}guli{\`e}re de la loi de $W$ sachant $\hat{W}$
est $(K(w,\cdot))_{w \in \wwf}$, o{\`u} $K(w,\cdot)$ est la loi de
$\varepsilon \cdot (w-\underline{w})$. 
\end{coro}

\begin{demo}
En appliquant la propri{\'e}t{\'e} pr{\'e}c{\'e}dente {\`a} $R = |W| =
\hat{W}-\underline{\hat{W}}$, on voit que $\varepsilon \cdot |W|$ et
un mouvement brownien. Donc $(W,|W|)$ a m{\^e}me loi que $(\varepsilon
\cdot |W|,|W|)$. Comme $\hat{W} = |W|-L$ est une
fonction de $|W|$, on en d{\'e}duit que $(W,\hat{W})$ a m{\^e}me loi que
$(\varepsilon \cdot (\hat{W}-\underline{\hat{W}}),\hat{W})$. 
\hfill $\square$
\end{demo}

\subsection{Construction de deux cha{\^\i}nes de Markov, sur  $\wwf$ et sur
  $\wwf_+$}~\label{construction de chaines de Markov}

Le corollaire pr{\'e}c{\'e}dent fournit un moyen simple de construire une
cha{\^\i}ne de Markov sur $\wwf$ stationnaire pour le noyau $K$ : on se
donne un mouvement brownien $W^{(0)}$ et une suite de variables al{\'e}atoires
$\varepsilon_1,\varepsilon_2,\ldots$ ind{\'e}pendantes et de loi uniforme sur
$E$, ind{\'e}pendante de $W$. On pose 
\begin{eqnarray*}
R^{(0)} &=& W^{(0)}-\underline{W^{(0)}},\\ 
W^{(1)} &=& \varepsilon_1 \cdot R^{(0)},\\ 
R^{(1)} &=& W^{(1)}-\underline{W^{(1)}},\\ 
W^{(2)} &=& \varepsilon_2 \cdot R^{(1)},
\end{eqnarray*}
et ainsi de suite. 

Notons $\wwf_+$ la partie de $\wwf$ form{\'e}e des trajectoires
positives et posons 
$$F(e,r) = e \cdot r - \underline{e \cdot r}\ \text{ pour }\
e \in E, r \in \wwf_+.$$ 
La relation de r{\'e}currence $R^{(n)} = F(\varepsilon_n,R^{(n-1)})$
montre que la suite $(R^{(n)})_{n \in \nnf}$ est une cha{\^\i}ne de Markov
sur le sous-ensemble $\wwf_+$. Nous allons voir qu'il est plus commode de
travailler avec cette cha{\^\i}ne plut{\^o}t qu'avec la cha{\^\i}ne $(W^{(n)})_{n \in \nnf}$.  

\subsection{N{\'e}cessit{\'e} de contr{\^o}ler les z{\'e}ros}

Pour $f \in \wwf$, $t>0$ et $h>0$, notons $V_t(f,\rho)$ la boule de
centre $f$ et de rayon $\rho$ pour la norme de la convergence uniforme
sur $[0,t]$ : 
$$V_t(f,\rho) = \{g \in V : ||g-f||_{[0,t]} < \rho\}.$$ 
Les boules $V_t(f,\rho)$ pour $f$ polyn{\^o}me {\`a} coefficients
rationnels, $t \in \nnf$ et $\rho \in \qqf_+^*$ forment une base
d{\'e}nombrable d'ouverts pour la topologie de la convergence 
uniforme sur les compacts. 

Pour montrer la densit{\'e} presque s{\^u}re des orbites sous l'action de la
transformation de L{\'e}vy, il suffit d'apr{\`e}s le corollaire~\ref{Condition
  suffisante de densite des orbites} de montrer que pour tout $f \in \wwf$,
$t>0$ et $h>0$, il existe $n \in \nnf$ tel que $K^n(W^{(0)},V_t(f,h))
> 0$ presque s{\^u}rement, autrement dit que, 
$$\sup_{n \in \nnf} P[W^{(n)} \in V_t(f,\rho)|W^{(0)}] > 0\
\text{ p.s.}$$
Mais pour $n \in \nnf^*$, $R^{(n-1)}=|W^{(n)}|$ donc 
$$\big|\big|\ R^{(n)} - |f|\ \big|\big|_{[0,t]} \le \big|\big|\ W^{(n)}
- f\ \big|\big|_{[0,t]}.$$ 
Pour que $W^{(n)}$ soit proche de $f$ pour la norme de la convergence uniforme
sur $[0,t]$, il est donc n{\'e}cessaire que $R^{(n)}$ soit proche de $|f|$
pour la norme de la convergence uniforme sur $[0,t]$. Mais le fait que
$R^{(n)}$ soit proche de $|f|$ ne garantit pas qu'on puisse rendre
$W^{(n)}$ proche de $f$ par un choix convenable de la famille de
signes $\varepsilon_n$.  
En effet, la trajectoire $W^{(n)} = \varepsilon_n \cdot R^{(n)}$ ne peut
changer de signe qu'en un z{\'e}ro de $R^{(n)}$. Si $f$ poss{\`e}de un z{\'e}ro
isol{\'e} $z_0$ et change de signe en $z_0$, il faut donc que $R^{(n)}$
poss{\`e}de un z{\'e}ro proche de $z_0$ pour que $W^{(n)}$ puisse approcher $f$.    

C'est pourquoi nous allons introduire une topologie prenant en compte
la distance entre les z{\'e}ros. 

\section{Topologie de la convergence
uniforme sur les compacts avec contr{\^o}le des z{\'e}ros}

\subsection{Construction d'{\'e}carts d{\'e}finissant la topologie}

Pour $f \in \wwf$ et $t>0$ on note $Z(f) = \{s \in \rrf_+ : f(s)=0\}$
et $Z_t(f) = Z(f) \cap [0,t]$. 

Un {\'e}cart naturel sur $\wwf$ permettant de s'assurer que les z{\'e}ros de
$f$ sur $[0,t]$ sont proches de ceux de $g$ et inversement, est la distance de
Hausdorff entre $Z_t(f)$ et $Z_t(g)$ : 
\begin{eqnarray*}
D_t(f,g) &=& d_H(Z_t(f),Z_t(g))\\ 
&=& \max\{\delta > 0 : Z_t(f) \subset Z_t(g) +
[-\delta,\delta]\ \text{ et }\ Z_t(g) \subset Z_t(f) + [-\delta,\delta]\}
\end{eqnarray*}
Cet {\'e}cart pr{\'e}sente cependant l'inconv{\'e}nient d'{\^e}tre sensible aux
ph{\'e}nom{\`e}nes de bord : pour $t > s > 0$, $D_t(f,g) > t-s$ lorsque
$f(t)=0$ et $g$ ne s'annule pas sur $[s,t]$, m{\^e}me si $g$ s'annule peu
apr{\`e}s $t$. On rem{\'e}die {\`a} cet inconv{\'e}nient en posant 
$$d_t^{CZ}(f,g) = \inf\{\delta > 0 : Z_{t-\delta}(f) \subset Z(g) +
]-\delta,\delta[\ \text{ et }\ Z_{t-\delta}(g) \subset Z(f) +
]-\delta,\delta[\}.$$
On remarque que l'inclusion $Z_{t-\delta}(f) \subset Z(g) +
]-\delta,\delta[$ est d'autant plus facile {\`a} r{\'e}aliser que $\delta$ est
grand, et que si elle est v{\'e}rifi{\'e}e pour un r{\'e}el $\delta>0$, on peut
trouver $\delta'<\delta$ pour lequel elle est encore v{\'e}rifi{\'e}e. 
On a donc l'{\'e}quivalence
$$d_t^{CZ}(f,g) < \delta \Longleftrightarrow Z_{t-\delta}(f) \subset Z(g) +
]-\delta,\delta[\ \text{ et }\ Z_{t-\delta}(g) \subset Z(f) +
]-\delta,\delta[.$$
Autrement dit, l'in{\'e}galit{\'e} $d_t^{CZ}(f,g) < \delta$ signifie que tout
z{\'e}ro de $f$ ant{\'e}rieur {\`a} $t-\delta$ est {\`a} distance inf{\'e}rieure {\`a}
$\delta$ d'un z{\'e}ro de $g$, et inversement. {\`A} l'aide de l'{\'e}quivalence
ci-dessus, on v{\'e}rifie facilement que la formule d{\'e}finit un {\'e}cart sur $\wwf$. 

Notons $d_t^{CU}$ associ{\'e} {\`a} la norme de
la convergence uniforme sur $[0,t]$ : pour $f,g \in \wwf$,  
$$d_t^{CU}(f,g) = ||f-g||_{[0,t]} = \max\{|f(s)-g(s)|\ ;\ s \in [0,t]\}.$$

\begin{defi} {\bf (Topologie CUCZ sur $\wwf$)}

On appelle topologie de la convergence uniforme sur les compacts
avec contr{\^o}le des z{\'e}ros la topologie associ{\'e}e aux {\'e}carts 
$d_t^{CU}$ et $d_t^{CZ}$ pour $t \ge 0$. 
\end{defi}

Comme ces {\'e}carts sont croissants par rapport {\`a} $t$, les ouverts 
$$V_t(f,\rho,\delta) = \{g \in V :
d^{CU}_t(f,g) < \rho\ ;\ d^{CZ}_t(f,g) < \delta\}.$$ 
pour $f \in \wwf$ fix{\'e} et $t>0$, $\rho>0$, $\delta>0$
forment une base de voisinages de $f$. 

\subsection{Propri{\'e}t{\'e}s de la topologie CUCZ}

Gr{\^a}ce {\`a} la croissance des {\'e}carts $d_t^{CU}$ et $d_t^{CZ}$ par rapport
{\`a} $t$, on voit que la topologie CUCZ est m{\'e}trisable. Nous allons
voir que $\wwf$ est s{\'e}parable pour cette topologie. Nous allons
m{\^e}me montrer un r{\'e}sultat plus pr{\'e}cis. 

\begin{lemm}~\label{Densite des fonctions continues affines par morceaux}
{\bf (Densit{\'e} des fonctions continues affines par morceaux)}

Les fonctions continues, nulles en $0$, affines par morceaux, ayant
des points de subdvision rationnels et prenant des valeurs
rationnelles en ces points forment une partie dense de $\wwf$ pour la
topologie CUCZ. 

Les fonctions continues, positives, nulles en $0$, 
affines par morceaux, ayant des points de subdvision rationnels et 
prenant des valeurs rationnelles en ces points forment une partie 
dense de $\wwf_+$ pour la topologie CUCZ. 
\end{lemm}

\begin{demo}
Soient $f \in \wwf$ et $t>0$, $\rho>0$, $\delta>0$. Quitte {\`a} r{\'e}duire $\delta$, 
on peut supposer que
$$\mbox{Osc}_{[0,t]}(f,\delta) := \sup\{|f(s_1)-f(s_2)|\ ;
\ (s_1,s_2) \in [0,t]^2, |s_1-s_2| \le \delta\} < \rho.$$ 
Choisissons un
nombre fini de z{\'e}ros de $f$ sur $[0,t]$, $0=z_0<\ldots<z_m$, tels que
$$Z_t(f) \subset \bigcup_{k=0}^m ]z_k-\delta,z_k+\delta[$$
et construisons une subdivision $0=t_0<\ldots<t_n=t$ de $[0,t]$ en
intervalles de longueur $\le \delta$, contenant les instants
$z_1<\ldots<z_m$. Soit $g$ l'application obtenue en interpolant
lin{\'e}airement $f$ sur chaque intervalle de subdivision et constante
{\'e}gale {\`a} $f(t)$ sur $[t,+\infty[$. On v{\'e}rifie facilement que
$d_t^{CU}(f,g)<\rho$ gr{\^a}ce {\`a} l'in{\'e}galit{\'e} 
$\mbox{Osc}(f \big|_{[0,t]},\delta) < \rho$. 
Par ailleurs, comme $g$ s'annule en $z_0<\ldots<z_m$, on a 
$$Z_{t-\delta}(f) \subset Z_t(f) \subset Z(g) + ]-\delta,\delta[.$$
Inversement, sur chaque intervalle de subdivision, $g$ est affine et
co{\"\i}ncide avec $f$ aux extr{\'e}mit{\'e}s. Pour que $g$ poss{\`e}de un z{\'e}ro sur un
intervalle de subdivision, il faut que $f$ poss{\`e}de aussi un z{\'e}ro sur cet
intervalle. Comme les intervalles de sudivision sont de longueur $\le
\delta$, on a donc 
$$Z_{t-\delta}(g) \subset Z_t(g) \subset Z(f) + ]-\delta,\delta[,$$
ce qui montre que $d_t^{CZ}(f,g)<\delta$. Ainsi $g \in
V_t(f,\rho,\delta)$. De plus $g$ est positive si $f$ l'est. 

Ces propri{\'e}t{\'e}s restent valables si l'on remplace
les instants de subdivision $t_1<\ldots<t_n$ et les valeurs non nulles
de $g$ {\`a} ces instants par des rationnels proches. 
\hfill $\square$
\end{demo}

\begin{coro} {\bf (Existence d'une base d{\'e}nombrable d'ouverts)}

La topologie CUCZ poss{\`e}de une base d{\'e}nombrable d'ouverts et engendre
la m{\^e}me tribu que la topologie de la convergence uniforme sur les
compacts. 
\end{coro}

\begin{demo}
Soit $D$ une partie partie d{\'e}nombrable dense de $\wwf$. On v{\'e}rifie
facilement que les ouverts $V_t(f,\rho,\delta)$ pour $f \in D$, $t \in
\nnf^*$, $\rho \in \qqf_+^*$, $\delta \in \qqf_+^*$ forment une base
d{\'e}nombrable d'ouverts. 

Il reste {\`a} v{\'e}rifier que ces ouverts sont des bor{\'e}liens pour la
topologie de la convergence uniforme sur les compacts. Cela se voit en
{\'e}crivant que $d_t^{CZ}(f,g)<\delta$ si et seulement si
$$\inf\{|g(s)|\ ;\ s \in [0,t-\delta]\ ;\ d(s,Z(f)) \ge
\delta\} > 0$$ 
et
$$\exists \delta' \in \qqf \cap ]0,\delta[, \sup_{z \in
  Z_{t-\delta}(f)} \inf \{|g(s)|\ ;\ s \in [z-\delta',z+\delta']\} = 0,$$
et en remarquant que les bornes inf{\'e}rieures sur $s$ et sup{\'e}rieure sur
$z$ ci-dessus se ram{\`e}nent par continuit{\'e} {\`a} des bornes sur des
ensembles d{\'e}nombrables denses. 
\hfill $\square$
\end{demo}

Donnons un exemple d'ouvert utile pour la suite. 

\begin{lemm}~\label{Exemple d'ouvert de la topologie CUCZ}
{\bf (Exemple d'ouvert de la topologie CUCZ)}

Pour $b>a \ge 0$, l'ensemble des trajectoires de $\wwf$ poss{\`e}dant au moins 
un z{\'e}ro dans $]a,b[$ est un ouvert de la topologie CUCZ. 
\end{lemm}

\begin{demo}
En effet si $f \in \wwf$ s'annule en $z
\in ]a,b[$, alors toute application $g \in \wwf$ telle que
$d_b^{CZ}(f,g) < \min(z-a,b-z)$ poss{\`e}de au moins un z{\'e}ro dans $]a,b[$. 
\hfill $\square$
\end{demo}

\subsection{Approximation d'une fonction {\`a} partir d'une approximation de 
sa valeur absolue}

Soient $f \in \wwf$. Le lemme ci-dessous montre que si une trajectoire
$r \in \wwf_+$ approche $|f|$ pour la topologie CUCZ, alors en
imposant un nombre fini de signes de $e \in E$, on obtient une
trajectoire $e \cdot r$ approchant $f$.

\begin{lemm} 
{\bf (Comment approcher $f$ {\`a} partir d'une approximation de $|f|$)}

Soient $f \in \wwf$ et $t>0$, $\rho>0$, $\delta>0$ tel que
$$\osc_{[0,t]}(f,\delta) := \sup\{|f(s_1)-f(s_2)|\ ;\ s_1,s_2 \in [0,t]
\ \rm{et}\ |s_1-s_2| \le \delta\} \le \rho.$$ 
Soit $r \in \wwf_+ \cap V_t(|f|,\rho,\delta)$. Alors :
\begin{enumerate}
\item Si $s_1<s_2$ sont deux instants de $[0,t]$ tels que $f(s_1)$ et
  $f(s_2)$ sont de signes oppos{\'e}s et de valeurs absolues strictement
  plus grandes que $\rho$, alors $r$ poss{\`e}de un z{\'e}ro dans
  $]s_1,s_2[$. Par cons{\'e}quent, si $I$ est un intervalle d'excursion de
  $r$, le signe de $f$ est constant sur $\{s \in I \cap [0,t] : |f(s)|
  > \rho\}$.
\item Soit $e \in E$. Pour que $||e \cdot r - f||_{[0,t]} < 5\rho$,
  il suffit que pour tout $I$ intervalle d'excursion de $r$ commenc{\'e}e avant
  $t$ tel que $\sup_{I \cap [0,t]} |r(s)| > 2\rho$, $e(q(I))$ soit
  {\'e}gal au signe de $f$ sur l'ensemble non vide $\{s \in I \cap
  [0,t] : |f(s)| > \rho\}$.
\end{enumerate} 
\end{lemm}

\begin{demo} Montrons les deux points.
\begin{enumerate}
\item 
Montrons le r{\'e}sultat par contraposition. Soient $s_1<s_2$ deux
instants de $[0,t]$ tels que $f(s_1)$ et $f(s_2)$ soient de signes
oppos{\'e}s et de valeurs absolues strictement plus grandes que $\rho$. 
D'apr{\`e}s le th{\'e}or{\`e}me des valeurs interm{\'e}diaires, $f$ 
poss{\`e}de au
moins un z{\'e}ro $z \in ]s_1,s_2[$. Comme $\osc_{[0,t]}(f,\delta) \le \rho$, on
a $s_1 + \delta \le z \le s_2 - \delta \le t - \delta $. Mais comme
$d_t^{CZ}(|f|,r) <\delta$, la trajectoire $r$ poss{\`e}de un z{\'e}ro dans
l'intervalle $]z - \delta , z + \delta[$, qui est inclus dans
$]s_1,s_2[$. Donc $s_1$ et $s_2$ appartiennent {\`a} des intervalles
d'excursion de $r$ disjoints.  
\item Supposons que $e$ est choisi comme ci-dessus. Soit $s \in
  [0,t]$. De deux choses l'une : 
\begin{itemize}
\item soit $r(s)>2\rho$, et alors $|f(s)| > 2\rho - ||f-r||_{[0,t]}
  > \rho$. Par cons{\'e}quent, si $q$ est le rationnel num{\'e}rotant
  l'excursion de $r$ enjambant $s$, le signe $e(q)$ est celui de
  $f(s)$, si bien que 
$$\big|(e \cdot r)(s)- f(s)\big| = \big|r(s)-|f(s)|\big| < \rho.$$ 
\item soit $r(s) \le 2\rho$, et alors $|f(s)| \le 2\rho + ||f-r||_{[0,t]}
  < 3 \rho$, d'o{\`u}  
$$\big|(e \cdot r)(s)- f(s)\big| \le r(s) + |f(s)| < 5\rho.$$ 
\end{itemize}
Dans tous les cas, $\big|(e \cdot r)(s)- f(s)\big| < 5\rho$.  
\end{enumerate} 
\hfill $\square$
\end{demo}

\begin{coro}~\label{Ce qu'il faudra montrer}
{\bf (Passage de la cha{\^\i}ne $(R_n)_{n \in \nnf}$ {\`a} la cha{\^\i}ne $(W_n)_{n \in \nnf}$)}

Pour montrer que l'orbite de preque toute trajectoire $w \in \wwf$
est dense $\wwf$ pour la topologie CUCZ, il suffit de montrer
que pour tout ouvert $U$ de $\wwf_+$ pour la topologie CUCZ,
$$\sup_{n \in \nnf} P[R^{(n)} \in U|R^{(0)}] > 0\ \text{ p.s.}$$
\end{coro}

\begin{demo}
En effet, gr{\^a}ce au 
corollaire~\ref{Condition suffisante de densite des orbites} 
et {\`a} l'existence d'une base d{\'e}nombrable d'ouverts, il
suffit de montrer que pour tout $f \in \wwf_+$,
$t>0$, $\rho>0$, $\delta>0$,  
$$\sup_{n \in \nnf} P[W^{(n)} \in V_t(f,5\rho,\delta)|W^{(0)}] > 0\
\text{ p.s.}$$
Quitte {\`a} r{\'e}duire $\delta$, on peut supposer que
$\osc_{[0,t]}(f,\delta) \le \rho$. 

Soit $n \in \nnf$. Alors $W^{(n+1)} = \varepsilon_{n+1} \cdot
R^{(n)}$. Notons $\ic(R^{(n)})$ l'ensemble des
intervalles $I$ d'excursion de $R^{(n)}$ commen{\c c}ant avant $t$ tels que
$\sup_{I \cap [0,t]} |R^{(n)}(s)| > 2\rho$. L'ensemble al{\'e}atoire
$\ic(R^{(n)})$ est fini par continuit{\'e} des trajectoires de $R^{(n)}$. 

D'apr{\`e}s le lemme pr{\'e}c{\'e}dent, pour que  l'{\'e}v{\'e}nement $[W^{(n+1)} \in
V_t(f,5\rho,\delta)]$ soit r{\'e}alis{\'e}, il suffit que $R^{(n)} \in
V_t(|f|,\rho,\delta)$ et que pour tout $I \in \ic(R^{(n)})$,
$\varepsilon_{n+1}(q(I))$ soit {\'e}gal au signe de $f$ sur $\{s \in I
\cap [0,t] : |f(s)| > \rho\}$. Comme $\varepsilon_{n+1}$ est
ind{\'e}pendante de $R^{(n)}$ et suit la loi uniforme sur
$\{-1,1\}^{\qqf_+^*}$, on a donc 
$$P[W^{(n+1)} \in V_t(f,\rho,\delta) | R^{(n)}] \ge {}^\iif {[R^{(n)} \in
V_t(|f|,\rho,\delta)]}\ \Big(\frac{1}{2}\Big)^{\card\ \ic(R^{(n)})}.$$
Par cons{\'e}quent, si $P[R^{(n)} \in V_t(|f|,\rho,\delta) | R^{(0)}] >
0$ presque s{\^u}rement, alors 
$$P[W^{(n+1)} \in V_t(f,5\rho,\delta) | R^{(0)}]>0 \mbox{ presque
  s{\^u}rement.}$$ 
On obtient ainsi le r{\'e}sultat voulu en remarquant que $\sigma(R^{(0)})
= \sigma(W^{(0)})$ puisque $R^{(0)}
= |W^{(1)}|$ et $W^{(0)} = \widehat{W^{(1)}}$ (voir les rappels du
paragraphe 2.1). \hfill $\square$ 
\end{demo}

Un des int{\'e}r{\^e}ts de travailler avec la cha{\^\i}ne de Markov $(R^{(n)})_{n
  \in \nnf}$ plut{\^o}t qu'avec la cha{\^\i}ne de Markov $(W^{(n)})_{n
  \in \nnf}$ est qu'il est possible de pr{\'e}server, ou presque, le d{\'e}but
des trajectoires par un choix convenable des familles de signes. 

En effet, la famille de signes ${\bf 1} \in E$, constante {\'e}gale {\`a} $1$,
v{\'e}rifie $F({\bf 1},r) = r - \underline{r} = r$
pour tout $r \in \wwf_+$. La difficult{\'e} est qu'on ne peut
imposer qu'un nombre fini de signes dans les familles
$\varepsilon_1,\varepsilon_2,\ldots$ pour avoir des probabilit{\'e}s
strictement positives. Nous allons donc utiliser le fait que $F(e,r)$
est proche de $r$ en norme uniforme sur un intervalle $[0,t]$ d{\`e}s que $e(q)=1$
pour tout rationnel $q$ num{\'e}rotant une excursion de $r$ de hauteur
$\ge \eta$ commen{\c c}ant avant $t$, avec $\eta > 0$ petit. 

Cela nous am{\`e}ne {\`a} {\'e}tablir des propri{\'e}t{\'e}s de continuit{\'e} presque partout. 

\section{R{\'e}sultats de continuit{\'e} presque partout}

Sauf mention explicite du contraire, on munira toujours $E =
\{-1,1\}^{\qqf_+^*}$ de la topologie produit et $\wwf_+$ de la
topologie induite par la topologie de la convergence uniforme 
sur les compacts avec contr{\^o}le des z{\'e}ros. 

\subsection{Continuit{\'e} de l'action de $E$ sur $\wwf$}

Dans ce paragraphe, nous allons montrer la continuit{\'e} de l'application
$(e,w) \mapsto e \cdot w$ de $E \times \wwf$ dans $\wwf$ pour la
topologie CUCZ pourvu qu'on se restreigne aux trajectoires sans z{\'e}ro
rationnel. 

\begin{lemm}~\label{continuite de la numerotation}{\bf (Continuit{\'e} de
    la num{\'e}rotation des excursions)}

Soit $w_0 \in \wwf$ une trajectoire sans z{\'e}ro rationnel autre que $0$. 
Soit $t>0$
tel que $w_0(t) \ne 0$. Il existe deux r{\'e}els $T>t$ et $\delta > 0$ tels
que pour tout $w \in \wwf$, 
$$d^{CZ}_T(w_0,w) < \delta \Rightarrow Q_t(w) = Q_t(w_0).$$ 
\end{lemm}

\begin{demo}
Notons $q_0 = Q_t(w_0)$ et $]a,b[ = I_t(w_0)$ l'intervalle d'excursion de 
$w_0$ enjambant $t$. Soit $F$ l'ensemble (fini) des {\'e}l{\'e}ments de 
$\qqf_+^*$ pr{\'e}c{\'e}dant
$q_0$ pour l'ordre introduit au d{\'e}but du paragraphe 3.2, y
compris $q_0$. Pour que $Q_t(w)=q_0$, il faut et il suffit que $q_0$
soit le seul {\'e}l{\'e}ment de $F$ dans $I_t(w)$. Cela se produit d{\`e}s que
$d^{CZ}_{b+\delta}(w_0,w) < \delta$ o{\`u} 
$$\delta = \min_{q \in F \cup\{t\}} \min(|q-a|,|q-b|).$$
En effet, l'in{\'e}galit{\'e} $d^{CZ}_{b+\delta}(w_0,w) < \delta$
assure que les bornes de l'intervalle d'excurion $I_t(w)$ sont dans
$]a-\delta,a+\delta[$ et $]b-\delta,b+\delta[$. \hfill $\square$ 
\end{demo}

\begin{lemm}{\bf (Comparaison des modules de continuit{\'e} de $w$ et $e
    \cdot w$)}

Soient $e \in E$, $w \in \wwf$ et $t>0$. Pour tout $\delta > 0$, notons
$${\rm Osc}_{[0,t]}(w,\delta) = \sup\{|w(s_1)-w(s_2)|\ ;\ s_1,s_2 \in
[0,t]\ \rm{et}\ |s_1-s_2| \le \delta\}.$$
Alors 
${\rm Osc}_{[0,t]}(e \cdot w,\delta) \le
2\ {\rm Osc}_{[0,t]}(w,\delta)$.
\end{lemm}

\begin{demo}
Il suffit de remarquer que $|e \cdot w| = |w|$ et d'utiliser les 
in{\'e}galit{\'e}s
$${\rm Osc}(|w|) \le {\rm Osc}(w) \le {\rm Osc}(w_+)  + {\rm Osc}(w_-) 
\le 2{\rm Osc}(|w|)$$
dans lesquelles on a omis d'{\'e}crire l'intervalle $[0,t]$ et la variable 
$\delta$.
%
%
\hfill $\square$ 
\end{demo}

\begin{prop}~\label{Continuite de l'action} 
{\bf (Continuit{\'e} de l'action de $E$ sur $\wwf$)}

Si l'on munit $E$ de la topologie produit et $\wwf$ de la topologie CUCZ,  
l'application $(e,w) \mapsto e \cdot w$ de $E \times
\wwf$ dans $\wwf$ est continue en tout couple $(e_0,w_0)$ o{\`u} $e_0 \in
E$ et $w_0 \in \wwf$ est une trajectoire sans z{\'e}ro rationnel.
\end{prop}

\begin{demo}
Soient $(e_n)_{n \ge 1}$ une suite convergant vers $e_0$ pour la
topologie produit et $(w_n)_{n \ge 1}$ une suite convergant vers $w_0$ pour la
topologie CUCZ. 
Comme les trajectoires $e_n \cdot w_n$ et $w_n$ ont les m{\^e}mes z{\'e}ros,
il suffit de montrer que la suite $(e_n \cdot w_n)_{n \ge 1}$ converge
uniform{\'e}ment sur les compacts. 

Comme la suite $(w_n)_{n \ge 1}$ converge vers $w_0$ uniform{\'e}ment sur
tout segment $[0,t]$, elle est uniform{\'e}ment {\'e}quicontinue sur
$[0,t]$. D'apr{\`e}s le lemme pr{\'e}c{\'e}dent, il en est de m{\^e}me pour la suite
$(e_n \cdot w_n)_{n \ge 1}$. Il suffit donc de v{\'e}rifier que $(e_n
\cdot w_n)(t) \to (e_0 \cdot w_0)(t)$ pour $t \ge 0$ fix{\'e}. 

Si $w_0(t)=0$, alors on remarque simplement que $|(e_n \cdot w_n)(t)| =
|w_n(t)| \to |w_0(t)| = 0$.

Si $w_0(t) \ne 0$, on montre la convergence de $(e_n(Q_t(w_n)))_{n \ge
  1}$. Cette convergence d{\'e}coule de la convergence ponctuelle de
$(e_n)_{n \ge 1}$ et du fait que $Q_t(w_n) = Q_t(w_0)$ {\`a} partir d'un
certain rang, gr{\^a}ce au lemme~\ref{continuite de la
  numerotation}. \hfill $\square$
\end{demo}

\subsection{Continuit{\'e} de l'application $w \mapsto w-\underline{w}$ de
$\wwf$ dans $\wwf_+$}

Pour $w \in \wwf$, on note $N(w)$ l'ensemble des instants de records
n{\'e}gatifs de $w$, qui est aussi l'ensemble des z{\'e}ros de $w-\underline{w}$ 
$$N(w) = \{t \in \rrf_+ : \underline{w}(t) = w(t)\} =
Z(w-\underline{w}).$$

\begin{prop}~\label{Difference avec le minimum courant} 
{\bf (Continuit{\'e} de l'application $w \mapsto w-\underline{w}$ de
$\wwf$ dans $\wwf_+$)}

Soit $w_0 \in \wwf$. Si pour tout $z \in N(w_0)$ et $\delta>0$,  
$\underline{w_0}(z-\delta) > \underline{w_0}(z+\delta)$ avec la
convention $\underline{w_0}(t) = 0$ pour $t < 0$,  
alors  l'application $w \mapsto w-\underline{w}$, de $\wwf$
muni de la topologie de la convergence uniforme sur les compacts dans 
$\wwf_+$ muni de la topologie CUCZ, est continue en $w_0$. 
\end{prop}

\begin{demo}
Soient $w_0 \in \wwf$ v{\'e}rifiant les hypoth{\`e}ses et $w \in \wwf$. Notons
$r_0 = w_0 - \underline{w_0}$ et $r = w-\underline{w}$. 
Soient $t>0$, $\rho>0$ et $\delta \in ]0,t[$. 
Pour avoir $r \in V_t(r_0,\rho,\delta)$, il suffit que 
\begin{enumerate}
\item $||r-r_0||_{[0,t]} < \rho$ ;
\item $Z_{t-\delta}(r) \subset Z_t(r_0)\ +\ ]-\delta,\delta[$ ;
\item $Z_{t-\delta}(r_0) \subset Z_t(r)\ +\ ]-\delta,\delta[$.
\end{enumerate}
Montrons que ces trois conditions sont v{\'e}rifi{\'e}es si
$||w-w_0||_{[0,t]}$ est suffisament petit. 

On commence par remarquer que 
$||\underline{w}-\underline{w_0}||_{[0,t]} \le ||w-w_0||_{[0,t]}$,
d'o{\`u} $$||r-r_0||_{[0,t]} \le 2 ||w-w_0||_{[0,t]}.$$
Donc la condition $1$ est r{\'e}alis{\'e}e d{\`e}s que $||w-w_0||_{[0,t]} < \rho/2$. 

Par continuit{\'e} de $r_0$ et compacit{\'e} de $K = \{s \in
  [0,t] : d(s,Z_t(r_0)) \ge \delta\}$, la borne 
$\alpha = \inf\{r_0(s)\ ; s \in K\}$ est atteinte et strictement
positive. La condition $2$ est r{\'e}alis{\'e}e d{\`e}s
que $||w-w_0||_{[0,t]} < \alpha/2$ puisque cette in{\'e}galit{\'e} entra{\^\i}ne 
$||r-r_0||_{[0,t]} < \alpha$, ce qui interdit {\`a} $r$ de s'annuler sur
$K$.   
 
Enfin, l'hypoth{\`e}se faite sur $w_0$, la continuit{\'e} de $\underline{w_0}$ et
compacit{\'e} de $Z_t(r_0)$, assurent que la borne 
$$\beta = \inf\{\underline{w_0}(z-\delta)-\underline{w_0}(z+\delta)\ ; 
z \in Z_t(r_0)\}$$ 
est strictement positive. La condition $2$ est r{\'e}alis{\'e}e d{\`e}s
que $||w-w_0||_{[0,t]} < \beta/2$ puisque cette in{\'e}galit{\'e} entra{\^\i}ne,
pour tout $z \in Z_t(r_0)$,
$$\underline{w}(z-\delta)-\underline{w}(z+\delta) >
\underline{w_0}(z-\delta)-\underline{w_0}(z+\delta) - \beta \ge 0,$$
ce qui implique l'existence d'un record n{\'e}gatif de $w$, donc d'un z{\'e}ro
de $r$ dans l'intervalle $]z-\delta,z+\delta[$. \hfill $\square$
\end{demo}

Par composition, on obtient ainsi la continuit{\'e} presque partout de
$F$. 
 
\begin{prop}{\bf (Continuit{\'e} presque partout de $F$)}

L'application $F : (e,r) \mapsto e \cdot r - \underline{e \cdot
  r}$ de $E \times \wwf_+$ dans $\wwf_+$ est continue presque partout  
pour $P_\varepsilon \otimes P_R$ o{\`u} $P_\varepsilon$ est la loi uniforme sur
$E = \{-1,1\}^{\qqf_+^*}$ et $P_R$ la loi du mouvement brownien 
r{\'e}fl{\'e}chi. 
\end{prop}

\subsection{D{\'e}finition de la transformation $F$ apr{\`e}s un instant $a \ge 0$}

Pour pouvoir pr{\'e}server le d{\'e}but des trajectoires, nous allons d{\'e}finir
une transformation $F_a$ de $E \times \wwf_+$ dans $\wwf_+$ qui se
comporte comme l'identit{\'e} de $\wwf_+$ avant un instant $a \ge 0$ et
comme $F$ ensuite. Nous avons besoin d'introduire quelques notations.   

\begin{defi}~\label{hybridation} {\bf (Familles hybrides de signes)}

Pour $t \ge 0$ et $e_1,e_2 \in E = \{-1,1\}^{\qqf_+^*}$, on d{\'e}finit
la famille hybride \og $e_1$ puis $e_2$ {\`a} l'instant $t$\fg not{\'e}e $\displaystyle{e_1
\mathop{\smile}^t e_2}$ par  
$$(e_1 \mathop{\smile}^t e_2) (q) = e_1(t) \mbox{ si } q \le t,$$ 
$$(e_1 \mathop{\smile}^t e_2) (q) = e_2(t) \mbox{ si } q > t .$$
\end{defi}

\begin{defi} {\bf (D{\'e}finition de $D_a$ et $F_a$)}

Pour $a \ge 0$ et $w \in \wwf$, on note $D_a(w)= \inf(Z(w) \cap [a,\infty[)$ 
le premier z{\'e}ro de $w$ apr{\`e}s l'instant $a$. 

On d{\'e}finit la transformation \og $F$ apr{\`e}s $a$\fg de $E \times \wwf_+$
dans $\wwf_+$ en posant 
$$F_a(e,r) = F({\bf 1} \mathop{\smile}^{D_a(r)} e,r).$$
Autrement dit, 
\begin{eqnarray*} 
&\mbox{pour}\ t \le D_a(r),\ F_a(e,r)(t) = r(t),\\ 
&\mbox{pour}\ t > D_a(r),\ F_a(e,r)(t) = (e \cdot r)(t) -
\min_{[D_a(r),t]} (e \cdot r). 
\end{eqnarray*}
\end{defi}

On remarque que $F_0=F$ puisque $D_0(r)=0$ pour tout $r \in \wwf_+$. 
Voyons quelques propri{\'e}t{\'e}s imm{\'e}diates de la transformation $F_a$. 

\begin{lemm}~\label{Preservation de $D_a$ par $F_a$} 
{\bf (Pr{\'e}servation de $D_a$ par $F_a$)}

Quels que soient $e \in E$ et $r \in \wwf_+$, $D_a(F_a(e,r))=D_a(r)$. 
\end{lemm}

\begin{demo}
Par construction, la trajectoire $F_a(e,r)$ co{\"\i}ncide avec $r$ sur
$[0,D_a(r)]$, donc $D_a(r)$ est le premier z{\'e}ro de $F_a(e,r)$ sur
l'intervalle $[a,+\infty[$.
\hfill $\square$.
\end{demo}

Nous allons voir que la transformation $F_a$ dans laquelle les signes
sont choisis au hasard se comporte apr{\`e}s l'instant $D_a$ comme la
transformation $F$. Pour donner un sens pr{\'e}cis {\`a} cette affirmation,
nous avons besoin d'introduire l'op{\'e}rateur de translation $\theta_{D_a}$. 

\begin{defi} {\bf (Op{\'e}rateur de translation $\theta_{D_a}$)}

Pour tout $r \in \wwf_+$ tel que $D_a(r)<+\infty$, on note
$\theta_{D_a}(r)$ la trajectoire de $\wwf_+$ d{\'e}finie par
$$\theta_{D_a}(r)(t) = r(D_a(r)+t) \mbox{ pour } t \ge 0.$$ 
\end{defi}

Par r{\'e}currence du mouvement brownien r{\'e}fl{\'e}chi, l'op{\'e}rateur
$\theta_{D_a}$ est d{\'e}fini $P_R$-presque partout et la propri{\'e}t{\'e} forte
de Markov nous dit que $\theta_{D_a}$ pr{\'e}serve la loi $P_R$. 

\begin{lemm}~\label{Lien entre $F$ et $F_a$} 
{\bf (Lien entre $F$ et $F_a$)}

Soient $r \in \wwf_+$ tel que $D_a(r)<+\infty$ et $\varepsilon$ une
variable al{\'e}atoire de loi uniforme sur $E$. Alors 
$$\theta_{D_a}(F_a(\varepsilon,r)) \mbox{ a m{\^e}me loi que }
F(\varepsilon,\theta_{D_a}(r)).$$
\end{lemm}

\begin{demo}
Par d{\'e}finition de $F$, $F_a$ et de $\theta_{D_a}$, on a pour tout $t
\ge 0$ : 
$$F(\varepsilon,\theta_{D_a}(r))(t) = \varepsilon(Q_t(\theta_{D_a}(r)))
\ r(D_a(r)+t),$$
et gr{\^a}ce au fait que $D_a(F_a(e,r))=D_a(r)$,
$$\theta_{D_a}(F_a(\varepsilon,r)) (t) = \varepsilon(Q_{D_a(r)+t}(r))
\ r(D_a(r)+t).$$
Il suffit donc de montrer que  
$$(\varepsilon(Q_{D_a(r)+t}(r)))_{t \ge 0} \mbox{ a m{\^e}me loi que
} (\varepsilon(Q_t(\theta_{D_a}(r))))_{t \ge 0}.$$
Soit $s$ une permutation (d{\'e}terministe) de $\qqf_+$ telle que
$s(0)=0$ et telle que pour tout rationnel num{\'e}rotant une excursion de
$\theta_{D_a}(r)$, $s(q)$ soit le rationnel num{\'e}rotant l'excursion
correspondante de $r$. Par construction de $s$, on a donc
$Q_{D_a(r)+t}(r) = s(Q_t(\theta_{D_a}(r)))$ pour tout $t \ge
0$. Le r{\'e}sultat d{\'e}coule du fait que $\varepsilon \circ s$ a m{\^e}me loi
que $\varepsilon$. \hfill $\square$
\end{demo}

Le r{\'e}sultat ci-dessous montre que la transformation $F_a$ o{\`u} les signes sont
choisis au hasard pr{\'e}serve la loi du mouvement brownien r{\'e}fl{\'e}chi. 

\begin{coro}~\label{preservation de la loi du mouvement brownien reflechi}
{\bf (Effet de la transformation $F_a$)}

Soient $R$ un mouvement brownien r{\'e}fl{\'e}chi et $\varepsilon$ une
variable al{\'e}atoire ind{\'e}pendante de $R$ et de loi uniforme sur
$E=\{-1,1\}^{\qqf_+^*}$. Alors 
\begin{enumerate}
\item Le processus $\theta_{D_a}(F_a(\varepsilon,R))$ est ind{\'e}pendant
  de $(R_t)_{0 \le t \le D_a(R)}$ et a m{\^e}me loi que
  $F(\varepsilon,\theta_{D_a}(R))$. 
\item Le processus $F_a(\varepsilon,R)$ est un mouvement brownien
r{\'e}fl{\'e}chi qui co{\"\i}ncide avec $R$ jusqu'{\`a} l'instant $D_a(R)$. 
\end{enumerate}
\end{coro}

\begin{demo}
D'apr{\`e}s le lemme pr{\'e}c{\'e}dent et par ind{\'e}pendance de $\varepsilon$ et
de $R$, 
$$\lc \big( \theta_{D_a}(F_a(\varepsilon,R)) \big| R \big) = \lc \big(
F(\varepsilon,\theta_{D_a}(R)) \big| R \big).$$
On en d{\'e}duit, gr{\^a}ce {\`a} la propri{\'e}t{\'e} de Markov et {\`a} la pr{\'e}servation de
la loi du mouvement brownien r{\'e}fl{\'e}chi par $F(\varepsilon,\cdot)$ et
$\theta_{D_a}$,  
$$\lc \big( \theta_{D_a}(F_a(\varepsilon,R)) \big| (R_t)_{0 \le t \le
  D_a(R)} \big) = \lc \big( F(\varepsilon,\theta_{D_a}(R)) \big) = P_R,$$
ce qui montre le point 1. 

On en d{\'e}duit le point 2 par la propri{\'e}t{\'e} de Markov en remarquant que
le processus $F_a(\varepsilon,R)$ s'obtient en concat{\'e}nant $(R_t)_{0
  \le t \le D_a(R)}$ et $\theta_{D_a}(F_a(\varepsilon,R))$.
\hfill $\square$
\end{demo}

\subsection{Continuit{\'e} presque partout de la transformation $F_a$}

Nous allons montrer que $F_a$ est continue presque partout sur $E
\times \wwf_+$. 

\begin{lemm}~\label{continuite de D}{\bf (Continuit{\'e} de
    $D_a$ et de $\theta_{D_a}$)}

Soit $w_0 \in \wwf_+$. Si l'instant $a$ n'est pas le d{\'e}but d'une 
excursion de $w_0$, l'application 
$D_a : \wwf_+ \to [0,+\infty]$ est continue en $w_0$. 
Si de plus $D_a(w_0)<+\infty$, alors la fonction 
$\theta_{D_a} : \wwf_+ \to \wwf_+$ est continue en $w_0$. 
\end{lemm}
\begin{demo}
Pour montrer la continuit{\'e} de $D_a$, on distingue trois cas.  
\begin{enumerate}
\item Cas o{\`u} $D_a(w_0)=a$ (autrement dit $w_0(a)=0$). 

Soit $b>a$. Comme $a$ n'est pas un d{\'e}but d'excursion de $w_0$, $w_0$
poss{\`e}de au moins un z{\'e}ro dans $]a,b[$. L'ensemble 
$O_{a,b} = \{w \in \wwf_+ : w \mbox{ poss{\`e}de un z{\'e}ro dans } ]a,b[\}$ 
est donc un voisinage de $w_0$ 
(gr{\^a}ce au lemme~\ref{Exemple d'ouvert de la topologie CUCZ}) et pour 
tout $w \in O_{a,b}$, on a $a \le D_a(w) < b$. 

\item Cas o{\`u} $a<D_a(w_0)<+\infty$. 

Notons $d = D_a(w_0)$ et notons $g$
  le dernier z{\'e}ro de $w_0$ avant l'instant $a$. Alors $g<a<d$. 
Soit $\delta \in ]0,\min(a-g,d-a)]$. Pour tout $w \in \wwf$
tel que $d_{d+\delta}^{CZ}(w_0,w)<\delta$, $w$ poss{\`e}de un z{\'e}ro dans
$]d-\delta,d+\delta[$ mais n'en poss{\`e}de pas dans
$[g+\delta,d-\delta]$, donc $d-\delta < D_a(w) < d+\delta$. 

\item Cas o{\`u} $D_a(w_0)=+\infty$. 

Notons $g$ le dernier z{\'e}ro de $w_0$ avant l'instant $a$. Alors $g<a$. 
Soit $t>a$. Pour tout $w \in \wwf$
tel que $d_{t+a-g}^{CZ}(w_0,w)<a-g$, $w$ ne poss{\`e}de pas de z{\'e}ro dans
$[a,t]$, donc $D_a(w) > t$. 
\end{enumerate}
Dans les trois cas, $D_a$ est continue en $w_0$.   

Pour montrer la continuit{\'e} de $\theta_{D_a}$, on {\'e}tend la topologie CUCZ {\`a} 
$\tilde{\wwf} = {\cal C}(\rrf_+,\rrf)$ en posant 
$Z(w) = \{t \in \rrf_+ : w(t)=0\} \cup \{0\}$ pour tout $w \in \wwf$ 
puis en d{\'e}finissant les {\'e}carts $d_t^{CU}$ et $d_t^{CU}$ comme sur $\wwf$. 
Le fait de mettre syst{\'e}matiquement $0$ dans $Z(w)$ ne change rien pour 
les trajectoires qui s'annulent en $0$ mais permet d'{\'e}viter les probl{\`e}mes de 
bord en $0$. 

Sur l'ensemble des trajectoires $w$ telles que $a$ n'est pas un d{\'e}but 
d'excursion de $w$ et $D_a(w)<+\infty$, on {\'e}crit alors $\theta_{D_a}$ 
comme la compos{\'e}e de l'application $w \mapsto (\theta_{D_a}(w),w)$ 
et de $(t,w) \mapsto \theta_t(w)=w(t+\cdot)$. La premi{\`e}re de ces 
applications est continue par continuit{\'e} de $D_a$. La continuit{\'e} de 
la seconde d{\'e}coule des in{\'e}galit{\'e}s suivantes pour 
$0 \le a \le b \le T$, $u,v \in \wwf$ et $t>0$ :
\begin{eqnarray*}
d_t^{CU}(\theta_a(u),\theta_b(v)) 
&\le& d_t^{CU}(\theta_a(u),\theta_a(v))\ +\ d_t^{CU}(\theta_a(v),\theta_b(v))\\
&\le& d_{T+t}^{CU}(u,v)\ +\ \osc_{[0,T+t]}(v,b-a)
\end{eqnarray*}
et
\begin{eqnarray*}
d_t^{CZ}(\theta_a(u),\theta_b(v)) 
&\le& d_t^{CZ}(\theta_a(u),\theta_a(v))\ +\ d_t^{CZ}(\theta_a(v),\theta_b(v))\\
&\le& d_{T+t}^{CZ}(u,v)\ +\ b-a.
\end{eqnarray*}
D'o{\`u} le r{\'e}sultat annonc{\'e}. \hfill $\square$.
\end{demo}

Notons $W_{+,a}$ l'ensemble des trajectoires de $\wwf$ qui sont
positives sur $[0,a]$. De la continuit{\'e} de $D_a$ et de la continuit{\'e}
de l'action de $E$ sur $\wwf$ (proposition~\ref{Continuite de l'action}), on
d{\'e}duit imm{\'e}diatement le corollaire suivant. 

\begin{coro} {\bf (Continuit{\'e} de l'action apr{\`e}s $D_a$ de $E$ sur $\wwf$)}

Soient $e_0 \in E$ et $w_0 \in \wwf$. Si l'instant $a$ n'est pas le
d{\'e}but d'une excursion de $w_0$ et $D_a(w_0) \notin \qqf$,
les applications 
$$\displaystyle{(e,r) \mapsto {\bf 1} \mathop{\smile}^{D_a(r)} e}
\mbox{ de } E \times \wwf_+ \mbox{ dans } E,$$ 
$$\displaystyle{(e,r) \mapsto ({\bf 1} \mathop{\smile}^{D_a(r)} e)
  \cdot r} \mbox{ de } E \times \wwf_+ \mbox{ dans } \wwf_{+,a}$$ 
sont continues en $(e_0,w_0)$.
\end{coro}

La proposition ci-dessous se d{\'e}montre comme la
proposition~\ref{Difference avec le minimum courant}.  

\begin{prop} {\bf (Continuit{\'e} de l'application $w \mapsto w-\underline{w}$ de
$\wwf_{+,a}$ dans $\wwf_+$)}

Soit $w_0 \in \wwf_{+,a}$ telle que pour tout $z \in N(w_0) \cap
[a,+\infty[$ et $\delta>0$, $\underline{w_0}(z-\delta) >
\underline{w_0}(z+\delta)$ avec la convention $\underline{w_0}(t) = 0$
pour $t < 0$. L'application $w \mapsto w-\underline{w}$ de $\wwf_{+,a}$ 
muni de la topologie de la convergence uniforme sur les compacts dans 
$\wwf_+$ muni de la topologie CUCZ est continue en $w_0$. 
\end{prop}

Par composition, on obtient finalement le continuit{\'e} presque partout
de $F_a$. 

\begin{prop} {\bf (Continuit{\'e} presque partout de $F_a$)} 

L'application 
$$F_a : (e,r) \mapsto \displaystyle{(1
  \mathop{\smile}^{D_a(r)} e) \cdot r - \underline{(1
  \mathop{\smile}^{D_a(r)} e) \cdot r}}$$ 
de $E \times \wwf_+$ dans $\wwf_+$ est continue presque partout.  
\end{prop}

\section{Accessibilit{\'e}s d'ouverts}

Dans toute cette partie, nous introduisons la notion d'accessibilit{\'e}
par $F_a$ qui est un outil essentiel de la d{\'e}monstration. 

On garde les notations introduites dans la partie~\ref{construction de
chaines de Markov}. 
On notera $P_R$ la loi du mouvement brownien r{\'e}fl{\'e}chi (qui est donc
une probabilit{\'e}s sur $\wwf_+$) et $P_\varepsilon$ la loi uniforme sur
$E = \{-1,1\}^{\qqf_+^*}$. 

\begin{defi} {\bf (Accessibilit{\'e} par $F_a$)}

Soient $a > 0$, $r \in \wwf_+$ et $B$ une partie mesurable de $\wwf_+$. 
On dira que $B$ est accessible par $F_a$ depuis $r$ s'il existe $n \in
\nnf$ tel que 
$$P[F_a(\varepsilon_n,\cdot) \circ \cdots \circ
F_a(\varepsilon_1,\cdot) (r) \in B]>0.$$
\end{defi}

Le corollaire~\ref{Ce qu'il faudra montrer} nous dit que pour
d{\'e}montrer le th{\'e}or{\`e}me~\ref{Densite pour la topologie CUCZ}, il
suffit de v{\'e}rifier que tout ouvert non-vide de $\wwf_+$ est accessible
par $F_0=F$ depuis $P_R$-presque toute trajectoire $r \in \wwf_+$ . 


Nous allons d{\'e}montrer que si $b \ge a \ge 0$, l'accessibilit{\'e}
d'un ouvert (ou m{\^e}me d'un presque-ouvert, notion d{\'e}finie {\`a} la 
sous-section~\ref{presque-ouvert}) par la transformation $F_b$ 
entra{\^\i}ne son accessibilit{\'e} par $F_a$.

\subsection{Cons{\'e}quences de la continuit{\'e} presque partout de $F_a$} 

Pour tout $a \ge 0$ et tout bor{\'e}lien $B$ de $\wwf_+$, notons 
$A_{a,1}(B)$ l'ensemble des trajectoires de $\wwf_+$ d'o{\`u} l'on acc{\`e}de par
$F_a$ en un coup avec probabilit{\'e} strictement positive : 
$$A_{a,1}(B) = \{r \in \wwf_+ : P[F_a(\varepsilon_1,r) \in B] > 0\}.$$
Notons $F_a^{-1}(B)$ l'image r{\'e}ciproque de $B$ par $F_a$ : 
$$F_a^{-1}(B) = \{(e,r) \in E \times  \wwf_+ : F_a(e,r) \in B\}.$$ 
Notons $O_{a,1}(B)$ la projection de l'int{\'e}rieur de $F_a^{-1}(B)$ sur $\wwf_+$ : 
$$O_{a,1}(B) = \{r \in \wwf_+ : \exists e \in E, (e,r) \in
\big(F_a^{-1}(B)\big)^\circ\}.$$

\begin{lemm}~\label{Consequence de la continuite}
{\bf (Cons{\'e}quence de la continuit{\'e} de $F_a$)}

Fixons $a \ge 0$. Soit $V$ un ouvert de $\wwf_+$ pour la topologie CUCZ. 
Alors $O_{a,1}(V)$ est un ouvert contenu dans $A_{a,1}(V)$ et $P_R[A_{a,1}(V)
\setminus O_{a,1}(V)] = 0$. 
\end{lemm}

\begin{demo}
Le fait que $O_{a,1}(V)$ soit un ouvert contenu dans $A_{a,1}(V)$ est
imm{\'e}diat. En effet, si $r_0 \in O_{a,1}(V)$, alors on peut choisir
$e_0 \in E$ tel que $(e_0,r_0) \in \big(F_a^{-1}(V)\big)^\circ$ et :
\begin{itemize}
\item pour tout $r$ dans un certain voisinage de $r_0$, $(e_0,r) \in
\big(F_a^{-1}(V)\big)^\circ$ donc $r \in O_{a,1}(V)$ ; 
\item l'ensemble des $e \in E$ tels que $(e,r_0) \in F_a^{-1}(V)$ est un
voisinage de $e_0$, donc est de mesure positive pour $P_\varepsilon$.   
\end{itemize}
Par ailleurs, notons 
$$A = \{(e,r) \in E \times (A_{a,1}(V) \setminus O_{a,1}(V)) :
F_a(e,r) \in V\}.$$  
Alors par d{\'e}finition de $O_{a,1}(V)$, 
$$A \subset F_a^{-1}(V) \setminus \big(F_a^{-1}(V)\big)^\circ.$$
Donc $A$ est contenu dans l'ensemble des points de discontinuit{\'e} de
$F_a$. Comme $F_a$ est continue presque partout sur $E \times \wwf_+$, on
a donc 
$$0 = (P_\varepsilon \otimes P_R) (A) = \int_{A_{a,1}(V) \setminus O_{a,1}(V)}
P[F_a(\varepsilon_1,r) \in V]\ P_R(dr).$$
Mais $P[F_a(\varepsilon,r) \in V] > 0$ pour tout $r \in A_{a,1}(V)$,
donc $P_R[A_{a,1}(V) \setminus O_{a,1}(V)] = 0$. \hfill $\square$ 
\end{demo}

Signalons quelques propri{\'e}t{\'e}s imm{\'e}diates de l'application qui {\`a} un 
bor{\'e}lien $B$ associe le bor{\'e}lien $A_{a,1}(B)$. Ces propri{\'e}t{\'e}s nous 
seront utiles par la suite.

\begin{lemm}~\label{Proprietes de $A_{a,1}$} 
{\bf (Propri{\'e}t{\'e}s de $A_{a,1}$)}
\begin{enumerate}
\item Pour toute suite $(B_n)_{n \in \nnf}$ de bor{\'e}liens de $\wwf$, 
$$A_{a,1} \Big( \bigcup_{n \in \nnf} B_n \Big) = \bigcup_{n \in \nnf}
A_{a,1}(B_n)$$

\item Si $P_R(B)=0$, alors $P_R(A_{a,1}(B))=0$. 

\item Si $B_1 \subset B_2$ $P_R$-presque s{\^u}rement, alors
$A_{a,1}(B_1) \subset A_{a,1}(B_2)$ $P_R$-presque s{\^u}rement. 
\end{enumerate}
\end{lemm}

\begin{demo}
Montrons les diff{\'e}rents points 

Le premier point d{\'e}coule imm{\'e}diatement des in{\'e}galit{\'e}s 
$$\sup_{n \in \nnf} P[F_a(\varepsilon_1,r) \in B_n] \le
P[F_a(\varepsilon_1,r) \in \bigcup_{n \in \nnf} B_n] \le \sum_{n \in
  \nnf} P[F_a(\varepsilon_1,r) \in B_n].$$ 

Le deuxi{\`e}me point vient du fait que $F_a(\varepsilon_1,R^{(0)})$ a m{\^e}me loi que
$R^{(0)}$ (corollaire~\ref{preservation de la loi du mouvement
  brownien reflechi}), d'o{\`u} 
$$\int_{\wwf_+} P[F_a(\varepsilon_1,r) \in B]\ P_R(dr) =
P[F_a(\varepsilon_1,R^{(0)}) \in B] = P_R(B) = 0.$$
Donc $P[F_a(\varepsilon_1,r) \in B] = 0$ pour $P_R$-presque tout $r \in
\wwf_+$, c'est-{\`a}-dire $P_R(A_{a,1}(B))=0$. 

Le troisi{\`e}me point se d{\'e}duit des deux premiers en remarquant que 
$$A_{a,1}(B_2) \cup A_{a,1}(B_1 \setminus B_2) = A_{a,1}(B_1 \cup B_2)
\supset A_{a,1}(B_1)$$
et que $P_R(A_{a,1}(B_1 \setminus B_2))=0$. \hfill $\square$ 
\end{demo}

\subsection{Ensembles presque ouverts dans $\wwf_+$}~\label{presque-ouvert}

Nous introduisons ici une notion commode pour la suite. 

\begin{defi}~\label{Parties presque ouvertes} 
{\bf (Parties presque ouvertes)}

Soit $V$ une partie de $\wwf_+$. On dit que $V$ est presque ouvert dans 
$\wwf_+$ et que $V$ est un presqu'ouvert de $\wwf_+$ lorsque 
$V \setminus V^\circ$ est n{\'e}gligeable pour $P_R$. De fa{\c c}on {\'e}quivalente, 
un presqu'ouvert de $\wwf_+$ est la r{\'e}union d'un ouvert de $\wwf_+$ 
et d'un n{\'e}gligeable pour $P_R$.  
\end{defi}

Autrement dit, $V$ est presque ouvert dans 
$\wwf_+$ si $P_R$-presque tout point de $V$ est int{\'e}rieur {\`a} $V$. 
Voyons quelques propri{\'e}t{\'e}s des presque-ouverts. 

\begin{lemm}~\label{Proprietes de stabilite} {\bf (Propri{\'e}t{\'e}s 
imm{\'e}diates de stabilit{\'e} des presque-ouverts)}
\begin{itemize}
\item Toute union d{\'e}nombrable de presque-ouverts est un presque-ouvert.
\item Toute intersection finie de presque-ouverts est un presque-ouvert.
\item L'image r{\'e}ciproque d'un ouvert par une application continue 
presque partout de $\wwf_+$ dans $\wwf_+$ est un presque-ouvert. 
\item L'image r{\'e}ciproque d'un presque-ouvert par une application continue 
presque partout de $\wwf_+$ dans $\wwf_+$ et pr{\'e}servant la mesure $P_R$
est un presque-ouvert. 
\end{itemize}
\end{lemm}


\begin{prop}{\bf (Autres propri{\'e}t{\'e}s de stabilit{\'e})}

Soient $V$ un presque-ouvert de $\wwf_+$ et $a>0$. 
Les ensembles ci-dessous sont presque ouverts. 
\begin{itemize} 
\item le translat{\'e} $\theta_{D_a}^{-1}(V)$ ;
\item l'ensemble des 
trajectoires d'o{\`u} l'on peut acc{\'e}der {\`a} $V$ par $F_a$ en un coup : 
$$A_{a,1}(V) = \{r \in \wwf_+ : P[F_a(\varepsilon_1,r) \in V] > 0\} ;$$
\item l'ensemble des trajectoires d'o{\`u} l'on peut acc{\'e}der {\`a} $V$ 
par $F_a$ en $n$ coups, avec $n \in \nnf$ :  
$$A_{a,n}(V) = \{r \in \wwf_+ : P[F_a(\varepsilon_n,\cdot) \circ \cdots \circ
F_a(\varepsilon_1,\cdot) (r) \in V]>0\} ;$$
\item le domaine d'accessibilit{\'e} de $V$ par $F_a$ :
$$A_a(V) = \bigcup_{n \in \nnf} A_{a,n}(V).$$
\end{itemize}   
\end{prop} 

\begin{demo}
Le fait que $\theta_{D_a}^{-1}(V)$ est presque ouvert vient de ce que 
$\theta_{D_a}$ est continue presque partout sur $\wwf_+$ et pr{\'e}serve la 
mesure $P_R$. 

Par ailleurs, $A_{a,1}(V)$ est la r{\'e}union de $A_{a,1}(V^\circ)$, qui 
est presque-ouvert d'apr{\`e}s le lemme~\ref{Consequence de la continuite}
(Cons{\'e}quence de la continuit{\'e} de $F_a$), et de $A_{a,1}(V \setminus V^\circ)$, 
qui est n{\'e}gligeable d'apr{\`e}s le lemme~\ref{Proprietes de $A_{a,1}$} 
(Propri{\'e}t{\'e}s de $A_{a,1}$). Donc $A_{a,1}(V)$ est presque ouvert. 

Pour montrer que pour tout $n \in \nnf$, $A_{a,n}(V)$ est presque ouvert, 
il suffit d'{\'e}tablir la relation de r{\'e}currence 
$A_{a,n}(V) = A_{a,1}(A_{a,n-1}(V))$ pour $n \ge 2$. On {\'e}tablit cette 
{\'e}galit{\'e} en remarquant que pour tout $r \in \wwf_+$, 
$$P[F_a(\varepsilon_n,\cdot) \circ \cdots \circ F_a(\varepsilon_1,\cdot) 
(r) \in V] = \int_E P[F_a(\varepsilon_n,\cdot) \circ \cdots \circ 
F_a(\varepsilon_2,\cdot) (F_a(e,r)) \in V] P_\varepsilon(de),$$ 
d'o{\`u}
$$r \in A_{a,n}(V) \Longleftrightarrow P_\varepsilon \{e \in E : F_a(e,r) 
\in A_{a,n-1}(V)\} > 0 \Longleftrightarrow r \in A_{a,1}(A_{a,n-1}(V)),$$ 
ce qu'il fallait d{\'e}montrer.
\hfill $\square$
\end{demo}

Les presque-ouverts se pr{\^e}tent bien {\`a} l'{\'e}tude de leur accessibilit{\'e}. Mais  
l'int{\'e}r{\^e}t principal des presque-ouverts appara{\^\i}t dans les deux paragraphes 
suivants.  

\subsection{Comparaison des accessibilit{\'e}s pour diff{\'e}rentes valeurs de $a$}

Commen{\c c}ons par d{\'e}montrer le lemme suivant. 

\begin{lemm}~\label{Comparaison des accessibilites en un coup} 
{\bf (Comparaison des accessibilites en un coup)}

Soient $b > a \ge 0$ et $V$ un presque-ouvert de $\wwf_+$. Pour 
$P_R$-presque tout $r \in \wwf_+$,  
$$P[F_b(\varepsilon,r) \in V]>0 \Longrightarrow P[F_a(\varepsilon,r)
\in V]>0.$$ 
Autrement dit, avec les notations du 
lemme~\ref{Consequence de la continuite}, 
$A_{b,1}(V) \subset A_{a,1}(V)$ presque s{\^u}rement.
\end{lemm}

\begin{demo}
Gr{\^a}ce aux propri{\'e}t{\'e}s de $A_{a,1}$ (lemme~\ref{Proprietes de $A_{a,1}$}, 
points 1 et 2), on peut se contenter de montrer le 
r{\'e}sultat dans le cas o{\`u} $V$ est ouvert. 

Fixons une trajectoire $r \in \wwf_+$ sans z{\'e}ro isol{\'e} et sans z{\'e}ro
rationnel telle que $P[F_b(\varepsilon,r) \in V]>0$. Notons
$d_t=D_t(r)$ pour $t \ge 0$. Nous allons d{\'e}montrer que 
$$P[F_a(\varepsilon,r) \in V\ |\ {\bf 1} \mathop{\smile}^{d_b} \varepsilon]>0
\text{ avec probabilit{\'e} strictement positive}.$$ 

Par ind{\'e}pendance des signes, la loi conditionnelle
de $\varepsilon$ sachant $\displaystyle{{\bf 1} \mathop{\smile}^{d_b}
  \varepsilon}$ admet la version r{\'e}guli{\`e}re donn{\'e}e par  
$$\lc(\varepsilon\ |\ {\bf 1} \mathop{\smile}^{d_b} \varepsilon = e) =
\lc(\varepsilon \mathop{\smile}^{d_b} e).$$
pour toute famille $e$ de l'ensemble $E_{d_b} = \{e \in E : \forall q
\le d_b, e(q)=1\}$. 

Par hypoth{\`e}se, on a avec probabilit{\'e} positive 
$$F({\bf 1} \mathop{\smile}^{d_b}
  \varepsilon,r) = F_b(\varepsilon ,r) \in V,$$
et presque s{\^u}rement 
$$\liminf_{t \to d_b+} ({\bf 1} \mathop{\smile}^{d_b} \varepsilon)(Q_t(r)) =
\liminf_{t \to d_b+} \varepsilon(Q_t(r)) = -1.$$
gr{\^a}ce au fait que $r$ poss{\`e}de une infinit{\'e} d'excursions
imm{\'e}diatement apr{\`e}s l'instant $d_b$. 
  
Il suffit donc de montrer que si $e \in E_{d_b}$ v{\'e}rifie $F(e,r) \in
V$ et $\liminf_{t \to d_b+} e(Q_t(r)) = -1$, alors
$$P[F_a(\varepsilon,r) \in V\ |\ {\bf 1} \mathop{\smile}^{d_b}
\varepsilon = e]>0,$$ 
autrement dit
$$P[F_a(\varepsilon \mathop{\smile}^{d_b} e,r) \in V]>0.$$ 

Soient $\rho>0$ et $\delta>0$ tels que $V_\infty(F(e,r),\rho,\delta)
\subset V$, o{\`u} l'on note
$$V_\infty(f,\rho,\delta) = \{g \in V : ||f-g||_{[0,\infty[} < \rho\ ;\ Z(f) 
\subset Z(g) + ]-\delta,\delta[\ ;\ Z(g) \subset Z(f) +
]-\delta,\delta[\}.$$ 
Choisissons un instant $t_0$
r{\'e}alisant le minimum d'une excursion n{\'e}gative de $e \cdot r$ de
hauteur $<\rho/4$ contenue dans $]d_b,d_b +\delta[$. 

Gr{\^a}ce aux hypoth{\`e}ses sur $r$, on peut choisir {\'e}galement des instants
$t_1> \ldots >t_n$ de $]d_a,d_b[$, r{\'e}alisant les maxima d'excursions
de $r$, tels que $r(t_n)< \ldots < r(t_1) < r(t_0)$ et 
$$Z(r) \cap ]d_a,d_b[ \subset \bigcup_{i=1}^n ]t_i - \delta , t_i + \delta[.$$
Notons $q_0>\ldots>q_n$ les rationnels num{\'e}rotant ces excursions. 
Soit $C$ l'ensemble des familles de signes $f \in E$ telles que 
$$\begin{array}{l}
\!\! f(q)=-1 \mbox{ si } q \in \{q_0,\ldots,q_n\},\\ 
\!\! f(q)=1 \mbox{ si } q \notin \{q_0,\ldots,q_n\}
\mbox{ et $q$ num{\'e}rote une excursion de hauteur} 
\ge r(t_n) \mbox{ avant } t_0,\\
\!\! f(q)=1 \mbox{ si } q \notin \{q_0,\ldots,q_n\}
\mbox{ et $q$ num{\'e}rote une excursion de longueur} 
\ge \delta \mbox{ avant } d_b +\delta.
\end{array}$$
L'appartenance {\`a} $C$ ne d{\'e}pend que d'un nombre fini de signes, donc
l'{\'e}v{\'e}nement $\{\varepsilon \in C\}$ est de probabilit{\'e} strictement
positive. Si $\varepsilon \in C$, on peut faire les observations suivantes.
\begin{enumerate}
\item Pour tout rationnel $q$ num{\'e}rotant une excursion de hauteur $\ge
  \rho/4$, 
$$({\bf 1} \mathop{\smile}^{d_a} \varepsilon \mathop{\smile}^{d_b}
e)(q) = e(q).$$
Donc
$$||({\bf 1} \mathop{\smile}^{d_a} \varepsilon \mathop{\smile}^{d_b}
e) \cdot r - e \cdot r||_{[0,+\infty[} < \rho/2$$
d'o{\`u}
$$||F_a(\varepsilon \mathop{\smile}^{d_b} e,r) - F(e,r)||_{[0,+\infty[} =
||F({\bf 1} \mathop{\smile}^{d_a} \varepsilon
\mathop{\smile}^{d_b} e,r) - F(e,r)||_{[0,+\infty[} < \rho.$$   

\item Les trajectoires $F_a(\varepsilon,r)$ et $F(e,r)$ co{\"\i}ncident
  avec $r$ sur $[0,d_a]$ donc ont les m{\^e}mes z{\'e}ros sur $[0,d_a]$. 
Par ailleurs, pour tout $s \ge t_0$, 
$$\displaystyle{ \underline{ ({\bf 1} \mathop{\smile}^{d_a}
  \varepsilon \mathop{\smile}^{d_b} e) \cdot r} (s) = \min_{[d_a,s]}
\big( ({\bf 1} \mathop{\smile}^{d_a} \varepsilon \mathop{\smile}^{d_b}
e) \cdot r \big) = \min_{[t_0,s]} (e \cdot r) = \underline{e \cdot r}(s), }$$
donc $\displaystyle{F_a(\varepsilon
  \mathop{\smile}^{d_b} e,r) = F({\bf 1} \mathop{\smile}^{d_a}
  \varepsilon \mathop{\smile}^{d_b} e,r)}$ et $F(e,r)$ ont aussi les
m{\^e}mes z{\'e}ros sur $[t_0,\infty[$ puisque leurs z{\'e}ros sont 
les instants de records
n{\'e}gatifs de $\displaystyle{({\bf 1} \mathop{\smile}^{d_a} \varepsilon
  \mathop{\smile}^{d_b} e) \cdot r}$ et de $e \cdot r$.

\item Les instants $t_n<\ldots<t_0$ sont des instants de
records n{\'e}gatifs de $\displaystyle{({\bf 1} \mathop{\smile}^{d_a}
  \varepsilon \mathop{\smile}^{d_b} e) \cdot r}$ donc des z{\'e}ros de
$\displaystyle{F_a(\varepsilon \mathop{\smile}^{d_b} e,r) = F({\bf 1}
  \mathop{\smile}^{d_a} \varepsilon \mathop{\smile}^{d_b}
  e,r)}$. Comme $F(e,r)$ co{\"\i}ncide avec $r$ sur $[0,d_b]$ et comme $d_b
< t_0 < d_b + \delta$, 
\begin{eqnarray*}
Z \big( F(e,r) \big)\ \cap\ ]d_a,t_0[\ 
&\subset&  \big( Z(r)\ \cap\ ]d_a,d_b[ \big)\ \cup\ [d_b,t_0[\\ 
&\subset&  \bigcup _{i=0}^n ]t_i - \delta , t_i + \delta[.\\ 
&\subset&  Z \big( F_a(\varepsilon \mathop{\smile}^{d_b} e,r) \big)\ 
+\ ]-\delta , \delta[.
\end{eqnarray*}
Compte tenu du point 2, on a donc
$$Z \big( F(e,r) \big) \subset Z \big( F_a(\varepsilon,r) \big)\
+\ ]-\delta , \delta[.$$

\item Inversement, si $s$ est un z{\'e}ro de
  $\displaystyle{F_a(\varepsilon \mathop{\smile}^{d_b} e,r)}$ dans
  $]d_a,t_0[$, alors :
  \begin{itemize}
     \item soit $d_b \le s < t_0$ et alors $0 \le s-d_b <\delta$ ; 
     \item soit $s < d_b$ ; comme $s$ est dans une excursion n{\'e}gative de
  $\varepsilon \cdot r$, cette excursion est de longueur $< \delta$
  puisque $\varepsilon \in C$.  
   \end{itemize}
Dans tous les cas, $s$ est {\`a} distance $<\delta$ d'un z{\'e}ro de $r$
ant{\'e}rieur {\`a} $d_b$, donc d'un z{\'e}ro de $F(e,r)$, ce qui compte tenu du point 2
montre l'inclusion 
$$Z \big( F_a(\varepsilon \mathop{\smile}^{d_b} e,r) \big) \subset\ Z \big(
F(e,r) \big) +\ ]-\delta , \delta[.$$
\end{enumerate}
Sur l'{\'e}v{\'e}nement de probabilit{\'e} strictement positive
$\{\varepsilon \in C\}$, on a ainsi
$$F_a(\varepsilon \mathop{\smile}^{d_b} e,r) \in
V_\infty(F(e,r),\rho,\delta) \subset V,$$
ce qui ach{\`e}ve la d{\'e}monstration. \hfill $\square$ 
\end{demo}

Signalons un corollaire int{\'e}ressant bien qu'il ne soit pas utilis{\'e} 
dans la suite.

\begin{coro} {\bf (Accessibilit{\'e} en un coup d'un presque-ouvert 
depuis lui-m{\^e}me)}

Sous les hypoth{\`e}ses du lemme~\ref{Consequence de la continuite}, on a
$A_{a,1}(V) \supset V$ presque s{\^u}rement. 
\end{coro}

\begin{demo}
Gr{\^a}ce au lemme~\ref{Proprietes de $A_{a,1}$}, il suffit de v{\'e}rifier
l'inclusion pour un ouvert de la forme $V_t(f,\rho,\delta)$ puisque
$V$ peut s'{\'e}crire comme r{\'e}union d{\'e}nombrable de tels ouverts. Mais si
$V=V_t(f,\rho,\delta)$ avec $t > 0$, $\rho>0$, $\delta>0$, alors $V$
est accessible en un coup par $F_{\max(a,t)}$ depuis tout $w \in V$
puisque $F_{\max(a,t)}$ pr{\'e}serve les trajectoires sur $[0,t]$. Donc
$V$ est accessible en un coup par $F_a$ depuis presque tout $r \in V$. 
\hfill $\square$
\end{demo}

Nous pouvons enfin d{\'e}montrer la proposition ci-dessous, qui sera
extr{\^e}mement utile par la suite. 

\begin{prop}~\label{Comparaison des accessibilites} 
{\bf (Comparaison des accessibilit{\'e}s)}

Soient $b \ge a \ge 0$ et $V$ un presque-ouvert de $\wwf_+$. Pour
$P_R$-presque tout $r \in \wwf_+$, si $V$ est accessible par $F_b$
depuis $r$, alors $V$ est accessible par $F_a$
depuis $r$.
\end{prop}

\begin{demo}
Montrons par r{\'e}currence que pour tout $n \in \nnf$, $A_{b,n}(V) \subset
A_{a,n}(V)$ presque s{\^u}rement. Le r{\'e}sultat est trivial pour $n=0$ et d{\'e}j{\`a} 
{\'e}tabli pour $n=1$ (lemme~\ref{Comparaison des accessibilites en un coup}). 
Une fois l'inclusion $A_{b,n}(V) \subset A_{a,n}(V)$ 
presque s{\^u}rement {\'e}tablie pour un entier $n \in \nnf$, on {\'e}crit 
 $$A_{b,n+1}(V) = A_{b,1}(A_{b,n}(V)) \mathop{\subset}_{\mbox{p.s.}} 
A_{b,1}(A_{a,n}(V)) \mathop{\subset}_{\mbox{p.s.}} A_{a,1}(A_{a,n}(V)) 
= A_{a,n+1}(V),$$
gr{\^a}ce au fait que $A_{b,n}(V)$ et $A_{a,n}(V)$ sont des presque-ouverts et 
gr{\^a}ce au point 3 du lemme~\ref{Proprietes de $A_{a,1}$}.
\hfill $\square$
\end{demo}

\subsection{Bor{\'e}liens stables par $F_a$}

L'int{\'e}r{\^e}t de la transformation $F_a$ est de pr{\'e}server les trajectoires 
jusqu'{\`a} l'instant $D_a$, ce qui nous am{\`e}ne {\`a} nous int{\'e}resser aux 
bor{\'e}liens stables par $F_a$. 

\begin{defi} {\bf (Parties stables par $F_a$)}

Soit $a \ge 0$. On dit qu'une partie $B$ de $\wwf_+$ est stable par $F_a$ 
si pour tout $e \in E$ et $r \in B$, on a $F_a(e,r) \in B$.  
\end{defi}


\begin{lemm}~\label{Exemples de parties stables par $F_a$} 
{\bf (Exemples de parties stables par $F_a$)}

Est stable par $F_a$ :
\begin{itemize}
\item tout bor{\'e}lien ant{\'e}rieur {\`a} $D_a$ dans la filtration naturelle 
canonique de $\wwf_+$ ; 
\item l'ouvert $O_{a,b} = \{w \in \wwf_+ : w \mbox{ poss{\`e}de un z{\'e}ro dans } ]a,b[\}$ pour $b>a$ ; 
\item toute intersection de parties stables par $F_a$. 
\end{itemize}
\end{lemm}

\begin{demo}
Montrons les deux premiers points. 

Remarquons que $D_a$ est un temps d'arr{\^e}t pour la filtration naturelle 
$\fc^0$ associ{\'e}e au processus canonique sur $\wwf_+$. Si $B \in \fc^0_{D_a}$ 
et $r \in B$, alors pour tout $e \in B$, la trajectoire $F_a(e,r)$ co{\"\i}ncide 
avec $r$ jusqu'{\`a} l'instant $D_a(r)$, donc $F_a(e,r) \in B$ gr{\^a}ce au 
crit{\`e}re de Galmarino (voir \cite{Revuz-Yor}, chapitre I, exercice 4.21). 

Montrons que pour $b>a$, $O_{a,b}$ est stable par $F_a$. Soient 
$r \in O_{a,b}$ et $e \in E$. Alors $D_a(r)<b$. De deux choses l'une : 
\begin{itemize}
\item soit la trajectoire $e \cdot r$ reste positive ou nulle sur 
$]D_a(r),b[$. Dans ce cas $e \cdot r = r$ sur $]D_a(r),b[$, d'o{\`u} 
$F_a(e,r)=r$ sur $]D_a(r),b[$, et m{\^e}me sur $[0,b]$. En particulier, 
$F_a(e,r)$ poss{\`e}de un z{\'e}ro dans $]a,b[$.
\item soit la trajectoire $e \cdot r$ prend des valeurs strictement 
n{\'e}gatives sur $]D_a(r),b[$ et alors tout instant de record strictement 
n{\'e}gatif sur $]D_a(r),b[$ est un z{\'e}ro de $F_a(e,r)$. 
\end{itemize}
Dans tous les cas, $F_a(e,r) \in O_{a,b}$, ce qu'il fallait d{\'e}montrer.
\hfill $\square$
\end{demo}

%
%


L'int{\'e}r{\^e}t de la notation de bor{\'e}lien stable par $F_a$ appara{\^\i}t dans la 
proposition ci-dessous. 

\begin{prop}~\label{Accessibilites successives et intersection} 
{\bf (Accessibilit{\'e}s successives et intersection)}

Soient $b > a \ge 0$ et $B_0,B_1,B_2$ des bor{\'e}liens de $\wwf_+$ tels que 
\begin{itemize}
\item $B_1$ est accessible par $F_a$ depuis presque toute trajectoire de $B_0$ ;
\item $B_1$ est stable par $F_b$ ;
\item $B_2$ est accessible par $F_b$ depuis presque toute trajectoire de $B_1$.
\end{itemize}
Alors  $B_1 \cap B_2$ est accessible par $F_b$ depuis presque tout 
$r \in B_1$. Si de plus, $B_1$ et $B_2$ sont presque ouverts, alors 
$B_1 \cap B_2$ est accessible par $F_a$ depuis presque tout $r \in B_0$. 
\end{prop}

\begin{demo}
Soit $r \in B_1$. Comme $B_1$ est stable par $F_b$, on a pour tout $n \in \nnf$,
$F_b(\varepsilon_n,\cdot) \circ \cdots \circ F_b(\varepsilon_1,\cdot) (r) 
\in B_1\ \text{ s{\^u}rement},$ 
d'o{\`u} 
$$P[F_b(\varepsilon_n,\cdot) \circ \cdots \circ F_b(\varepsilon_1,\cdot) (r) 
\in B_1 \cap B_2]
= P[F_b(\varepsilon_n,\cdot) \circ \cdots \circ F_b(\varepsilon_1,\cdot) 
(r) \in B_2].$$
Donc $B_1 \cap B_2$ est accessible par $F_b$ depuis $r$ puisque
$B_2$ l'est, ce qui montre la premi{\`e}re affirmation. 

Si de plus $B_1$ et $B_2$ sont presque ouverts, alors $B_1 \cap B_2$ 
aussi, donc par comparaison des accessibilit{\'e}s
(proposition~\ref{Comparaison des accessibilites}), $B_1 \cap B_2$ 
est accessible par $F_a$ depuis presque tout $r \in B_1$.
Par accessibilit{\'e}s successives
(proposition~\ref{Accessibilites successives}), 
$B_1 \cap B_2$ est accessible par $F_a$ depuis presque tout $r \in B_0$. 
\hfill $\square$
\end{demo}

\subsection{Accessibilit{\'e} et translation}

Nous allons maintenant relier l'accessibilit{\'e} d'un ouvert $V$ par $F$ 
{\`a} l'accessibilit{\'e} d'un translat{\'e} $\theta_{D_a}^{-1}(V)$ par $F_a$. 

\begin{lemm} {\bf (Translation et accessibilit{\'e})}

Soient $a \ge 0$, et $V$ un presque-ouvert de $\wwf_+$. 
Il y a {\'e}quivalence entre  
\begin{enumerate}
\item $V$ est accessible depuis presque toute trajectoire de $\wwf_+$. 
\item $\theta_{D_a}^{-1}(V)$ est accessible par $F_a$ depuis presque 
toute trajectoire de $\wwf_+$.
\end{enumerate}
\end{lemm}

\begin{demo}
Soit $r \in \wwf_+$ telle que $D_a(r)<+\infty$.
Une application r{\'e}p{\'e}t{\'e}e du lemme~\ref{Lien entre $F$ et $F_a$} 
(Lien entre $F$ et $F_a$) montre que pour tout $n \in \nnf$,
$$\theta_{D_a} \circ F_a(\varepsilon_n,\cdot) \circ \cdots \circ 
F_a(\varepsilon_1,\cdot) (r)\ \mbox{ a m{\^e}me loi que }\ 
F(\varepsilon_n,\cdot) 
\circ \cdots \circ F(\varepsilon_1,\cdot) \circ \theta_{D_a} (r),$$ 
d'o{\`u}, 
\begin{eqnarray*}
P[F_a(\varepsilon_n,\cdot) \circ \cdots \circ F_a(\varepsilon_1,\cdot) (r) 
\in \theta_{D_a}^{-1}(V)]
&=&P[\theta_{D_a} \circ F_a(\varepsilon_n,\cdot) \circ \cdots \circ 
F_a(\varepsilon_1,\cdot) (r) \in V]\\
&=&P[F(\varepsilon_n,\cdot) \circ \cdots \circ F(\varepsilon_1,\cdot) 
\circ \theta_{D_a} (r) \in V],
\end{eqnarray*}
Par cons{\'e}quent, $\theta_{D_a}^{-1}(V)$ est accessible par $F_a$ depuis $r$
si et seulement si $V$ est accessible par $F_a$ depuis $\theta_{D_a}(r)$.
Comme $D_a(r)<+\infty$ pour presque tout $r \in \wwf_+$, on a donc
$$A_a(\theta_{D_a}^{-1}(V)) = \theta_{D_a}^{-1}(A_0(V)) \mbox{ p.s.}$$
Comme $\theta_{D_a}$ pr{\'e}serve la loi $P_R$, on a ainsi
$$P_R[A_a(\theta_{D_a}^{-1}(V))] = P_R(A_0(V)),$$ 
ce qui entra{\^\i}ne l'{\'e}quivalence annonc{\'e}e.
\hfill $\square$
\end{demo}

Voyons maintenant un exemple d'application de ce lemme et de la 
proposition~\ref{Accessibilites successives et intersection}, qui anticipe 
sur les r{\'e}sultats qui vont {\^e}tre obtenus dans la section suivante :
en utilisant le fait que pour tout $t>0$, 
$\rho>0$ et $\delta>0$, l'ouvert $V_t(0,\rho,\delta)$
est accessible par $F$ depuis presque toute trajectoire de $\wwf_+$
(d'apr{\`e}s la proposition~\ref{Approximation de zero en topologie CUCZ}), 
on obtient le r{\'e}sultat suivant qui montre qu'on peut \og remettre {\`a} z{\'e}ro\fg une 
trajectoire apr{\`e}s un instant $D_a$ en pr{\'e}servant la trajectoire avant cet 
instant. 

\begin{prop}~\label{Remise a zero apres l'instant $D_a$} 
{\bf (Remise {\`a} z{\'e}ro apr{\`e}s l'instant $D_a$)}

Soient $a \ge 0$, $B$ un  bor{\'e}lien de $\wwf_+$ ant{\'e}rieur {\`a} $D_a$ et
$t>0$, $\rho>0$ et $\delta>0$. 
Alors le bor{\'e}lien $B \cap \theta_{D_a}^{-1}(V_t(0,\rho,\delta))$, 
{\'e}gal {\`a}  
$$\{r \in B : D_a(r) < +\infty\ ;\ ||r||_{[D_a(r),D_a(r)+t]} < \rho \ ;\  
\max_{s \in [D_a(r),D_a(r)+t-\delta]} d(s,Z_r) < \delta)\}$$
est accessible par $F_a$ depuis presque toute trajectoire de $B$. 
\end{prop}


\section{Approximation de la trajectoire nulle}

Le but de cette partie est de montrer le r{\'e}sultat suivant. 

\begin{prop}~\label{Approximation de zero en topologie CUCZ} 
{\bf (Approximation de la trajectoire nulle en topologie CUCZ)}

  Soient $t>0$, $\rho>0$ et $\delta>0$. L'ouvert $V_t(0,\rho,\delta)$
  est accessible depuis presque toute trajectoire de l'ouvert $\wwf_+$. 
\end{prop}

La d{\'e}monstration se fait en deux grandes {\'e}tapes. La premi{\`e}re consiste 
{\`a} approcher z{\'e}ro uniform{\'e}ment sur le segment $[0,t]$, la seconde consiste 
{\`a} densifier les z{\'e}ros sur l'intervalle $[0,t]$. 

\subsection{Approximation uniforme de la trajectoire nulle}

Pour approcher la trajectoire nulle, l'id{\'e}e est la suivante : notons 
${\bf -1} \in E$ la
famille de signes dans laquelle tous les signes valent $-1$. La
transformation $F({\bf -1},\cdot)$, de $\wwf_+$ dans $\wwf_+$, associe {\`a}
toute trajectoire $r$ la trajectoire $(-r)-\underline{(-r)}
= \overline{r}-r$. Nous allons montrer que les images it{\'e}r{\'e}es de toute
trajectoire par cette transformation convergent uniform{\'e}ment sur les
compacts vers la  trajectoire nulle, puis utiliser des arguments de
continuit{\'e}. 

L'outil de la d{\'e}monstration est le comptage du nombre
d'oscillations de hauteur fix{\'e}e.  

\begin{defi} {\bf (Amplitude d'une application d'un intervalle dans $\rrf$)}

Si $f$ est une application d'un intervalle $I$ dans $\rrf$, on appelle
amplitude de $f$ sur $I$ le diam{\`e}tre de $f(I)$, c'est-{\`a}-dire la
diff{\'e}rence $\amp_I\ f = \sup_I f - \inf_I f$.
\end{defi}

\begin{defi}~\label{Nombre d'oscillations} {\bf (Nombre d'oscillations de
  hauteur $h$)}

Soit $h>0$ fix{\'e}. Pour $w \in \wwf$, notons $(T_n(w))$ la suite d'instants
d{\'e}finie par $T_0(w) = 0$ et 
$$T_{n+1}(w) = \inf\{t \ge T_n(w) : \amp_{[T_n(w),t]}\ w \ge h\}.$$ 
Pour tout $t \in \rrf_+$, on appelle nombre
d'oscillations de $w$ hauteur $h$ avant l'instant $t$ l'entier
$N_t(w,h) = \sup\{n \in \nnf : T_n(w) \le t\}$. 
\end{defi}

Remarquons que par continuit{\'e} de $w$, la suite $(T_n(w))$ est
strictement croissante et tend vers $+\infty$, si bien que le nombre
d'oscillations sur un segment $(0,t]$ est fini. 

\begin{lemm} {\bf (Caract{\'e}risation de $N_t(w,h)$)}

Quels que soient $w \in \wwf$ et $h>0$, $t \ge 0$ et $n \in \nnf$, 
$$N_t(w,h) \ge n+1 \Longleftrightarrow \ \sup_{0 = a_0 \le \ldots \le
  a_{n+1} = t}\ \min_{0 \le i \le n}\ \amp_{[a_i,a_{i+1}]}\ w \ge h.$$
\end{lemm}

\begin{demo}
On remarque que l'amplitude de $w$ sur chaque intervalle
$[a_i,a_{i+1}]$ d{\'e}pend contin{\^u}ment de $a_1,\ldots,a_n$, si bien que le
sup pour $0 = a_0 \le \ldots \le a_{n+1} = t$ est en fait un maximum. 

L'implication $\Rightarrow$ est {\'e}vidente : si $N_t(w,h) \ge n+1$, il
suffit de prendre $a_1,\ldots,a_n$ {\'e}gaux {\`a}
$T_1(w),\ldots,T_n(w)$. L'amplitude de $w$ sur chaque intervalle
$[a_i,a_{i+1}]$ est sup{\'e}rieure ou {\'e}gale {\`a} $h$. 

R{\'e}ciproquement, supposons qu'il existe une subdivision $0 = a_0 \le \ldots \le
  a_{n+1} = t$ telle que l'amplitude de $w$ sur chaque intervalle
$[a_i,a_{i+1}]$ soit sup{\'e}rieure ou {\'e}gale {\`a} $h$. Alors une r{\'e}currence
imm{\'e}diate montre que $T_i(w) \le a_i$ pour tout $i \in [0 \ldots
n+1]$. En particulier, $T_{n+1}(w) \le t$, d'o{\`u} $N_t(w,h) \ge n+1$. 
\hfill $\square$
\end{demo}

\begin{coro}\label{semi-continuite} {\bf (Semi-continuit{\'e} 
sup{\'e}rieure de $N_t(w,h)$ par rapport {\`a} $w$)}

Quels que soient $h>0$, $t \ge 0$ et $n \in \nnf$, l'ensemble $\{w \in
\wwf : N_t(w,h) \le n\}$ est un ouvert pour la topologie de la
convergence uniforme sur $[0,t]$.  
\end{coro}

\begin{demo}
D'apr{\`e}s le lemme pr{\'e}c{\'e}dent, pour tout $w \in \wwf$, 
$$N_t(w,h) \le n \Longleftrightarrow \ \sup_{0 = a_0 \le \ldots \le
  a_{n+1} = t}\ \min_{0 \le i \le n}\ \amp_{[a_i,a_{i+1}]}\ w < h.$$
Il suffit de remarquer que pour tout segment $[a,b] \subset [0,t]$
l'application $w \mapsto \amp_{[a,b]}\ w$ est lipschitzienne de
rapport 2, lorsque $\wwf$ est muni de la norme $||\cdot||_{[0,t]}$, et
que cette propri{\'e}t{\'e} est stable par passage aux bornes sup{\'e}rieure et
inf{\'e}rieure. \hfill $\square$ 
\end{demo}

Nous allons maintenant nous int{\'e}resser {\`a} l'effet de la transformation 
$F({\bf -1},\cdot)$ sur le nombre d'oscillations. Commen{\c c}ons par un 
lemme simple. 

\begin{lemm} {\bf (Effet de la transformation $w \mapsto \overline{w}-w$ 
sur l'amplitude)}

Pour tout $w \in \wwf$ et pour tout segment $[a,b] \subset \rrf_+$,
l'amplitude de $\overline{w}-w$ sur $[a,b]$ est major{\'e}e par celle de
$w$. 
\end{lemm}

\begin{demo} Il y a deux cas {\`a} consid{\'e}rer.  

Si $\sup_{[a,b]} w = \overline{w}(b)$, alors pour tout $t \in [a,b]$,    
$$0 \le \overline{w}(t)-w(t) \le\ \overline{w}(b) - \inf_{[a,b]} w =
\sup_{[a,b]} w - \inf_{[a,b]} w,$$ 
d'o{\`u} $\amp_{[a,b]} (\overline{w}-w) \le \amp_{[a,b]} w$. 

Si $\sup_{[a,b]} w < \overline{w}(b)$, alors $\overline{w}$ est
constant sur $[a,b]$, d'o{\`u} $\amp_{[a,b]} (\overline{w}-w) =
\amp_{[a,b]} w$.

Dans les deux cas, $\amp_{[a,b]} (\overline{w}-w) \le \amp_{[a,b]} w$. 
\hfill $\square$
\end{demo} 

Ce lemme montre que le nombre d'oscillations de hauteur
fix{\'e}e de $\overline{w}-w$ avant $t$ est au plus {\'e}gal {\`a} celui de
$w$. Pour les trajectoires positives, on a un r{\'e}sultat meilleur. 

\begin{prop} {\bf (R{\'e}duction du nombre d'oscillations)}

Fixons $h>0$ et reprenons les notations de la d{\'e}finition~\ref{Nombre
  d'oscillations}. Alors pour tout $r \in \wwf_+$ et $n \ge 1$, 
$T_n(\overline{r}-r) \ge T_{n+1}(r)$. Par cons{\'e}quent, pour tout $t \ge 0$, 
$$N_t(\overline{r}-r,h) \le \max(N_t(r,h)-1,0).$$
\end{prop}

\begin{demo}
L'in{\'e}galit{\'e} $T_n(\overline{r}-r) \ge T_{n+1}(r)$ se d{\'e}montre par
r{\'e}currence. 

On commence par remarquer que pour tout $r \in \wwf_+$, et $t \ge 0$,  
$\amp_{[0,t]}\ r = \overline{r}(t)$, si bien que $T_1(r)$ est le temps
d'atteinte du niveau $h$ par $\overline{r}$ et par $r$. De m{\^e}me,
$T_1(\overline{r}-r)$ est le temps d'atteinte du niveau $h$ par
$\overline{r}-r$. 

Pour montrer que $T_1(\overline{r}-r) \ge T_2(r)$, il suffit donc de
v{\'e}rifier que si $0 \le t < T_2(r)$, alors $\overline{r}(t)-r(t) <
h$. Il y a trois cas {\`a} consid{\'e}rer :
\begin{itemize}
\item si $0 \le t < T_1(r)$, alors $0 \le r(t) \le \overline{r}(t) <
  h$~; 
\item si $t = T_1(r)$, alors $r(t) = \overline{r}(t) = h$~; 
\item si $T_1(r) < t < T_2(r)$, alors $\overline{r}(t) =
  \sup_{[T_1(r),t]} r$ et $r(t) \ge \inf_{[T_1(r),t]} r$, d'o{\`u} par
  diff{\'e}rence $\overline{r}(t)-r(t)
\le \amp_{[T_1(r),t]}\ r< h$.
\end{itemize}
Dans tous les cas $\overline{r}(t)-r(t) < h$, ce qui montre que
$T_1(\overline{r}-r) \ge T_2(r)$. 

Soit $n \ge 1$ tel que $T_n(\overline{r}-r) \ge
T_{n+1}(r)$. Montrons que $T_{n+1}(\overline{r}-r) \ge
T_{n+2}(r)$. De deux choses l'une :
\begin{itemize}
\item si $T_n(\overline{r}-r) \ge T_{n+2}(r)$, il n'y a rien {\`a} montrer puisque
$T_{n+1}(\overline{r}-r) \ge T_n(\overline{r}-r)$~; 
\item si $T_n(\overline{r}-r) < T_{n+2}(r)$, on remarque que pour
tout $t$ tel que $T_{n+1}(r) \le t < T_{n+2}(r)$,  
$$\amp_{[T_n(\overline{r}-r),t]}\ (\overline{r}-r) \le
\amp_{[T_n(\overline{r}-r),t]}\ r \le \amp_{[T_{n+1}(r),t]}\
r < h.$$
\end{itemize} 
Dans tous les cas, $T_{n+1}(\overline{r}-r) \ge T_{n+2}(r)$, ce qui
ach{\`e}ve la r{\'e}currence.  \hfill $\square$
\end{demo} 

\begin{coro} {\bf (Des ouverts en cascade)}

Pour tout $h>0$, $t \ge 0$ et $n \in \nnf$, l'ensemble
$V_{t,n}(h) = \{r \in \wwf_+ : N_t(r,h) \le n\}$ est un ouvert de $\wwf_+$ 
pour la topologie de la convergence uniforme sur $[0,t]$ 
et est accessible par F en un coup depuis toute trajectoire de $V_{t,n+1}(h)$. 
\end{coro} 

\begin{demo}
Le fait que $V_{t,n}(h)$ est ouvert d{\'e}coule imm{\'e}diatement de la 
semi-conti\-nuit{\'e} sup{\'e}rieure de $N_t(w,h)$ par rapport {\`a} $w$ 
(corollaire~\ref{semi-continuite}).

Soit $r \in V_{t,n+1}(h)$. D'apr{\`e}s la proposition, $F({\bf -1},r)
\in V_{t,n}(h)$. Mais l'application $F(\cdot,r) : e \mapsto F(e,r) = e
\cdot r - \underline{e \cdot r}$ de $E$ (muni de la topologie produit)
dans $\wwf$ (muni de la topologie de la convergence uniforme sur les
compacts) est continue. L'ensemble des $e \in E$ tels que $F(e,r)
\in V_{t,n}(h)$ est donc un ouvert non vide de $E$, de probabilit{\'e}
strictement positive pour $P_\varepsilon$. \hfill $\square$
\end{demo} 

Nous savons maintenant comment approcher la trajectoire nulle 
uniform{\'e}ment sur un segment $[0,t]$. 

\begin{prop}\label{Approximation uniforme de la trajectoire nulle}
{\bf (Approximation uniforme de la trajectoire nulle)}

Tout voisinage de la trajectoire nulle pour la topologie de la
convergence uniforme sur les compacts est accessible depuis toute
trajectoire de $\wwf_+$. 
\end{prop}

\begin{demo}
Comme les ouverts $(V_{t,n}(h))_{n \in \nnf}$ recouvrent $\wwf_+$, le
corollaire pr{\'e}c{\'e}dent et le lemme~\ref{Accessibilites successives} 
(Accessibilit{\'e}s successives) montrent que $V_{t,0}(h)$ est accessible 
par $F$ depuis toute trajectoire de $\wwf_+$. Mais $V_{t,0}(h)$ n'est 
autre que la trace sur
$\wwf_+$ de la boule de centre $0$ et de rayon $h$ pour la norme
$||\cdot||_{[0,t]}$. D'o{\`u} le r{\'e}sultat. 
\end{demo}

\subsection{Approximation de la trajectoire nulle en topologie CUCZ}

Le but de cette partie est de montrer que tout voisinage de la
trajectoire nulle (pour la topologie CUCZ) est accessible depuis
presque toute trajectoire pour la cha{\^\i}ne de Markov $(R^{(n)})_{n \in
  \nnf}$. En utilisant l'accessibilit{\'e} de tout voisinage de la
trajectoire nulle pour la topologie de la convergence uniforme et 
le lemme~\ref{Accessibilites successives} (Accessibilit{\'e}s successives), 
il nous suffit de montrer le r{\'e}sultat suivant. 

\begin{prop}~\label{Densification des zeros} 
{\bf (Densification des z{\'e}ros)}

  Soient $t>0$, $\rho>0$ et $\delta>0$. L'ouvert $V_t(0,\rho,\delta)$
  est accessible depuis presque toute trajectoire de l'ouvert
  $V_t(0,\rho)$. 
\end{prop}

L'id{\'e}e de la d{\'e}monstration est la suivante : pour $\delta>0$ fix{\'e},
notons $f_\delta$ la transformation de $\wwf_+$ dans $\wwf_+$ d{\'e}finie par 
$$f_\delta(r) = F({\bf 1} \mathop{\smile}^{g(r)} (-{\bf 1}),r),$$
o{\`u} $g(r)$ est le d{\'e}but de la premi{\`e}re excursion de $r$ de longueur
$\ge \delta$ et $\displaystyle{{\bf 1} \mathop{\smile}^{g(r)} (-{\bf 1})}$ 
la famille 
hybride de signes valant $1$ jusqu'{\`a} l'instant $g(r)$ et $-1$
apr{\`e}s (voir d{\'e}finition~\ref{hybridation}). 

La transformation $f_\delta$ pr{\'e}serve le d{\'e}but de la
trajectoire $r$ jusqu'{\`a} l'instant $g(r)$ et agit ensuite comme la
transformation $F({\bf -1},\cdot)$ introduite dans la pr{\'e}c{\'e}dente partie. 
On v{\'e}rifie facilement que pour tout $r \in \wwf_+$ et pour tout $t
\in \rrf_+$, $||f_\delta(r)||_{[0,t]} \le ||r||_{[0,t]}$ et
$g(f_\delta(r)) > g(r)$. Si $r$ ne poss{\`e}de pas d'intervalle de
constance, on montre que la suite d'instants $(g(f_\delta^n(r)))_{n
  \in \nnf}$ tend vers $+\infty$. Les images successives d'une
trajectoire proche de $0$ en norme uniforme sur un segment $[0,t]$
restent donc proches de $0$ pour
l'{\'e}cart $d_t^{CU}$ et finissent donc par {\^e}tre {\`a} distance inf{\'e}rieure {\`a}
$\delta/2$ de $0$ pour l'{\'e}cart $d_t^{CZ}$. 

La d{\'e}monstration comporte toutefois deux difficult{\'e}s suppl{\'e}mentaires
par rapport {\`a} la pr{\'e}c{\'e}dente. D'une part, on ma{\^\i}trise moins facilement
le nombre d'it{\'e}rations n{\'e}cessaires pour que la premi{\`e}re excursion de
longueur $\ge \delta$ commence apr{\`e}s un instant $t$ fix{\'e}. D'autre
part, pour une trajectoire \og typique\fg $r_0 \in \wwf_+$, la trajectoire 
$r_1 = f_\delta(r_0)$ poss{\`e}de de nombreux maxima locaux de m{\^e}me hauteur
(provenant de z{\'e}ros de $r_0$ entre lesquels le maximum courant depuis
$g(r_0)$ n'a pas vari{\'e}), si bien que l'application $F$ 
n'est pas continue en $\displaystyle{({\bf 1} \mathop{\smile}^{g(r_1)}
(-{\bf 1}),r_1)}$. 

C'est pourquoi la d{\'e}monstration que nous proposons passe par la
construction d'une suite $(e_n)_{n \ge 1}$ de familles de signes telle
que les trajectoires $F(e'_n,\cdot) \circ \cdots \circ F(e'_1,\cdot) (r)$
soient proches de $f_\delta^n(r)$ pour toute famille $(e'_n)_{n \ge
  1}$ suffisamment proche de $(e_n)_{n \ge 1}$. 
 
Commen{\c c}ons par {\'e}tablir un lemme simple. Par commodit{\'e}, nous notons 
$F_e = F(e,\cdot)$ l'application  de $\wwf_+$ dans $\wwf_+$ obtenue 
{\`a} partir de $F$ en fixant une famille de signes $e \in E$.

\begin{lemm} {\bf (Obtention d'une partie presque s{\^u}re stable presque 
tous les $F_e$)}

Soit $\Lambda$ une partie de $\wwf_+$ de probabilit{\'e} 1 pour la loi du
mouvement brownien r{\'e}fl{\'e}chi. Alors la partie 
$$\tilde{\Lambda} = \{r \in \wwf_+ : \forall n \in \nnf,\ \ {\rm pour\
  presque\ tout}\ \ e_1,\ldots,e_n \in E, F_{e_n}
\circ \cdots \circ F_{e_1} (r) \in \Lambda\}$$
est encore de probabilit{\'e} 1 pour la loi du mouvement brownien r{\'e}fl{\'e}chi.  
De plus, si $r \in \tilde{\Lambda}$, alors $F_e(r) \in
\tilde{\Lambda}$ pour presque tout $e \in E$
\end{lemm}

\begin{demo}
Soient $R$ un mouvement brownien r{\'e}fl{\'e}chi et $(\varepsilon_n)_{n \ge
  1}$ une suite de variables al{\'e}atoires ind{\'e}pendantes et de loi uniforme
sur $E$, ind{\'e}pendante de $R$. Comme pour tout $n \in \nnf$, la variable
al{\'e}atoire $F_{\varepsilon_n} \circ \cdots \circ F_{\varepsilon_1}(R)$
a m{\^e}me loi que $R$, on a 
$$1 = P[\forall n \in \nnf,\ \ F_{\varepsilon_n}
\circ \cdots \circ F_{\varepsilon_1} (r) \in \Lambda] = \int_{\wwf_+}
P[\forall n \in \nnf,\ \ F_{\varepsilon_n} \circ \cdots \circ
F_{\varepsilon_1} (r) \in \Lambda]\ P_R(dr),$$
d'o{\`u} $P[\forall n \in \nnf,\ \ F_{\varepsilon_n} \circ
\cdots \circ F_{\varepsilon_1} (r) \in \Lambda] = 1$ pour
$P_R$-presque tout $r \in \wwf_+$. \hfill $\square$
\end{demo}

Voyons maintenant la d{\'e}monstration de la 
proposition~\ref{Densification des zeros}. 

\begin{demo}
Nous allons appliquer le lemme pr{\'e}c{\'e}dent {\`a} la partie $\Lambda$ de
$\wwf_+$ form{\'e}e des trajectoires $r$ sans z{\'e}ro isol{\'e}, sans intervalle de
constance et telles que l'application $e \mapsto F_e(r)$ est continue
presque partout sur $E$. 

Soit $r_0 \in \tilde{\Lambda}$ v{\'e}rifiant $\max_{[0,t]} r_0 < \rho$. 
Nous allons montrer que $V_t(0,\rho,\delta)$ est accessible depuis $r_0$.
La d{\'e}monstration comporte trois {\'e}tapes. 

\paragraph{Premi{\`e}re {\'e}tape :} on construit par r{\'e}currence une suite
$(e_n)_{n \ge 1}$ bien choisie. 

Notons $\eta = \rho - \max_{[0,t]} r_0 > 0$. Notons $g_0$ le
d{\'e}but de la premi{\`e}re excursion de longueur $\ge \delta$. Soient $0 =
z_0 < \ldots < z_d = g_0$ des z{\'e}ros de $r_0$ espac{\'e}s de moins de
$\delta$. Posons $t_0=0$. Soit $t_{d+1} \in ]g_0,g_0+\delta[$ un
instant tel que $r_0(t_{d+1}) = \max_{[g_0,t_{d+1}]} r_0$.  

Comme $r_0 \in \Lambda$, on peut choisir des petites excursions de
$r_0$ dont les maxima soient r{\'e}alis{\'e}s en des instants $t_1 < \ldots <
t_d$ proches de $z_1 < \ldots < z_d$, espac{\'e}s de moins de $\delta$,
tels que $t_d+\delta>z_d=g_0$ et dont les hauteurs v{\'e}rifient  
$$r_0(t_1) < \ldots < r_0(t_d) <
  \min(\eta,r_0(t_{d+1})).$$

\begin{figure}[hbtp]\label{densification}
  \begin{center}
    \leavevmode \input{densification.pstex_t} 
  \end{center}
  Figure~\ref{densification}. --- Passage de $r_0$ \`a $r_1=F(e_1,r_0)$. 
\end{figure}

Soit $C(r_0) \subset E$ l'ensemble des familles de signes
telles que l'action $e \mapsto e \cdot r_0$ affecte  
\begin{itemize}
\item du signe + les excursions de hauteur $\ge \eta$ ant{\'e}rieures {\`a}
  $g_0$ ; 
\item du signe -- les excursions de hauteur $\ge \eta$ post{\'e}rieures {\`a}
  $g_0$ ; 
\item du signe -- l'excursion commen{\c c}ant {\`a} $g_0$ ; 
\item du signe -- les excursions enjambant $t_1,\ldots,t_{d+1}$ ; 
\item du signe + les excursions ant{\'e}rieures {\`a} un instant $t_i$ et de
  hauteur $> r(t_i)$. 
\end{itemize}
La partie $C(r_0)$ est un ouvert de $E$ de probabilit{\'e} strictement positive. 
On peut donc choisir $e_1 \in C(r_0)$ tel que $F$ est continue en
$(e_1,r_0)$ et $r_1 := F_{e_1}(r_0) = e_1 \cdot r_0 - \underline{e_1
  \cdot r_0} \in \tilde{\Lambda}$. 

Par construction, les instants $t_1 < \ldots < t_{d+1}$ sont des records
n{\'e}gatifs de $e_1 \cdot r_0$  et donc des z{\'e}ros de $r_1$ ; par cons{\'e}quent
la premi{\`e}re excursion de $r_1$ de longueur $\ge \delta$ d{\'e}bute {\`a} un
instant $g_1>g_0$. Par ailleurs, $r_1$ est major{\'e}e strictement par
$\rho$ sur $[0,t]$. En effet,  
\begin{itemize}
\item si $s \in [0,g_0]$, $(e_1 \cdot r_0)(s) \le \max r_0([0,t])$ et
  $-\underline{e_1 \cdot r_0}(s) < \eta$ ;
\item si $s \in [g_0,t]$, $(e_1 \cdot r_0)(s) < \eta$ et $-\underline{e_1
    \cdot r_0}(s) \le  \max r_0([0,t])$. 
\end{itemize}
Dans tous les cas, $r_1(s) < \max r_0([0,t]) + \eta = \rho$.

Comme $r_1$ v{\'e}rifie les m{\^e}mes hypoth{\`e}ses que $f$, on peut donc
continuer la construction en choisissant $e_2 \in C(r_1)$ tel que $F$
est continue en $(e_2,r_1)$ et $r_2 := F_{e_2}(r_1) \in
\tilde{\Lambda}$, puis $e_3 \in C(r_2)$ tel que $F$ est continue en
$(e_3,r_2)$ et $r_3 := F_{e_3}(r_2) \in \tilde{\Lambda}$, etc... 

\paragraph{Deuxi{\`e}me {\'e}tape :} on montre que le d{\'e}but de la premi{\`e}re
excursion de longueur $\ge \delta$ des trajectoires $r_n = F_{e_n}
\circ \cdots \circ F_{e_1} (r)$ tend vers $+\infty$. 

Notons $g_n$ le d{\'e}but de la premi{\`e}re excursion de $r_n$ de longueur
$\ge \delta$. La suite croissante $(g_n)_{n \in \nnf}$ poss{\`e}de une
limite $g_\infty$. Montrons que $g_\infty = +\infty$. 

On raisonne par l'absurde en supposant que $g_\infty <
+\infty$. Fixons alors $m \in \nnf$ tel que $g_m > g_\infty -
\delta$. Pour tout $n \ge m$, $e_{n+1} \cdot r_n$ est de signe
constant sur $[g_n,g_n+\delta]$ et {\it a fortiori} sur le
sous-intervalle $[g_{n+1},g_m+\delta]$, n{\'e}gatif puisque $g_{n+1}$ est
un instant de record n{\'e}gatif de $e_{n+1} \cdot r_n$. Ce record dure au
moins une dur{\'e}e $\delta$ puisque $r_{n+1}$ ne s'annule pas entre $g_{n+1}$ et
$g_{n+1}+\delta$. Donc pour tout $t \in [g_{n+1},g_m+\delta]$, 
$$r_{n+1}(t) = (e_{n+1} \cdot r_n)(t) - \underline{(e_{n+1} \cdot
  r_n)}(t) = - r_n(t) + r_n(g_{n+1}).$$ 
Cette {\'e}galit{\'e} a deux cons{\'e}quences : pour tout $n \ge m$,
\begin{enumerate}
\item l'application $r_{n+1} + r_n$ est constante
  sur $[g_{n+1},g_m+\delta]$ ; 
\item pour tout $t \in [g_{n+1},g_m+\delta]$, $r_n(t) \le
  r_n(g_{n+1})$.  
\end{enumerate}
On en d{\'e}duit par r{\'e}currence que pour tout $n \ge m$ et $t \in
[g_\infty,g_m+\delta]$, 
\begin{eqnarray*}
r_m(t) \ge r_m(g_n) &\text { si } n-m \text { est pair }\ \ \ \\
r_m(t) \le r_m(g_n) &\text { si } n-m \text { est impair}.
\end{eqnarray*}
En faisant tendre $n$ vers l'infini, on voit que cela entra{\^\i}ne que
$r_m$ est constante ({\'e}gale {\`a} $r_m(g_\infty)$) sur
$[g_\infty,g_m+\delta]$, ce qui contredit le fait que $r_m \in
\tilde{\Lambda}$.  

\paragraph{Troisi{\`e}me {\'e}tape :} on utilise la continuit{\'e} de $F$ aux
points $(e_n,r_{n-1})$. 

Les deux premi{\`e}res {\'e}tapes montrent que $r_n \in V_t(0,\rho,\delta)$ {\`a}
partir d'un certain rang. Mais pour tout $n \in \nnf$
l'application $G_n : (e'_n,\ldots,e'_1,r) \mapsto F_{e'_n} \circ
\cdots \circ F_{e'_1} (r)$ de $E^n \times
\wwf_+$ est continue au point $(e_n,\ldots,e_1,r_0)$. Cela se montre
par r{\'e}currence en remarquant que c'est {\'e}vident pour $n=0$ et que pour
tout $n \ge 1$, $G_n$ est la compos{\'e}e de 
$$(e'_n,\ldots,e'_1,r) \mapsto (e'_n,G_{n-1}(e'_{n-1},\ldots,e'_1,r))$$
et de $F$ qui continue en $(e_n,r_{n-1})$. Comme $r_n \in
V_t(0,\rho,\delta)$ pour $n$ suffisamment grand, $F_{e'_n}
\circ \cdots \circ F_{e'_1} (r_0) \in V_t(0,\rho,\delta)$ pour tout
$(e'_1,\ldots,e'_n)$ dans un certain voisinage de $(e_n,\ldots,e_1)$,
ce qui montre que $V_t(0,\rho,\delta)$ est accessible depuis $r_0$.
\hfill $\square$ 
\end{demo}

\section{Construction d'excursions de hauteur prescrite 
et lemme du v{\'e}rin}

Dans cette partie nous montrons comment, {\`a} partir d'une 
trajectoire proche la trajectoire nulle en topologie CUCZ, 
construire une trajectoire ayant une excursion de hauteur voulue 
localis{\'e}e pr{\`e}s d'un instant donn{\'e}. De telles excursions sont 
l'outil permettant de soulever la trajectoire brownienne 
d'une hauteur donn{\'e}e entre deux instants donn{\'e}s, ce qui est 
l'objet du lemme du v{\'e}rin. 
Dans toute la suite, nous noterons pour $b > a \ge 0$, 
$$O_{a,b} = \{r \in \wwf_+ : \exists t \in ]a,b[, r(t)=0\}$$ 
l'ensemble des trajectoires de $\wwf_+$ poss{\'e}dant au moins un 
z{\'e}ro dans $]a,b[$. 
Cet ensemble est un ouvert de $\wwf_+$ d'apr{\`e}s le 
lemme~\ref{Exemple d'ouvert de la topologie CUCZ}, stable par $F_a$ 
d'apr{\`e}s le lemme~\ref{Exemples de parties stables par $F_a$}. 

\subsection{Pr{\'e}liminaires}

\begin{lemm}~\label{Continuite des parametres d'une excursion} 
{\bf (Continuit{\'e} des param{\`e}tres d'une excursion)}

Pour $t>0$ et $r \in \wwf_+$, notons $G_t(r)$, $D_t(r)$, $H_t(r)$ 
l'extr{\'e}mit{\'e} gauche, l'extr{\'e}mit{\'e} droite et la hauteur 
de l'excursion de $r$ enjambant $t$ :  
$$G_t(r) = \sup(Z(r) \cap [0,t]),\ D_t(r) = \inf(Z(r) \cap ]t,+\infty[),$$ 
$$H_t(r) = \max\{r(s)\ ;\ s \in [G_t(r),D_t(r)]\}.$$ 
Les fonctionnelles $G_t$, $D_t$ et $H_t$ sont continues en toute trajectoire 
$r_0 \in \wwf_+$ telle que $r_0(t)>0$. 

Plus g{\'e}n{\'e}ralement, si $b \ge a \ge 0$, les fonctionnelles qui 
{\`a} une trajectoire $r \in \wwf_+$ associent son maximum sur 
les intervalles 
$[G_a(r),G_b(r)]$, $[G_a(r),b]$, $[G_a(r),D_b(r)]$, $[a,G_b(r)]$, etc...   
(avec la convention $[c,d] = [d,c]$ si $c>d$) sont continues en toute 
trajectoire $r_0 \in \wwf_+$ telle que $r_0(a)>0$ et $r_0(b)>0$. 
\end{lemm}

\begin{demo} 
Soit $r_0 \in \wwf_+$ tels que $r_0(t)>0$. Notons $g=G_t(r_0)$, $d=D_t(r_0)$. 
Par hypoth{\`e}se $g<t<d$. Fixons $\delta>0$ tel que $\delta < \min(t-g,d-t)$. 
Pour tout $r \in \wwf_+$ tel que $d_{d+\delta}^{CZ}(r_0,r)<\delta$,  
$G_t(r) \in ]g-\delta,g+\delta[$ et $D_t(r) \in ]d-\delta,d+\delta[$, 
ce qui montre la continuit{\'e} de $G_t$ et $D_t$. 

On en d{\'e}duit les autres points par continuit{\'e} de l'application 
$(r,s,t) \mapsto \max_{[s,t]} r$
de $\wwf_+ \times \rrf_+^2$ dans $\rrf_+$. 
\hfill $\square$ 
\end{demo}

Le lemme ci-dessous reprend un lemme similaire de l'article de 
Malric~\cite{Malric 3, Malric 4}. 

\begin{lemm}~\label{Somme des hauteurs des excursions} 
{\bf (Somme des hauteurs des excursions)}

Soient $b > a \ge 0$ fix{\'e}s. Pour presque tout $r \in O_{a,b}$, la
somme des hauteurs des excursions de $r$ commenc{\'e}es apr{\`e}s 
$a$ et achev{\'e}es avant $b$ est infinie.
\end{lemm}

\begin{demo}
Notons $N_h(r)$ le nombre d'excursions de $r$ de hauteur $>h$ pendant 
l'intervalle $[a,b]$ et $(H_n(r))_{n \ge 1}$ la suite des hauteurs rang{\'e}es
par ordre d{\'e}croissant. Notons $L_t(r)$ le temps local de $r$ en $0$ {\`a} 
l'instant $t$. Alors pour presque tout $r \in O_{a,b}$, 
$$h N_h(r) \to L_b(r)-L_a(r) \mbox{ quand } h \to 0 \mbox{ et } 
L_b(r)-L_a(r) \in \rrf_+^*,$$ 
Par cons{\'e}quent, pour presque tout $r \in O_{a,b}$, $H_n(r) \to 0$ quand 
$n \to +\infty$ et
$$\sum_{n \ge 1} H_n(r) = \int_0^\infty N_h(r)\ dh = +\infty.$$
ce qui montre le r{\'e}sultat annonc{\'e}.
\hfill $\square$
\end{demo}

\subsection{Obtention d'une trajectoire dont le maximum sur un segment 
d{\'e}passe une hauteur fix{\'e}e.}

Dans ce paragraphe et le suivant, nous montrons deux lemmes faisant 
intervenir les excursions. La preuve de ces lemmes s'inspire de 
la d{\'e}monstration de la densit{\'e} {\`a} un temps par 
Malric~\cite{Malric 2, Malric 3, Malric 4} en la 
simplifiant {\`a} l'aide des r{\'e}sultats vus sur l'accessibilit{\'e} 
des ouverts. 

\begin{lemm}~\label{Maximum > $h$} {\bf (Obtention d'un maximum d{\'e}passant $h$)}

Soient $b > a \ge 0$ et $h>0$. Alors l'ouvert 
$\{r \in O_{a,b} : \max_{[D_a(r),b]} r > h\}$ est accessible par $F_a$ 
depuis presque tout $r \in O_{a,b}$. 
\end{lemm}

\begin{demo}
Pour $n \ge 1$ et $t_1<\ldots<t_n$ rationnels de $]a,b[$, notons 
$V_{t_1,\ldots,t_n}$ l'ensemble des trajectoires $r \in O_{a,b}$ ne 
s'annulant pas en $t_1,\ldots,t_n$ et telles que  
$$a < G_{t_1}(r) < D_{t_1}(r) < ... < G_{t_n}(r) < b$$
et
$$\sum_{k=1}^{n-1} H_{t_k}(r) + \max_{[G_{t_n}(r),D_{t_n}(r) \wedge b]}r > h.$$
Notons $U_n$ la r{\'e}union des $V_{t_1,\ldots,t_n}$ pour $t_1<\ldots<t_n$ 
rationnels de $]a,b[$. Autrement dit, $U_n$ est l'ensemble des trajectoires 
poss{\'e}dant $n$ excursions commen{\c c}ant dans $]a,b[$ dont les intervalles ferm{\'e}s 
d'excursion sont disjoints et dont la somme des hauteurs d{\'e}passe $h$, en ne 
comptant que la hauteur maximale avant $b$ si la derni{\`e}re excursion se 
termine apr{\`e}s l'instant $b$. En particulier, 
$U_1 = \{r \in O_{a,b} : \max_{[D_a(r),b]} r > h\}$.

D'apr{\`e}s le lemme~\ref{Continuite des parametres d'une excursion} , 
tous ces ensembles sont des presque-ouverts de $\wwf_+$. D'apr{\`e}s le 
lemme~\ref{Somme des hauteurs des excursions} et le fait que presque 
s{\^u}rement, les z{\'e}ros ne sont pas isol{\'e}s, $O_{a,b}$ est presque 
s{\^u}rement {\'e}gal {\`a} la r{\'e}union des $U_n$ pour $n \ge 1$. En effet, pour 
presque tout $r \in O_{a,b}$, on peut choisir un nombre fini 
d'excursions de $r$ compl{\`e}tement r{\'e}alis{\'e}es dans l'intervalle $]a,b[$ 
et dont la somme des hauteurs d{\'e}passe $h$. 

Il suffit donc de montrer que pour tout $n \ge 2$, $U_{n-1}$ est accessible 
par $F_a$ depuis toute trajectoire $r \in U_n$. Pour ce faire, nous allons 
d{\'e}montrer que si $t_1<\ldots<t_n$ sont des rationnels de $]a,b[$, alors 
$U_{n-1}$ est accessible par $F_{t_{n-2}}$ (avec la convention $t_0=a$ si 
$n=2$) depuis toute trajectoire de $V_{t_1,\ldots,t_n}$.

Soient $r \in V_{t_1,\ldots,t_n}$, $e \in E$ et 
$\widetilde{r} = F_{t_{n-2}}(e,r)$. Alors $\widetilde{r} \in V_{n-1}$ d{\`e}s que 
l'action de $e$ sur $r$ affecte du signe $+$ l'excursion enjambant $t_n$ et du 
signe $-$ l'excursion enjambant $t_{n-1}$. En effet, $\widetilde{r}$ 
a les m{\^e}mes excursions que $r$ avant l'instant $D_{t_{n-2}}$, tandis que 
pour tout $s \in [G_{t_n}(r),D_{t_n}(r)]$,  
$$\widetilde{r}(s) = r(s) - \min_{[t_{n-2},s]} (e \cdot r) \ge r(s) + 
H_{t_{n-1}}(r) > 0,$$
donc $[G_{t_n}(\widetilde{r}),D_{t_n}(\widetilde{r})] \supset 
[G_{t_n}(r),D_{t_n}(r)]$ et 
$$\max_{[G_{t_n}(\widetilde{r}),D_{t_n}(\widetilde{r}) \wedge b]} \widetilde{r} 
\ge \max_{[G_{t_n}(r),D_{t_n}(r) \wedge b]} \widetilde{r} 
\ge \max_{[G_{t_n}(r),D_{t_n}(r) \wedge b]} r + H_{t_{n-1}}(r).$$
Ainsi, $\widetilde{r} \in V_{t_1,\ldots,t_{n-2},t_n}$.
\hfill $\square$ 
\end{demo}

\subsection{Obtention d'une trajectoire dont la plus grande excursion sur un segment 
approche une hauteur fix{\'e}e}

\begin{lemm}~\label{Obtention d'une plus grande excursion de hauteur voulue} 
({\bf Obtention d'une plus grande excursion de hauteur voulue})

Soient deux instants $b > a \ge 0$, une hauteur $h>0$, une pr{\'e}cision 
$\Delta h \in ]0,h/2[$ et un r{\'e}el $\delta > 0$. Alors l'ensemble 
des trajectoires de $O_{b,b+\delta}$ telles que 
$$\max_{[G_b(r),D_b(r)]} r < \max_{[D_a(r),G_b(r)]} r\ \in\ ]h-\Delta h,h[$$
est accessible par $F_a$ depuis presque tout $r \in O_{a,b}$. 

Autrement dit, de presque toute trajectoire de $O_{a,b}$, on peut 
acc{\'e}der par $F_a$ {\`a} une trajectoire poss{\'e}dant un z{\'e}ro entre 
$b$ et $b+\delta$ et, pendant l'intervalle $]a,b[$, une excursion compl{\`e}te 
de hauteur dans $]h-\Delta h,h[$ qui r{\'e}alise 
le maximum sur l'intervalle $[D_a,D_b]$. 
\end{lemm}

\begin{demo}
La d{\'e}monstration de l'accessibilit{\'e} se d{\'e}compose en plusieurs {\'e}tapes. 

\paragraph{Premi{\`e}re {\'e}tape :} 
Obtention d'une trajectoire dont le maximum sur $[D_a,b]$ est $<h$, 
et existence d'un nombre fini d'excursions dont la derni{\`e}re r{\'e}alise 
un record et dont la somme des hauteurs approche $h$ par d{\'e}faut. 

Par remise {\`a} z{\'e}ro apr{\`e}s l'instant $D_a$, 
l'ensemble des trajectoires $r \in O_{a,b}$ telles que 
$\max_{[D_a(r),b]} r < \Delta h$ est accessible par $F_a$ depuis presque tout  
$r \in O_{a,b}$. 

Pour $n \ge 1$ et $t_1<\ldots<t_n$ rationnels de $]a,b[$, notons 
$V_{t_1,\ldots,t_n}$ l'ensemble des trajectoires $r \in O_{a,b}$ ne 
s'annulant pas en $t_1,\ldots,t_n$ et telles que
$$a < G_{t_1}(r) < D_{t_1}(r) < ... < G_{t_n}(r) < D_{t_n}(r) < b,$$
$$\sum_{k=1}^{n} H_{t_k}(r)\ \in\ ]h-\Delta h,h[,$$
$$H_{t_n}(r) > \max_{[D_a(r),G_{t_n}(r)]} r.$$

Alors l'ensemble de trajectoires $r \in O_{a,b}$ telles que 
$\max_{[D_a(r),b]} r < \Delta h$ est presque s{\^u}rement inclus dans  
la r{\'e}union des $V_{t_1,\ldots,t_n}$ pour $n \ge 1$ et $t_1<\ldots<t_n$ 
rationnels de $]a,b[$. En effet, si $r \in O_{a,b}$ et 
$\max_{[D_a(r),b]} r < \Delta h$, alors pour $n$ bien choisi, la somme 
des hauteurs de la plus haute excursion de $r$ sur $]a,b[$ et des $n-1$
plus hautes excursions de $r$ sur $]a,b[$ qui la pr{\'e}c{\`e}dent  
est dans $]h-\Delta h,h[$, gr{\^a}ce au 
lemme~\ref{Somme des hauteurs des excursions} et au fait que les 
hauteurs des excursions sont $< \Delta h$. 

\paragraph{Deuxi{\`e}me {\'e}tape :} Pour $n \ge 2$ et $t_1<\ldots<t_n$ rationnels 
de $]a,b[$, $V_{t_1,\ldots,t_{n-2},t_n}$ est accessible par $F_{t_{n-2}}$ 
depuis presque toute trajectoire de $V_{t_1,\ldots,t_n}$, avec la convention 
$t_0=a$. 
Cette {\'e}tape est illustr{\'e}e par la figure~\ref{Excursion de hauteur prescrite}

Comme $V_{t_1,\ldots,t_n}$ est ant{\'e}rieur {\`a} $D_{t_n}$, le 
lemme~\ref{Maximum > $h$} appliqu{\'e} {\`a} l'intervalle $[t_n,b]$ et {\`a} la hauteur 
$H_{t_{n-1}(r)}$ montre que le presque-ouvert 
$\{r \in V_{t_1,\ldots,t_n} : \max_{[D_{t_n}(r),b]} r > H_{t_{n-1}(r)}\}$ 
est accessible par $F_{t_n}$ depuis presque toute trajectoire de 
$V_{t_1,\ldots,t_n}$. 

Il suffit donc de montrer que $V_{t_1,\ldots,t_{n-2},t_n}$ 
est accessible par $F_{t_{n-2}}$ depuis presque toute trajectoire de 
$V_{t_1,\ldots,t_n}$ telle que $\max_{[D_{t_n}(r),b]} r > H_{t_{n-1}(r)}$.

Soient $r \in V_{t_1,\ldots,t_n}$ telle que $\max_{[D_{t_n}(r),b]} r > 
H_{t_{n-1}(r)}$, $e \in E$ et $\widetilde{r} = F_{t_{n-2}}(e,r)$. 
Supposons que l'action de $e$ sur $r$ affecte 
\begin{itemize}
\item du signe $-$ l'excursion enjambant $t_{n-1}$ ; 
\item du signe $+$ l'excursion enjambant $t_n$ ;
\item du signe $-$ la premi{\`e}re excursion de hauteur 
$\ge H_{t_{n-1}(r)}$ apr{\`e}s l'instant $D_{t_n}(r)$.
\item du signe $+$ les autres excursions de hauteur 
$\ge H_{t_{n-1}}(r)$ entre $D_{t_{n-2}}(r)$ et $G_{t_n}(r)$.  
\end{itemize}
Nous allons montrer que $\widetilde{r} \in V_{t_1,\ldots,t_{n-2},t_n}$.

Notons $g$ le dernier instant r{\'e}alisant le maximum de l'excursion 
enjambant $t_{n-1}$ et $d = \min\{t \ge D_{t_n}(r) : r(t) \ge H_{t_{n-1}}(r)\}$. 
Alors $t_{n-2} < g < G_{t_n}(r) < D_{t_n}(r) < d < b$. Pour tout 
$t \in [D_{t_{n-2}}(r),d]$, $(e \cdot r)_t \ge - H_{t_{n-1}}(r)$,
avec {\'e}galit{\'e} si et seulement si $t \in \{g,d\}$ et in{\'e}galit{\'e} stricte si 
$t \in ]g,d[$. 
Donc $[g,d]$ est un intervalle d'excursion de $\widetilde{r}$. 

De plus, pour tout $t \in [g,d]$, 
$\widetilde{r}(t) = (e \cdot r)(t) + H_{t_{n-1}}(r)$, 
et le maximum de $e \cdot r$ sur $[g,d]$ est atteint sur 
$[G_{t_n}(r),D_{t_n}(r)]$ et vaut $H_{t_n}(r)$. Donc 
$H_{t_n}(\widetilde{r}) = H_{t_n}(r) + H_{t_{n-1}}(r)$. 

Comme $\widetilde{r}$ co{\"\i}ncide avec $r$ jusqu'{\`a} l'instant $D_{t_{n-2}}$ 
et comme pour tout $t \in [D_{t_{n-2}}(r),g]$, 
$\widetilde{r}(t) \le (e \cdot r)_t + H_{t_{n-1}}(r) 
< H_{t_n}(r) + H_{t_{n-1}}(r)$, on a ainsi 
$\widetilde{r} \in V_{t_1,\ldots,t_{n-2},t_n}$.

Les deux premi{\`e}res {\'e}tapes montrent que la r{\'e}union des 
$V_t$ pour $t$ 
rationnel de $]a,b[$ est accessible par $F_a$ depuis presque toute 
trajectoire de $O_{a,b}$. Une fois qu'on a obtenu une trajectoire 
appartenant {\`a} l'un des $V_t$, on souhaite obtenir une trajectoire 
qui en outre s'annule entre $b$ et $b+\delta$ et telle que l'excursion 
enjambant $t$ r{\'e}alise le maximum sur $[D_a,D_b]$. Nous allons obtenir 
ce r{\'e}sultat 
par remise {\`a} z{\'e}ro apr{\`e}s l'instant $D_t$. 

\begin{figure}[hbtp]
\label{Excursion de hauteur prescrite}
  \begin{center}
    \leavevmode \input{passage.pstex_t} 
  \end{center}
  Figure~\ref{Excursion de hauteur prescrite}. --- 
Passage de $V_{t_1,\ldots,t_n}$ {\`a} $V_{t_1,\ldots,t_{n-2},t_n}$ par $F_{t_{n-2}}$. 
\end{figure}

\paragraph{Troisi{\`e}me {\'e}tape :} Pour tout rationnel $t \in ]a,b[$, l'ensemble 
des trajectoires de $O_{b,b+\delta}$ 
telles que $\max_{[G_b(r),D_b(r)]} r < \max_{[D_a(r),G_b(r)]} r 
\in ]h-\Delta h,h[$ est accessible par $F_a$ depuis presque toute 
trajectoire de $V_t$. 

Soient $t \in ]a,b[$ et $r \in V_t$. Alors l'excursion enjambant $t$ est 
enti{\`e}rement contenue dans l'intervalle $]a,b[$, de hauteur 
$H_t(r) \in ]h-\Delta h,h[$ v{\'e}rifiant $H_t(r) > \max_{[D_a(r),G_t(r)]} r$. 

Par remise {\`a} z{\'e}ro apr{\`e}s l'instant 
$D_t$, le presque-ouvert 
$V_t \cap \theta_{D_t}^{-1}(V_{b-t+\delta}(0,H_t(r),\delta))$ 
est accessible par $F_t$ donc par $F_a$ depuis $r$. Les trajectoires de cet 
ensemble v{\'e}rifient les in{\'e}galit{\'e}s voulues. 
\hfill $\square$
\end{demo}

\subsection{Lemme du v{\'e}rin}

Dans la suite, on s'int{\'e}resse {\`a} l'accessibilit{\'e} d'ouverts de la forme 
$$U_{a,b}(f,\rho,\delta) = \{r \in O_{a,a+\delta} \cap O_{b-\delta,b} : 
||r-f||_{[a,b]} < \rho\}$$ 
pour $b > a \ge 0$, $f \in \wwf_+$ telle que $f(a)=f(b)=0$, $\rho>0$ 
et $0 < \delta < b-a$. 

Le lemme du v{\'e}rin montre comment soulever 
une trajectoire brownienne d'une hauteur $h$ sur un intervalle fix{\'e}.  

\begin{lemm}~\label{Lemme du verin} {\bf (Lemme du v{\'e}rin)}

Consid{\'e}rons des instants $b > b' > a' > a \ge 0$, une hauteur $h>0$ et 
une trajectoire $g \in \wwf_+$, nulle hors de $[a',b']$. D{\'e}finissons une 
trajectoire $f \in \wwf_+$ par 
$$f(t) = \left| \begin{array}{ll} 
\displaystyle{h \times \frac{t-a}{a'-a}} & \mbox{ si } a \le t \le a'\\
\displaystyle{ g(t) + h }               & \mbox{ si } a' \le t \le b'\\
\displaystyle{h \times \frac{b-t}{b-b'}} & \mbox{ si } b' \le t \le b\\
\displaystyle{ 0 }                      & \mbox{ si } t \notin [a,b].
\end{array} \right.$$
Soient $\rho>\rho'>0$ et $\delta>0$ v{\'e}rifiant 
$\osc(f,\delta)<\rho'/4$ et $\delta < \min(a-a',b'-a',b-b')$.  

{\parindent 0cm 
Si l'ouvert $U_{a',b'}(g,\rho',\delta)$ est accessible par $F_{a'}$ depuis 
presque tout $r \in O_{a',a'+\delta}$ tel que $||r||_{[a',D_{a'}(r)]}<\rho'$, 
alors l'ouvert $U_{a,b}(f,\rho,\delta)$ est accessible par $F_a$ depuis 
presque tout $r \in O_{a,a+\delta}$ tel que $||r||_{[a,D_a(r)]}<\rho$.}
\end{lemm}

\begin{demo}
On commence par remarquer que gr{\^a}ce aux hyptoh{\`e}ses faites sur $\delta$, 
$$a < a+\delta < a' < b'-\delta < b' < b-\delta < b.$$ 

\paragraph{1.} D{\'e}monstration dans le cas o{\`u} $h<\rho/2$. 

Fixons $\Delta h \in ]0,\rho-\rho'[$ tel que $\Delta h < h/2$ et notons
\begin{eqnarray*}
V_0&=&\{ r \in O_{a,a+\delta} : ||r||_{[a,D_a(r)]}<\rho \},\\
V_1&=&\{ r \in  O_{a',a'+\delta} : ||r||_{[G_{a'}(r),D_{a'}(r)]} < 
||r||_{[D_a(r),G_{a'}(r)]} \in ]h-\Delta h,h[\},\\ 
V_2&=&\{ r \in O_{b'-\delta,b'} : 
||r-g||_{[D_{a'}(r),D_{b'-\delta}(r)]} < \rho' \},\\ 
V_3&=&\{ r \in O_{b-\delta,b} : 
||r||_{[D_{b'-\delta}(r),D_{b-\delta}(r)]} < \rho'/2 \},\\ 
V_4&=&\{ r \in O_{b-\delta,b} : h < ||r||_{[D_{b-\delta}(r),b]} < \rho/2 \}. 
\end{eqnarray*}
Alors $V_0$, $V_1$, $V_2$, $V_3$, $V_4$ sont des presque-ouverts 
(gr{\^a}ce au lemme~\ref{Continuite des parametres d'une excursion}) et :
\begin{itemize}
\item $V_1$ est accessible par $F_a$ depuis presque tout $r \in V_0$, 
gr{\^a}ce au lemme~\ref{Obtention d'une plus grande excursion de hauteur voulue} 
(obtention d'une plus grande excursion de hauteur voulue) ;  
\item $V_2$ est accessible par $F_{a'}$ depuis presque tout $r \in V_1$, 
d'apr{\`e}s l'hypoth{\`e}se du lemme et l'inclusion 
$U_{a',b'}(g,\rho',\delta) \subset V_2$ ;  
\item $V_3$ est accessible par $F_{b'-\delta}$ depuis presque tout $r \in V_2$, 
gr{\^a}ce au lemme~\ref{Remise a zero apres l'instant $D_a$} (on remet la 
trajectoire {\`a} z{\'e}ro apr{\`e}s l'instant $D_{b'-\delta}$).  
\item $V_4$ est accessible par $F_{b-\delta}$ depuis presque tout $r \in V_3$, 
gr{\^a}ce au lemme~\ref{Obtention d'une plus grande excursion de hauteur voulue} 
(obtention d'une plus grande excursion de hauteur voulue) ;  
 \end{itemize}
Comme $V_0$, $V_1$, $V_2$, $V_3$, sont stables par $F_a$, $F_{a'}$, 
$F_{b'-\delta}$, $F_{b-\delta}$ respectivement, le presque-ouvert
$V_0 \cap V_1 \cap V_2 \cap V_3 \cap V_4$ est accessible par $F_a$ 
depuis presque tout $r \in V_0$. 
Il suffit donc de montrer que $U_{a,b}(f,\rho,\delta)$ est accessible par 
$F_a$ depuis toute trajectoire de $V_0 \cap V_1 \cap V_2 \cap V_3 \cap V_4$.

Soient donc $r \in V_0 \cap V_1 \cap V_2 \cap V_3 \cap V_4$,
$e \in E$ et $\widetilde{r} = F_a(e,r)$. Nous allons montrer que 
$\widetilde{r} \in U_{a,b}(f,\rho,\delta)$ d{\`e}s que l'action de $e$ sur $r$ 
affecte 
\begin{itemize}
\item du signe $-$ la plus grande excursion sur $[D_a(r),G_{a'}(r)]$,  
\item du signe $-$ l'excursion r{\'e}alisant le maximum de $r$ sur 
$[D_{b-\delta}(r),b]$, 
\item du signe $+$ les autres excursions de hauteur 
$\ge \min(\rho - \rho' - \Delta h,\rho')/2$ avant $D_b(r)$.   
\end{itemize}
Pour cela, notons $h_1$ la hauteur de la plus grande excursion 
de $r$ sur $[D_a(r),G_{a'}(r)]$ et $h_2$ le maximum de $r$ sur 
$[D_{b-\delta}(r),b]$. Alors $h - \Delta h < h_1 < h < h_2$. 

On a bien s{\^u}r $\widetilde{r} \in O_{a,a+\delta}$ puisque 
$D_a(\widetilde{r})=D_a(r)$. Par ailleurs, 
$\widetilde{r} \in O_{b-\delta,b}$ puisque d'apr{\`e}s le choix des signes, 
l'instant r{\'e}alisant le maximum de $r$ sur 
$[D_{b-\delta}(r),b]$ est un instant de record n{\'e}gatif de $e \cdot r$ 
{\`a} partir de l'instant $D_a(r)$, et donc un z{\'e}ro de $\widetilde{r}$.  

Il reste {\`a} montrer que pour tout $t \in [a,b]$, 
$|\widetilde{r}(t)-f(t)|<\rho$. On distingue quatre cas.
\begin{enumerate}
\item Si $t \in [a,D_a(r)]$, alors $\widetilde{r}(t) = r(t) < \rho$ 
puisque $r \in V_0$, donc
$$|\widetilde{r}(t)-f(t)| \le \max(\widetilde{r}(t),f(t)) \le 
\max(r(t),h)<\rho.$$
\item Si $t \in [D_a(r),D_{a'}(r)]$, alors 
$$0 \le \widetilde{r}(t) \le 2||r||_{[D_a(r),D_{a'}(r)]} < 2h <\rho \text{ car } 
r \in V_1,$$
$$0 \le f(t) \le h+\osc(f,\delta) <\rho/2 + \rho'/4 <\rho,$$
donc
$$|\widetilde{r}(t)-f(t)| \le \max(\widetilde{r}(t),f(t)) <\rho.$$
\item Si $t \in [D_{a'}(r),b']$, alors $|r(t)-g(t)| < \rho'$. En effet, 
pour $t \in [D_{a'}(r),D_{b'-\delta}(r)]$ cela vient du fait que $r \in V_2$ ; 
pour $t \in [D_{b'-\delta}(r),b']$ cela vient du fait que $r \in V_3$ et 
$\osc(g,\delta)<\rho'/4$ d'o{\`u} $0 < r(t) < \rho'/2$ et $0 < g(t) < \rho'/4$.  
D'apr{\`e}s le choix des signes, on a donc
\begin{eqnarray*}
|\widetilde{r}(t)-f(t)| 
&=& |(e \cdot r)(t) + h_1 - g(t) - h|\\
&\le& |(e \cdot r)(t) - r(t)| + |r(t) - g(t)| + |h_1 - h|\\
&<& (\rho-\rho'-\Delta h) + \rho' + \Delta h\\
&=& \rho.
\end{eqnarray*} 
\item Si $t \in [b',b]$, alors $r(t) < \rho/2$. En effet,  
pour $t \in [b',D_{b-\delta}(r)]$ cela vient du fait que $r \in V_3$ et de 
l'in{\'e}galit{\'e} $\rho'/2<\rho/2$ ; pour $t \in [D_{b-\delta}(r),b]$ cela vient 
du fait que $r \in V_4$. Comme $e \cdot r$ est minor{\'e} par $-\rho/2$ sur 
$[a,D_b]$, on a donc
$$|\widetilde{r}(t)-f(t)| \le \max(\widetilde{r}(t),f(t)) \le 
\max(r(t)+\rho/2,h)<\rho.$$
\end{enumerate}  
Dans tous les cas, $|\widetilde{r}(t)-f(t)|<\rho$, ce qui ach{\`e}ve 
la preuve dans le cas o{\`u} $h<\rho/2$.  

\vfill\eject
\paragraph{2.} D{\'e}monstration dans le cas o{\`u} $h \ge \rho/2$. 

Choisissons $h' \in [\rho/4,\rho/2[$ et $n \ge 2$ entier tels que $h=nh'$. 
Soient $\rho_0=\rho'$ et $\rho_1,\ldots,\rho_n$ tels que
$$\max(\rho_0,2h') < \rho_1 < \ldots < \rho_n = \rho.$$
Soient 
$$a = a_n < \ldots < a_0 = a'.$$
$$b' = b_0 < \ldots < b_n = b.$$
les deux subdivisions r{\'e}guli{\`e}res de $[a,a']$ et $[b',b]$ en $n$ 
sous-intervalles. Pour $k \in [0 \ldots n]$, notons $f_k$ la 
trajectoire d{\'e}finie par
$$f_k(t) = [f(t) - (n-k)h' ]_+$$
et $\hc_k$ l'affirmation : \og l'ouvert $U_{a_k,b_k}(f_k,\rho_k,\delta)$ 
est accessible par $F_{a_k}$ depuis presque tout $r \in O_{a_k,a_k+\delta}$ 
tel que $||r||_{[a_k,D_{a_k}(r)]}<\rho_k$\fg. 
Par construction, $f_k$ est nulle en dehors de $[a_k,b_k]$.   
Comme $f_0=g$ et $f_n = f$, il s'agit de d{\'e}montrer que 
$\hc_0 \Rightarrow \hc_n$. 

Il suffit de montrer que pour tout $k \in [0 \ldots n-1]$, 
$\hc_k \Rightarrow \hc_{k+1}$. Pour cela, on remarque que  
$f_{k+1}$ se d{\'e}duit de $f_k$ par levage de la hauteur $h'<\eta/2$ sur 
l'intervalle $[a_k,b_k]$ et interpolation lin{\'e}aire sur les intervalles 
$[a_{k+1},a_k]$ et $[b_k,b_{k+1}]$. Comme 
$$\osc(f_k,\delta) \le \osc(f,\delta) < \rho'/4 < \rho_k/4,$$
$b_k-a_k \ge b'-a' >\delta$, $a_k-a_{k+1}>\delta$ et $b_{k+1} - b_k>\delta$ 
compte tenu des in{\'e}galit{\'e}s  
$$f(a_k)-f(a_{k+1}) = f(b_{k+1}) - f(b_k) = h' \ge \rho/4 > \osc(f,\delta),$$ 
il suffit d'appliquer le lemme du v{\'e}rin dans le cas o{\`u} il est d{\'e}j{\`a} d{\'e}montr{\'e}. 
\hfill $\square$
\end{demo}

\section{Approximation des fonctions continues affines par 
morceaux}

Dans cette partie, nous allons voir comment approcher les trajectoires 
affines par morceaux pour la topologie CUCZ. La d{\'e}monstration, r{\'e}sum{\'e}e 
par la figure~\ref{etapes de l'approximation} repose sur les propositions
~\ref{Accessibilites successives et intersection} 
(Accessibilit{\'e} successives et intersection)
et ~\ref{Remise a zero apres l'instant $D_a$} 
(Remise {\`a} z{\'e}ro apr{\`e}s l'instant $D_a$) 
et sur le lemme~\ref{Lemme du verin} (Lemme du v{\'e}rin). 

Le lemme suivant 
montre que pour approcher une trajectoire $f$ entre deux de ses z{\'e}ros, 
il suffit de savoir approcher les morceaux obtenus en d{\'e}coupant 
$f$ {\`a} un nombre finis de z{\'e}ros interm{\'e}diaires. 

Pour $b > a \ge 0$, $f \in \wwf_+$ telle que $f(a)=f(b)=0$, $\rho>0$ et 
$0 < \delta < b-a$, on note toujours
$$U_{a,b}(f,\rho,\delta) = \{r \in O_{a,a+\delta} \cap O_{b-\delta,b} : 
||r-f||_{[a,b]} < \rho\}.$$
Si $a=z_0 < \ldots < z_n=b$ et 
$\delta < \min(z_1-z_0,\ldots,z_n-z_{n-1})$, on note {\'e}galement 
$$U_{z_0,\ldots,z_n}(f,\rho,\delta) 
= \bigcap_{k=0}^{n-1} U_{z_k,z_{k+1}}(f,\rho,\delta),$$
On remarque que cet ouvert est contenu dans $U_{a,b}(f,\rho,\delta)$.

\begin{lemm}\label{Concatenation de ponts} 
{\bf (Concat{\'e}nation d'un nombre fini de ponts)}
  
Soit $f \in \wwf_+$. Soient $a = z_0 < \ldots < z_n = b$ des z{\'e}ros de $f$. 
Soient $\rho>0$ et $\delta>0$ tels que $\osc(f|_{[a,b]},\delta)<\rho$ et 
$\delta < \min(z_1-z_0,\ldots,z_n-z_{n-1})$. 

Supposons que pour tout $k \in [0 \ldots n-1]$, l'ouvert 
$U_{z_k,z_{k+1}}(f,\rho,\delta)$ 
est accessible par $F_{z_k}$ depuis presque toute trajectoire 
$r \in O_{z_k,z_k+\delta}$ telle que $||r||_{[z_k,D_{z_k}(r)]}<\rho$. 

Alors l'ouvert $U_{z_0,\ldots,z_n}(f,\rho,\delta)$ est accessible par $F_a$ 
depuis presque tout $r \in O_{a,a+\delta}$ tel que $||r||_{[a,D_a(r)]}<\rho$
(et donc depuis presque tout $r \in \wwf_+$ lorsque $a=0$). 
\end{lemm}

\begin{demo}
Le lemme se d{\'e}montre par r{\'e}currence sur le nombre $n$ de ponts. 

Pour $n=1$ il n'y a rien {\`a} montrer. 

Soit $n \ge 2$. Supposons la propri{\'e}t{\'e} {\'e}tablie pour $n-1$ ponts. 
Soient $f \in \wwf_+$, $a = z_0 < \ldots < z_n = b$, $\rho>0$ et $\delta>0$ 
comme dans l'{\'e}nonc{\'e}. 
On remarque que 
$$U_{z_0,\ldots,z_n}(f,\rho,\delta) = V_1 \cap V_2 \cap V_3$$
o{\`u} 
$$V_1 = U_{z_0,\ldots,z_{n-2}}(f,\rho,\delta) \cap 
\{r \in O_{[z_{n-2},z_{n-2}+\delta]} \cap O_{[z_{n-1}-\delta,z_{n-1}]} : 
||r||_{[z_{n-2},D_{z_{n-1}-\delta}(r)]} < \rho\},$$
$$V_2 = \{r \in O_{z_{n-1},z_{n-1}+\delta} : 
||r||_{[D_{z_{n-1}-\delta}(r),D_{z_{n-1}}(r)]}<\rho\}$$  
$$V_3 = U_{z_{n-1},z_n}(f,\rho,\delta)$$

Le presque-ouvert $V_1$ contient $U_{z_0,\ldots,z_{n-1}}(f,\rho,\delta)$, 
donc par hypoth{\`e}se de r{\'e}currence, $V_1$ est accessible par $F_a$ depuis 
presque tout $r \in O_{a,a+\delta}$ tel que $||r||_{[a,D_a(r)]}<\rho$. 
De plus, $V_1$ est stable par $F_{z_{n-1}-\delta}$. 

Par remise {\`a} z{\'e}ro apr{\`e}s l'instant $D_{z_{n-1}-\delta}$, l'ouvert 
$V_2$ est accessible par $F_{z_{n-1}-\delta}$ depuis presque toute 
trajectoire de $O_{[z_{n-1}-\delta,z_{n-1}]}$. De plus, $V_2$ est 
stable par $F_{z_{n-1}}$.

Enfin, par hypoth{\`e}se, l'ouvert $V_3$ est 
accessible par $F_{z_{n-1}}$ depuis presque toute trajectoire 
$r \in O_{z_{n-1},z_{n-1}+\delta}$ telle que 
$||r||_{[z_{n-1},D_{z_{n-1}}(r)]}<\rho$ et donc de presque toute 
trajectoire de $V_1 \cap V_2$. 

Ainsi, $U_{z_0,\ldots,z_n}(f,\rho,\delta) = V_1 \cap V_2 \cap V_3$ est 
accessible depuis presque tout $r \in O_{a,a+\delta}$ tel que 
$||r||_{[a,D_a(r)]}<\rho$. \hfill $\square$
\end{demo}

{\`A} l'aide du lemme pr{\'e}c{\'e}dent, du lemme de remise {\`a} z{\'e}ro et du lemme du v{\'e}rin, 
nous allons d{\'e}montrer l'accessibilit{\'e} des ouverts de la forme 
$U_{a,b}(f,\rho,\delta)$. 

\begin{prop} {\bf (Approximation d'un pont affine par morceaux)}

Soient $b > a \ge 0$, $f \in \wwf_+$ une trajectoire affine par morceaux 
sur $[a,b]$ telle que $f(a)=f(b)=0$. Soient $\rho>0$ et $\delta>0$.  
Alors l'ouvert $U_{a,b}(f,\rho,\delta)$ est accessible par $F_a$ depuis 
presque tout $r \in O_{a,a+\delta}$ tel que $||r||_{[a,D_a(r)]}<\rho$.
\end{prop}

\begin{demo}
Par hypoth{\`e}se, on peut trouver une subdivision $a=c_0<\ldots<c_n=b$ tels que 
$f \in \wwf_+$ soit affine sur chaque segment $[c_{k-1},c_k]$. 
De plus, l'image r{\'e}ciproque de tout r{\'e}el par la restriction de $f$ 
{\`a} $[a,b]$ est une union finie de singletons et d'intervalles de subdivision. 
Quitte {\`a} raffiner la subdivision, on peut donc supposer que pour tout point 
$c$ de subdivision, l'image r{\'e}ciproque de $f(c)$ par $f$ restreinte {\`a} $[a,b]$ 
est form{\'e}e uniquement de points de subdivision et d'intervalles de 
subdivision. Nous dirons alors que la subdivision de $[a,b]$ est compl{\`e}te 
relativement {\`a} $f$. 

Quitte {\`a} r{\'e}duire $\delta$, on peut supposer de plus que 
$\delta<\min(c_1-c_0,\ldots,c_n-c_{n-1})$ et que 
$\osc(f|_{[a,b]},\delta) < \rho/4$.  

Sous ces restrictions, on effectue alors une r{\'e}currence sur le nombre $n$ 
d'intervalles pour une subdivision compl{\`e}te. 

Remarquons d'abord que si $||f||_{[a,b]}<\rho$, alors 
$$\{ r \in \wwf_+ : ||r||_{[a,b]}<\rho \} \subset 
\{ r \in \wwf_+ : ||r-f||_{[a,b]}<\rho \}$$
puisque $||r-f||_{[a,b]} \le \max(||r||_{[a,b]},||f||_{[a,b]})$ pour tout 
$r \in \wwf_+$. Dans ce cas, il suffit d'appliquer le th{\'e}or{\`e}me de remise 
{\`a} z{\'e}ro apr{\`e}s l'instant $D_a$. 

\paragraph{1.} Le cas o{\`u} $n=1$ rentre dans ce cas particulier puisque $f$ est 
alors la fonction nulle. 

\paragraph{2.} Le cas o{\`u} $n=2$ est une application directe du lemme du v{\'e}rin.

En effet, supposons que $||f||_{[a,b]} \ge \rho$ (sans quoi il n'y a 
rien {\`a} montrer). Choisissons un r{\'e}el $\rho'$ tel que 
$4\osc(f|_{[a,b]},\delta) < \rho' < \rho$ et posons $h=f(c_1)-\rho'/2$. 
Alors $0 < \rho'/2 < \rho-\rho'/2 \le h<f(c_1)$. 
Notons $a'<b'$ les ant{\'e}c{\'e}dents de $h$ par $f$ et 
$g=(f-h)_+$. 

Comme $g$ est major{\'e}e par $f(c_1) - h = \rho'/2$, l'ouvert 
$U_{a',b'}(g,\rho',\delta)$ est accessible par $F_{a'}$ depuis 
presque tout $r \in O_{a',a'+\delta}$ tel que $||r||_{[a',D_{a'}(r)]}<\rho$, 
d'apr{\`e}s la remarque pr{\'e}liminaire. 

L'application $f$ s'obtient {\`a} partir de $g$ par levage 
de la hauteur $h$ et par interpolation lin{\'e}aire sur $[a,a']$ et $[b',b]$.
Par ailleurs $\delta$ v{\'e}rifie  
$\osc(g|_{[a',b']},\delta) \le \osc(f|_{[a,b]},\delta) < \rho'/4$ et
$\delta < \min(a'-a,c_1-a',b'-c_1,b-b')$ puisque 
$f(a')-f(a) = f(b')-f(b) = h > \rho'/4$ et 
$f(c_1)-f(a') = f(c_1)-f(b') = \rho'/2 > \rho'/4$.

Les hypoth{\`e}ses du lemme du v{\'e}rin sont satisfaites, ce qui montre que 
$U_{a,b}(f,\rho,\delta)$ est accessible par $F_a$ depuis presque tout 
$r \in O_{a,a+\delta}$ tel que $||r||_{[a,D_a(r)]}<\rho$.

\paragraph{3.} Montrons la propri{\'e}t{\'e} pour $n \ge 3$ en la supposant {\'e}tablie 
pour un nombre d'intervalles au plus {\'e}gal {\`a} $n-1$. Soit 
$h = \min\{f(c_k)\ ; 1 \le k \le n-1\}$. On distingue deux cas. 

Soit $h=0$. Dans ce cas, il existe $m \in [1 \ldots n-1]$ tel que $f(c_m)=0$. 
Les subdivisions $a=c_0<\ldots<c_m=z$ et $z=c_m<\ldots<c_n=b$ 
sont compl{\`e}tes relativement {\`a} $f$. Par hypoth{\`e}se de r{\'e}currence, l'ouvert 
$U_{a,z}(f,\rho,\delta)$ est accessible par $F_a$ depuis 
presque tout $r \in O_{a,a+\delta}$ tel que $||r||_{[a,D_a(r)]}<\rho$ et 
l'ouvert $U_{z,b}(f,\rho,\delta)$ est accessible par $F_z$ depuis 
presque tout $r \in O_{z,z+\delta}$ tel que $||r||_{[z,D_z(r)]}<\rho$. 
D'apr{\`e}s le lemme de concat{\'e}nation de ponts,
$U_{a,b}(f,\rho,\delta)$ est accessible par $F_a$ 
depuis presque tout $r \in O_{a,a+\delta}$ tel que $||r||_{[a,D_a(r)]}<\rho$. 

Soit $h>0$. Comme la subdivision $a=c_0<\ldots<c_n=b$ de 
$[a,b]$ est compl{\`e}te relativement {\`a} $f$, on a n{\'e}cessairement 
$f(c_1) = f(c_{n-1}) = h$. On pose $a'=c_1$, $b'=c_{n-1}$ et $g=(f-h)_+$. 
La trajectoire $g$ est affine par morceaux sur $[a',b']$ et 
v{\'e}rifie $g(a')=g(b')=0$. De plus, la subdivision 
$a' = c_1<\ldots<c_{n-1} = b'$ est compl{\`e}te relativement {\`a} $g$ 
et v{\'e}rifie $\delta<\min(c_2-c_1,\ldots,c_{n-1}-c_{n-2})$. Si l'on fixe 
$\rho' \in ]4\osc(f,\delta),\rho[$, on a alors
$$\osc(g|_{[a',b']},\delta) \le \osc(f|_{[a,b]},\delta) < \rho'/4.$$ 
Par cons{\'e}quent, on peut donc appliquer l'hypoth{\`e}se de r{\'e}currence {\`a} 
$g$, $\rho'$, $\delta$ et {\`a} la subdivision 
$a' = c_1 < \ldots <c_{n-1} = b'$, ce qui montre que l'ouvert 
$U_{a',b'}(g,\rho',\delta)$ est accessible par $F_{c_1}$ 
depuis presque tout $r \in O_{c_1,c_1+\delta}$ tel que 
$||r||_{[c_1,D_{c_1}(r)]}<\rho$. Mais $f$ s'obtient {\`a} partir de $g$ par levage 
de la hauteur $h$ et par interpolation lin{\'e}aire sur $[a,a']$ et $[b',b]$. 
Comme $\delta < \min(a'-a,b'-a',b-b')$, on peut appliquer le 
lemme du v{\'e}rin. Ainsi, l'ouvert $U_{a,b}(f,\rho,\delta)$ est accessible 
par $F_a$ 
depuis presque tout $r \in O_{a,a+\delta}$ tel que $||r||_{[a,D_a(r)]}<\rho$.
\hfill $\square$
\end{demo}

\begin{prop} {\bf (Approximation CUCZ d'une fonction affine par morceaux)}

Soit $f \in \wwf_+$, affine par morceaux. Pour tout $t>0$, $\rho>0$ 
et $\delta>0$, l'ouvert $V_t(f,\rho,\delta)$ est accessible par $F$ 
depuis presque toute trajectoire de $\wwf_+$.
\end{prop}

\begin{demo}
Quitte {\`a} modifier $f$ apr{\`e}s l'instant $t$ et {\`a} augmenter $t$, 
on peut se limiter au cas o{\`u} $f(t)=0$. Quitte {\`a} r{\'e}duire $\delta$, 
on peut supposer que $\delta$ est strictement inférieur à la longueur des 
excursions de $f$ sur $[0,t]$, puisque ces excursions sont en nombre fini. 
On choisit alors un nombre fini de z{\'e}ros 
de $f$, $0 = z_0 < \ldots < z_n = t$ de telle sorte que
$\delta < \min(z_1-z_0,\ldots,z_n-z_{n-1})$ et que tout z{\'e}ro 
de $f$ avant $t$ soit {\`a} distance $<\delta$ d'un de ces z{\'e}ros. 

Le r{\'e}sultat d'approximation des ponts appliqu{\'e} {\`a} $f$ sur chaque 
intervalle $[z_k,z_{k+1}]$ montre que l'ouvert 
$U_{z_k,z_{k+1}}(f,\rho,\delta)$ 
est accessible par $F_{z_k}$ depuis presque toute trajectoire 
$r \in O_{z_k,z_k+\delta}$ telle que $||r||_{[z_k,D_{z_k}(r)]}<\rho$. 

Par concat{\'e}nation d'un nombre fini de ponts, l'ouvert 
$$U_{z_0,\ldots,z_n}(f,\rho,\delta) 
= \bigcap_{k=0}^{n-1} U_{z_k,z_{k+1}}(f,\rho,\delta),$$
est accessible par $F_0$ depuis presque tout 
$r \in \wwf_+$. On termine en remarquant que 
$U_{z_0,\ldots,z_n}(f,\rho,\delta)$ 
est inclus dans $V_t(f,\rho,\delta)$.
\hfill $\square$
\end{demo}

Cette proposition ach{\`e}ve la d{\'e}monstration du 
th{\'e}or{\`e}me~\ref{Densite pour la topologie CUCZ}
compte tenu de la proposition~\ref{Ce qu'il faudra montrer}
(Passage de la cha{\^\i}ne $(R_n)_{n \in \nnf}$ {\`a} la cha{\^\i}ne $(B_n)_{n \in \nnf}$)
et de la densit{\'e} des fonctions continues, positives, nulles en $0$ 
et affines par morceaux dans $\wwf_+$ 
(lemme~\ref{Densite des fonctions continues affines par morceaux}). 

\begin{figure}[hbtp]\label{etapes de l'approximation}
  \begin{center}
    \leavevmode \input{approximation.pstex_t} 
  \end{center}
  Figure~\ref{etapes de l'approximation}. --- {\'E}tapes de l'approximation 
d'une fonction continue affine par morceaux. 
\end{figure}

\vfill
\eject

\section{Comparaison avec les outils de Malric}\label{comparaison}

\subsection{Parties atteignables}

Nous avons utilisé dans cet article la notion d'accessibilité. Malric 
utilise dans~\cite{Malric 5} une notion voisine, mais moins souple, 
dont nous adaptons la définition à toute transformation mesurable $T$ 
d'un espace probabilis{\'e} $(E,\ec,\pi)$ dans lui-m{\^e}me pr{\'e}servant 
la mesure. 

\begin{defi} {\bf (Atteignabilité)}
Une partie mesurable $G \in \ec$ est dite atteignable si pour 
tout $\eta>0$, il existe $n \in \nnf$ et une probabilité 
$\nu$ sur $(E,\ec)$ tels que : 
\begin{enumerate}
\item $\nu$ est absolument continue par rapport à $\pi$.
\item $T^n(\nu) = \pi$. 
\item $\nu(G) > 1-\eta$. 
\end{enumerate}  
\end{defi}

La proposition 1 de~\cite{Malric 5}, énoncée dans le cas de la 
transformation de Lévy, s'étend sans difficulté  au cas général : 
si $G$ est atteignable, alors pour $\pi$-presque tout $x \in E$, l'orbite 
$\{T^n(x)\ ;\ n \in \nnf\}$ visite $G$. 

Dans le cas où l'espace probabilisé $(E,\ec)$ est un espace polonais muni 
de la tribu borélienne, ce résultat peut être vu comme une conséquence de 
notre corollaire~\ref{Condition suffisante de recurrence}. En effet, 
l'atteignabilité de $G$ équivaut à une condition forte d'accessibilité que 
nous allons définir. 

\begin{defi} {\bf (Accessibilité, accessibilité forte)}

Soient $X$ une variable aléatoire de loi $\pi$, définie sur un certain 
espace probabilisé $(\Omega,\ac,P)$ et $K(\cdot,\cdot)$ une version 
régulière de la loi conditionnelle de $X$ sachant $T(X)$.
\begin{itemize}
\item On dit que $G$ est accessible depuis $x \in E$ s'il existe $n \in \nnf$
tel que $K^n(x,G)>0$. 
\item On dit que $G$ est fortement accessible si 
$\sup_{n \in \nnf} \pi\{x \in G : K^n(x,G)>0\} = 1$.
\end{itemize}
\end{defi}

L'accessibilité de $G$ depuis $\pi$-presque tout $x \in E$ constitue 
l'hypothèse de 
notre corollaire~\ref{Condition suffisante de recurrence}. 
Pour déduire de ce corollaire la proposition 1 de~\cite{Malric 5},
il suffit de démontrer les implications contenues dans la 
proposition ci-dessous.

\begin{prop} {\bf (Lien entre atteignabilité et accessibilité)}
\begin{enumerate}
\item $G$ est atteignable si et seulement si $G$ 
est fortement accessible.
\item Si $G$ est fortement accessible, alors $G$ est accessible 
depuis $\pi$-presque tout $x \in E$. 
\end{enumerate}
\end{prop}

\begin{demo}
Pour $n \in \nnf$, notons $A_n = \{x \in E : K^n(x,G) > 0\}$. 
Le point 2 découle immédiatement des inégalités 
$$\sup_{n \in \nnf} \pi(A_n) \le \pi(\bigcup_{n \in \nnf} A_n) \le 1.$$ 
Montrons le point 1. 

Supposons que $G$ est fortement accessible et montrons que 
$G$ est atteignable. Fixons $\eta>0$. Par hypothèse, on peut trouver 
$n \in \nnf$ tel que $\pi(A_n)>1-\eta$. Définissons un noyau de 
transition $L$ et une mesure $\nu$ sur $(E,\ec)$ par 
$$\begin{array}{l}
L(x,B) = K^n(x,B \cap G)/K^n(x,G) \text { si } x \in A_n,\\
L(x,B) = K^n(x,B) \text { sinon, }
\end{array}$$
Alors pour $\pi$-presque tout $x \in E$, la probabilité $L(x,\cdot)$ 
est absolument continue par rapport à $K^n(x,\cdot)$, qui est portée 
par $(T^n)^{-1}(\{x\})$. De plus $L(x,G)=1$ si $x \in A_n$, 
$L(x,G)=0$ sinon. 

Soit $\nu$ la mesure sur $(E,\ec)$ définie par 
$$\nu(B) = \int_E L(x,B)\ d\pi(x).$$
Alors 
\begin{enumerate}
\item Si $\nu(B)=0$, alors pour $\pi$-presque tout $x \in E$, $L(x,B)=0$ 
d'où $K^n(x,B)=0$, et donc 
$$\pi(B) = \int_E K^n(x,B)\ d\pi(x) = 0.$$
Donc $\nu$ est absolument continue par rapport à $\pi$.
\item Comme pour $\pi$-presque tout $x \in E$, la mesure image de 
$L(x,\cdot)$ par $T^n$ est $\delta_x$, on a $T^n(\nu) = \pi$. 
\item Enfin, $\nu(G) = \pi(A_n)> 1-\eta$. 
\end{enumerate}  
Cela montre que $G$ est atteignable.

Réciproquement, supposons que $G$ est atteignable. Fixons $\eta>0$. 
Choisissons $n \in \nnf$ et $\nu$ vérifiant les points 1, 2 et 3 
de la définition. Soient $Y$ une variable aléatoire de loi $\nu$ et 
$L_n(\cdot,\cdot)$ une version régulière de la loi conditionnelle de 
$Y$ sachant $T^n(Y)$. Alors 
$$\pi[G \cap (T^n)^{-1}(A_n^c)] = P[X \in G\ ;\ T^n(X) \in A_n^c] 
= \int_{A_n^c} K^n(x,G)\ d\pi(x) = 0,$$
donc par absolue continuité,
$$0 = \nu[G \cap (T^n)^{-1}(A_n^c)] = P[Y \in G\ ;\ T^n(Y) \in A_n^c] 
= \int_{A_n^c} L_n(x,G)\ d\pi(x).$$
Donc $L_n(x,G) = 0$ pour $\pi$-presque tout $x \in A_n^c$. 
Par conséquent, pour $\pi$-presque tout $x \in E$ 
$$\iif_{A_n}(x) \ge \iif_{[L_n(x,G) > 0]} \ge L_n(x,G),$$
d'où en intégrant par rapport à $\pi$,
\begin{eqnarray*}
\pi(A_n) \ge \int_E L_n(x,G)\ d\pi(x) = P[Y \in G] > 1-\eta.
\end{eqnarray*}
Comme $\eta>0$ est arbitraire, on en déduit que $G$ est fortement accessible.
\hfill $\square$
\end{demo}

Remarque : avec les notations de la démonstration précédente, 
on a pour tout $n \in \nnf$, $A_n \subset T^{-1}(A_{n+1})$ 
$\pi$-presque sûrement, d'où 
$\pi(A_n) \le \pi(T^{-1}(A_{n+1})) = \pi(A_{n+1})$. La suite 
$(\pi(A_n))_{n \in \nnf}$ est donc croissante. Pour montrer 
cette inclusion presque sûre, on écrit
$$\pi(A_n \setminus T^{-1}(A_{n+1})) = P[X \in A_n\ ;\ T(X) \in A_{n+1}^c] 
= \int_{A_{n+1}^c} K(x,A_n)\ d\pi(x).$$
Mais pour tout $x \in E$, si $K(x,A_n) > 0$, alors 
$$K^{n+1}(x,G) \ge \int_{A_n} K(x,dy) K^n(y,G) > 0$$ 
puisque $K^n(y,G)>0$ pour tout $y \in A_n$, donc $x \in A_{n+1}$.  
Donc $\pi(A_n \setminus T^{-1}(A_{n+1})) = 0$.

\subsection{Remontées de Lévy et remontées partielles}

Alors que nous utilisons des chaînes de Markov définies à l'aide de 
la transformation inverse de Lévy, la définition de l'atteignabilité 
conduit Malric à introduire les remontées de Lévy. 

Lorsque $U$ et $V$ sont deux variables aléatoires à valeurs dans $\wwf$, 
de loi absolument continue par rapport à la mesure de Wiener, telles que 
$\ttf(V) = U$, Malric dit que $V$ est un remonté de Lévy de $U$. Une façon 
simple de construire un remonté de Lévy de $U$ est la suivante 
(proposition 3 de \cite{Malric 5}) : retrancher 
à $U$ son minimum courant, multiplier les excursions hors de 0 
du processus obtenu par les signes fournis par un jeu de pile ou face 
indépendant dont on a modifié les signes pour ensemble fini (mais aléatoire) 
d'indices.  

Étant donné un ouvert $G$ de $\wwf$ et un mouvement brownien $W$, Malric 
cherche à construire des remontés de Lévy successifs de sorte que, 
pour un entier $n$ bien choisi, le 
$n$-ième remonté tombe dans $G$ avec probabilité proche de $1$. Cette méthode 
semble plus rigide que l'utilisation de notre 
corollaire~\ref{Condition suffisante de recurrence} puisque le nombre 
$n$ de remontées successives est déterministe. 

Malric contourne cette difficulté par l'introduction de la notion de 
remonté partiel de Lévy : si $U$ et $V$ sont deux variables aléatoires 
à valeurs dans $\wwf$, 
de loi absolument continue par rapport à la mesure de Wiener, $V$ 
est un remonté de partiel de Lévy de $U$ s'il existe un 
temps aléatoire $T$ à valeurs dans $[0,+\infty]$ tels que 
\begin{enumerate}
\item $T$ est une fin d'excursion de $U$ sur l'événement $[0<T<+\infty]$.
\item $U$ et $V$ coïncident sur l'intervalle de temps $[0,T]$.
\item Pour $t \ge T$, $|V_t| = U_t - \min\{U_s\ ; s \in [T,t]\}$. 
\end{enumerate}  

Un résultat-clé de Malric (proposition 5 de \cite{Malric 5}) affirme 
que toute suite finie de remontés partiels successifs de Lévy peut 
être approchée uniformément par une suite de remontés de Lévy.
Ce résultat permet de s'affranchir de la contrainte du 
nombre déterministe de remontées dans la définition de l'atteignabilité. 
Il permet aussi de préserver le début des trajectoires. Il joue finalement 
le même rôle que notre proposition~\ref{Comparaison des accessibilites} 
(Comparaison des accessibilités par $F_a$ pour différentes valeurs de $a$), 
dans laquelle les transformations $F_a$ jouent le rôle des remontées 
partielles de Lévy. 

Les autres différences marquantes, déjà notées dans l'introduction, 
sont que nous travaillons avec des transformations sur les mouvements 
browniens réfléchis en utilisant des propriétés de continuité dans la 
topologie de la convergence uniforme sur les compacts avec contrôle des zéros. 

\vfill
\eject

\vfill
\begin{flushleft}
Jean BROSSARD et Christophe LEURIDAN\\
INSTITUT FOURIER\\
Laboratoire de Math\'ematiques\\
UMR5582 (UJF-CNRS)\\
BP 74\\
38402 St MARTIN D'H\`ERES Cedex (France)\\
Christophe.Leuridan@ujf-grenoble.fr
\end{flushleft}

\end{document}